\documentclass[leqno,11pt]{amsart}
\usepackage{amsmath,amscd,amsthm}
\usepackage{epsfig,ifsym,color}
\usepackage{amssymb,latexsym}
\usepackage{mathrsfs}
\usepackage{hyperref}
\usepackage[all,cmtip]{xy}
\usepackage{graphicx}
\usepackage{caption}
\usepackage{subcaption}

\newtheorem{thm}{\bf Theorem}[section]
\newtheorem{df}[thm]{\bf Definition}
\newtheorem{prop}[thm]{\bf Proposition}
\newtheorem{cor}[thm]{\bf Corollary}
\newtheorem{lem}[thm]{\bf Lemma}
\newtheorem{rem}[thm]{\bf Remark}
\newtheorem{ex}[thm]{\bf Example}

\newtheorem{nono-theorem}{Theorem}[]
\newtheorem*{thm*}{Theorem}
\newtheorem{alg}{Algorithm}

\newcommand{\A}{\mathscr{A}}

\newcommand{\B}{\mathbf{B}}

\newcommand{\cP}{\mathscr{P}}

\newcommand{\pf}{\noindent{\bfseries Proof. }}
\newcommand{\ov}{\overline}

\newcommand{\hf}{\frac{1}{2}}

\newcommand{\gl}{\mathfrak{gl}}
\newcommand{\Z}{\mathbb{Z}}
\newcommand{\C}{\mathbb{C}}
\newcommand{\h}{\mathfrak{h}}
\newcommand{\te}{\widetilde{e}}
\newcommand{\tf}{\widetilde{f}}
\newcommand{\g}{\mathfrak{g}}
\newcommand{\td}{\widetilde}
\newcommand{\tE}{\mc{E}}
\newcommand{\tF}{\mc{F}}

\newcommand{\mc}{\mathcal}
\newcommand{\mf}{\mathfrak}

\newcommand{\J}{\mathbb{J}}
\newcommand{\psp}{\psi}
\newcommand{\psm}{\psi^\ast}

\newcommand{\om}{\omega}

\numberwithin{equation}{section}

\begin{document}
\title[ ]
{Super duality and Crystal bases for quantum ortho-symplectic superalgebras II}
\author{JAE-HOON KWON}
\address{ Department of Mathematics \\ Sungkyunkwan University \\ Suwon,  Republic of Korea}
\email{jaehoonkw@skku.edu} 

\thanks{This work was  supported by Samsung Science and Technology Foundation under Project Number SSTF-BA1501-01.}

\begin{abstract}
Let $\mc{O}^{int}_q(m|n)$ be a semisimple tensor category of modules over a quantum ortho-symplectic superalgebra  of type $B, C, D$ introduced in \cite{K13}. It is a natural counterpart of the category of finitely dominated integrable modules over a quantum group of type $B, C, D$ from a viewpoint of super duality. Continuing the previous work on type $B$ and $C$ \cite{K13}, we classify the irreducible modules in $\mc{O}^{int}_q(m|n)$, and prove the existence and uniqueness of their crystal bases in case of type $D$. A new combinatorial model of classical crystals of type $D$ is introduced, whose super analogue gives a realization of crystals for the highest weight modules in $\mc{O}^{int}_q(m|n)$. 
\end{abstract}

\maketitle
\setcounter{tocdepth}{1}

\section{Introduction}
This is a continuation of our previous work \cite{K13} on crystal bases for quantum ortho-symplectic superalgebras. 
In \cite{K13}, we constructed a semisimple tensor category $\mc{O}^{int}_q(m|n)$ of modules over an quantum superalgebra $U_q({\mf g}_{m|n})$, where $\g_{m|n}$ is an ortho-symplectic Lie superalgebra of type $B, C, D$ or $\g ={\mf b},  {\mf b}^\bullet, {\mf c}, {\mf d}$. This category is characterized by remarkably simple conditions similar to those for  the polynomial $U_q(\gl_{m|n})$-modules (cf. \cite{BR,CK}), while its irreducible modules are $q$-deformations of infinite-dimensional $\g_{m|n}$-modules appearing in a tensor power of a Fock space, which were studied in \cite{CKW} via Howe duality. 
Its semisimplicity is based on the fact that $\mc{O}^{int}_q(m|n)$ naturally corresponds  to the semisimple tensor category $\mc{O}^{int}_q(m+n)$ of finitely dominated integrable modules over a quantum enveloping algebra of the corresponding classical Lie (super)algebra $\g_{m+n}$ from a viewpoint of super duality \cite{CLW,CW} (more precisely, they are equivalent when $n=\infty$ and $q=1$). 

Motivated by the work on crystal bases of polynomial $U_q(\gl_{m|n})$-modules \cite{BKK}, we classified the irreducible modules in $\mc{O}^{int}_q(m|n)$ when $\g = {\mf b},  {\mf b}^\bullet$ and ${\mf c}$, and then proved the existence and uniqueness of their crystal bases, where the associated crystal is realized in terms of a new combinatorial object called ortho-symplectic tableaux of type $B$ and $C$, respectively \cite{K13}. 

In this paper, we establish the same result for $\g ={\mf d}$ (Theorem \ref{Existence of crystal base}). The strategy for the case $\g=\mf d$ is parallel to that of \cite{K13}. But the main ingredient of our proof different from \cite{K13} is to formulate the notion of an ortho-symplectic tableau of type $D$ (Definitions \ref{Definition of T_A(a)}, \ref{admissible}, and \ref{def:osptableaux}), where 
more technical difficulty enters compared to type $B$ and $C$.


An ortho-symplectic tableau of type $D$ is a sequence of two-column shaped skew tableaux of type $A$ with certain admissibility conditions on adjacent pairs similar to \cite{K13}. Then the crystal of an irreducible highest weight $U_q({\mf d}_{m|n})$-module in $\mc{O}^{int}_q(m|n)$ is realized as the set of ortho-symplectic tableaux of type $D$ associated to its highest weight, where the underlying tableaux of type $A$ are semistandard tableaux for  $\gl_{m|n}\subset \mf{d}_{m|n}$ \cite{BR}. Furthermore, when we replace the underlying tableaux of type $A$ with usual semistandard tableaux for $\gl_{m}\subset \mf{d}_{m}$ (by putting $n=0$), we obtain a new realization of crystals of integrable highest weight $U_q(\mf{d}_{m})$-modules in $\mc{O}^{int}_q(m)$ (type $D_{m}$) (Theorems \ref{crystal invariance of osp tableaux} and \ref{character formula for m+n}), which plays a crucial role in this paper.

We remark that the tableaux models for type $BCD$ introduced in this paper and \cite{K13} is based on the Fock space model (cf. \cite{CKW,Ha}), while the well-known Kashiwara-Nakashima tableaux (for non-spinor highest weights) \cite{KashNaka} are based on the crystals of natural representation and its tensor powers. We expect that our new combinatorial model for classical crystals is of independent interest and can be used for other interesting applications in the future.

Finally, combining with the results in \cite{BKK} for type $A$ and \cite{K13} for type $B$ and $C$, we conclude that the super duality, when restricted to the integrable modules over the classical Lie algebras, provides a natural semisimple tensor category for Lie  superalgebras of types $ABCD$, where a crystal base theory exists. We expect that this can be extended to a more general class of contragredient Lie superalgebras with isotropic simple roots, which includes simple finite-dimensional Lie superalgebras of exceptional types $F(3|1)$, $G(2)$ and $D(2|1,\alpha)$ ($\alpha\in\mathbb{Z}_{>0}$), and whose super duality has been established in \cite{CKW-2} recently.

The paper is organized as follows. In Section 2, we briefly recall the notations and results in \cite{K13}. In Section 3, we introduce our main combinatorial object called ortho-symplectic tableaux of type $D$. Then in Section 4, we prove that the set of ortho-symplectic tableaux associated to a given highest weight gives the character of the corresponding irreducible highest weight module over the Lie superalgebra $\mf{d}_{m|n}$. 
In Section 5, we classify the irreducible $U_q(\mf{d}_{m|n})$-modules in $\mc{O}^{int}_q(m|n)$, and prove the existence and uniqueness of their crystal bases, where the crystals are realized in terms of ortho-symplectic tableaux of type $D$. In Section 6, we give a proof of Theorem \ref{crystal invariance of osp tableaux}, which is a main result in Section 4.
 
\vskip 3mm

\section{Quantum superalgebra $U_q(\mf{d}_{m|n}$) and the category $\mc{O}^{int}_q(m|n)$}\label{Lie superalgebra}

\subsection{Notations} 
Throughout this paper, we assume that $m\in\mathbb{Z}_{>0}$ with $m\geq 2$ and $n\in \mathbb{Z}_{>0}\cup\{\infty\}$.
Let us recall the following notations for the classical Lie superalgebra $\mf{d}_{m|n}$ of type $D$  in \cite[Section 2]{K13}:

\begin{itemize}
\item[$\cdot$] $\J_{m|n}=\{\,\ov{m}< \ldots <\ov{2}< \ov{1} <\tfrac{1}{2}<\tfrac{3}{2}<\ldots<n-\tfrac{1}{2}\,\}$,

\item[$\cdot$] $P_{m|n}=\bigoplus_{a\in \J_{m|n}}\Z \delta_a \oplus \Z \Lambda_{\ov{m}}$ : the weight lattice, 

\item[$\cdot$] $P^\vee_{m|n}=\bigoplus_{a\in \J_{m|n}}\Z E_a \oplus \Z K'$ : the dual weight lattice,

\item[$\cdot$] $I_{m|n}=\{\,\ov{m},\ldots,\ov{1},0,\frac{1}{2},\ldots,n-\frac{3}{2}\,\}$,

\item[$\cdot$] $\Pi_{m|n}=\{\,\beta_i\,|\,i\in I_{m|n}\,\}$ : the set of simple roots,

\item[$\cdot$] $\Pi^\vee_{m|n}=\{\,\beta_i^\vee\,|\,i\in I_{m|n}\,\}$ : the set of simple coroots, where
{\allowdisplaybreaks
\begin{align*}
\beta_i&=
\begin{cases}
-\delta_{\ov{m}}-\delta_{\ov{m-1}}, & \text{if $i=\ov{m}$},\\
\delta_{\ov{k+1}}-\delta_{\ov{k}}, & \text{if $i=\ov{m-1},\ldots,\ov{1}$},\\
\delta_{\ov{1}}-\delta_{\hf}, & \text{if $i=0$},\\
\delta_i-\delta_{i+1}, & \text{if $i=\frac{1}{2},\ldots,n-\frac{3}{2}$},
\end{cases}\\
\beta_i^\vee&=
\begin{cases}
-E_{\ov{m}}-E_{\ov{m-1}}+K', & \text{if $i=\ov{m}$},\\
E_{\ov{k+1}}-E_{\ov{k}}, & \text{if $i=\ov{m-1},\ldots,\ov{1}$},\\
E_{\ov{1}}+E_{\hf}, & \text{if $i=0$},\\
E_i-E_{i+1}, & \text{if $i=\frac{1}{2},\ldots,n-\frac{3}{2}$},
\end{cases}
\end{align*}}

\item[$\cdot$] $I_{m|0}=\{\,\ov{m},\ldots,\ov{1}\,\}$ and $I_{0|n}=  \{\,\frac{1}{2},\ldots,n-\frac{3}{2}\,\}$.
\end{itemize}

\noindent Here, $\J_{m|n}$ is a $\Z_2$-graded set with $(\J_{m|n})_0=\{\,\ov{m},\ldots,\ov{1}\,\}$ and $(\J_{m|n})_1=\{\,1/2, \ldots, n-1/2\,\}$, and we write $|a|=\varepsilon$ for $a\in (\J_{m|n})_\varepsilon$ and $\epsilon\in \Z_2$. We assume that $\{\,\Lambda_{\ov{m}}, \delta_a\, (a\in \J_{m|n})\,\}$ and $\{\, K', E_a\, (a\in \J_{m|n}) \,\}$ are dual bases with respect to the natural pairing $\langle \,\cdot\, ,\,\cdot\,\rangle$ on $P^\vee_{m|n}\times P_{m|n}$, that is, 
\begin{equation*}
\langle E_b,\delta_a\rangle=\delta_{ab},\ \ \langle K',\delta_a\rangle=0,\ \ \langle E_a,\Lambda_{\ov{m}}\rangle=0,\ \ \langle K',\Lambda_{\ov{m}}\rangle=1,
\end{equation*}
for $a, b\in \J_{m|n}$, and  $\h_{m|n}^*:=\C\otimes_{\Z} P_{m|n}$ has a symmetric bilinear form $(\,\cdot\,|\,\cdot\,)$  given by
\begin{equation*}\label{bilinear form spo}
\begin{split}
&(\lambda|\delta_a)=\big\langle (-1)^{|a|}E_a -K, \lambda \big\rangle, \ \ \ (\Lambda_{\ov{m}}|\Lambda_{\ov{m}})=0,
\end{split}
\end{equation*}
for $a,b\in  \J_{m|n}$ and $\lambda\in   \h_{m|n}^*$. For $i\in I_{m|n}$, let $s_i=1$ for $i\in\{\ov{m},\ldots,\ov{1},0\}$, and $-1$ otherwise. Then $s_j \langle \beta_j^\vee, \mu \rangle = (\beta_j | \mu)$ for $j\in I_{m|n}$, $\mu\in \mf{h}^*_{m|n}$. Following \cite{Kac77}, the Dynkin diagram associated with the Cartan matrix $A=(a_{ij})=(\langle \beta_i^\vee,\beta_j\rangle)_{i,j\in I_{m|n}}$ is
\begin{center}
\hskip -3cm \setlength{\unitlength}{0.16in} \medskip
\begin{picture}(24,5.8)
\put(6,0){\makebox(0,0)[c]{$\bigcirc$}}
\put(6,4){\makebox(0,0)[c]{$\bigcirc$}}
\put(8,2){\makebox(0,0)[c]{$\bigcirc$}}
\put(10.4,2){\makebox(0,0)[c]{$\bigcirc$}}
\put(14.85,2){\makebox(0,0)[c]{$\bigcirc$}}
\put(17.25,2){\makebox(0,0)[c]{$\bigotimes$}}
\put(19.4,2){\makebox(0,0)[c]{$\bigcirc$}}
\put(21.5,2){\makebox(0,0)[c]{$\bigcirc$}}
\put(6.35,0.3){\line(1,1){1.35}} \put(6.35,3.7){\line(1,-1){1.35}}
\put(8.4,2){\line(1,0){1.55}} \put(10.82,2){\line(1,0){0.8}}
\put(13.2,2){\line(1,0){1.2}} \put(15.28,2){\line(1,0){1.45}}
\put(17.7,2){\line(1,0){1.25}} \put(19.8,2){\line(1,0){1.25}}
\put(21.95,2){\line(1,0){1.4}} 
\put(12.5,1.95){\makebox(0,0)[c]{$\cdots$}}
\put(24.5,1.95){\makebox(0,0)[c]{$\cdots$}}
\put(6,5){\makebox(0,0)[c]{\tiny $\beta_{\ov{m}}$}}
\put(6,-1.2){\makebox(0,0)[c]{\tiny $\beta_{\ov{m-1}}$}}
\put(8.2,1){\makebox(0,0)[c]{\tiny $\beta_{\ov{m-2}}$}}
\put(10.4,1){\makebox(0,0)[c]{\tiny $\beta_{\ov{m-3}}$}}
\put(14.9,1){\makebox(0,0)[c]{\tiny $\beta_{\ov{1}}$}}
\put(17.15,1){\makebox(0,0)[c]{\tiny $\beta_0$}}
\put(19.5,0.8){\makebox(0,0)[c]{\tiny $\beta_{\frac{1}{2}}$}}
\put(21.5,0.8){\makebox(0,0)[c]{\tiny $\beta_{\frac{3}{2}}$}}
\end{picture}\vskip 8mm
\end{center}

For $\Lambda=c\Lambda_{\ov{m}}+\sum_{a\in \J_{m|n}}\lambda_a\delta_a\in P_{m|n}$, we assume that the parity of  $\Lambda$ is $\sum_{a\geq \hf}\lambda_a$ $\pmod 2$, which we denote by $|\Lambda|$. In particular, we have $|\beta_i|=0$ for $i\neq 0$ and  $|\beta_0|=1$.

\subsection{The quantum superalgebra $U_q(\mf{d}_{m|n})$}
Let $q$ be an indeterminate. For $r\geq 0$, put $[r]=\frac{q^r-q^{-r}}{q-q^{-1}}$ and $[r]!=\prod_{k=1}^r[k]$. For $i\in I_{m|n}$, put $q_i = q^{{\ov s}_i}$, where ${\ov s}_i=-s_i$. The  quantum superalgebra $U_q({\mf d}_{m|n})$ is the
associative superalgebra (or $\Z_2$-graded algebra) with $1$ over $\mathbb{Q}(q)$ generated by
$e_i$, $f_i$ $(\,i\in I_{m|n}\,)$ and $q^h$ $(\,h\in P^{\vee}_{m|n}\,)$, which are
subject to the following relations \cite{Ya}:
{\allowdisplaybreaks
\begin{align*}
\ \ \ & {\rm deg}(q^h)=0,\ \  {\rm deg}(e_i)={\rm deg}(f_i)=|\beta_i|,\\
& q^0=1, \quad q^{h +h'}=q^{h}q^{h'}, \ \ q^h e_i=q^{\langle h,\beta_i\rangle}
e_i q^h, \quad q^h f_i=q^{-\langle h,\beta_i\rangle} f_i q^h, \\
& e_i f_j-(-1)^{|\beta_i||\beta_j|}f_j e_i =\delta_{ij}\frac{t_i-t_i^{-1}}{q_i-q^{-1}_i}, \\
& z_i z_j - (-1)^{|\beta_i||\beta_j|} z_j z_i =0,     \ \ \ \ \ \ \ \ \ \ \ \ \ \ \ \ \ \ \ \ \ \, \text{if $(\beta_i|\beta_j) =0$},
\\ & z_i^2z_j - (q+q^{-1})z_iz_jz_i +z_jz_i^2=0,  \ \ \ \ \ \ \ \ \ \ \  \text{if $i\neq 0$ and $(\beta_i|\beta_j) \neq 0$}, \\
& z_0 z_{\ov{1}} z_0 z_{\hf} + z_{\ov{1}} z_0
z_{\hf} z_{0} + z_{0} z_{\hf} z_{0} z_{\ov{1}}  + z_{\hf} z_{0} z_{\ov{1}} z_{0} -(q+q^{-1}) z_{0} z_{\ov{1}} z_{\hf} z_{0} =0,
\end{align*}
}
\hskip -1.5mm for $i,j\in I_{m|n}$, $h, h' \in P^{\vee}_{m|n}$ and $z=e,f$, where $t_i = q^{{\ov s}_i \beta_i^\vee}$.
Recall that there is a Hopf superalgebra structure on $U_q(\mf{d}_{m|n})$, where the comultiplication $\Delta$ is given by
$\Delta(q^h)=q^h\otimes q^h$, $\Delta(e_i)=e_i\otimes t_i^{-1} + 1\otimes e_i$, $\Delta(f_i) =f_i\otimes 1+ t_i\otimes f_i$, 
the antipode $S$ is given by
$S(q^h)=q^{-h}$, $S(e_i)=-e_it_i$, $S(f_i)=-t_i^{-1}f_i$,
and the counit $\varepsilon$ is given by  $\varepsilon(q^h)=1$, $\varepsilon(e_i)=\varepsilon(f_i)=0$ for $h\in P^\vee_{m|n}$ and  $i\in I_{m|n}$.

\subsection{The category $\mathcal{O}_q^{int}(m|n)$}

Let $\cP$ be the set of partitions and let 
\begin{equation*}
\cP(\mf{d})=\{\,(\lambda,\ell)\in \cP\times\Z_{> 0}\,|\,\ell-\lambda_1-\lambda_2\in \Z_{\geq 0}\,\},
\end{equation*}
where $\lambda=(\lambda_1,\lambda_2,\ldots)$.
For $(\lambda,\ell)\in \cP(\mf{d})$, let
\begin{equation*}
\Lambda_{m|\infty}(\lambda,\ell)=\ell\Lambda_{\ov{m}} + \lambda_1\delta_{\ov{m}}+\cdots+\lambda_m\delta_{\ov{1}}+\mu_{1}\delta_{\frac{1}{2}}+\mu_{2}\delta_{\frac{3}{2}}+\cdots,
\end{equation*}
where $\mu=(\mu_1,\mu_2,\ldots)=(\lambda_{m+1},\lambda_{m+2},\ldots)'$,  the conjugate partition of $(\lambda_{m+1},\lambda_{m+2},\ldots)$. Put $\cP(\mf{d})_{m|n}=\{\,(\lambda,\ell)\in \cP(\mf{d})\,|\,\Lambda_{m|\infty}(\lambda,\ell)\in P_{m|n}\,\}$. For $(\lambda,\ell)\in \cP(\mf{d})_{m|n}$, we write $\Lambda_{m|n}(\lambda,\ell)=\Lambda_{m|\infty}(\lambda,\ell)$. 

Let $\mathcal{O}_q^{int}(m|n)$ be  the category of $U_q(\mf{d}_{m|n})$-modules $M$ satisfying
\begin{itemize}
\item[(1)] $M=\bigoplus_{\gamma\in P_{m|n}}M_\gamma$ and $\dim M_\gamma <\infty$ for $\gamma\in P_{m|n}$,

\item[(2)] ${\rm wt}(M)\subset \bigcup_{i=1}^r\left(\ell_i\Lambda_{\ov{m}}+\sum_{a\in \J_{m|n}}\Z_{\geq 0}\delta_a \right)$ for some $r\geq 1$ and  $\ell_i\in \Z_{\geq 0}$,

\item[(3)] $f_{\ov{m}}$ acts locally nilpotently on $M$,
\end{itemize}
where ${\rm wt}(M)$ denotes the set of weights of $M$.
For $\Lambda\in P_{m|n}$, let $L_q(\mf{d}_{m|n},\Lambda)$ denote the irreducible highest weight $U_q(\mf{d}_{m|n})$-module with highest weight $\Lambda$.
By \cite[Theorems 4.2 and 4.3]{K13}, we have the following.

\begin{thm}\label{complete reducibility}
$\mathcal{O}^{int}_q(m|n)$ is a semisimple tensor category, and any highest weight module in $\mathcal{O}^{int}_q(m|n)$ is isomorphic to $L_q({\mf d}_{m|n},\Lambda_{m|n}(\lambda,\ell))$ for some $(\lambda,\ell)\in \cP(\mf{d})_{m|n}$.
\end{thm}

\subsection{Crystal base}\label{crystal base}
Let $M$ be a $U_q(\mf{d}_{m|n})$-module in $\mathcal{O}_q^{int}(m|n)$.
Let $i\in  I_{m|n}$ be given.
Suppose that $i\in I_{m|n}\setminus\{0\}$. 
For $u\in M$ of weight $\gamma$, we have a unique expression
\begin{equation*}
u=\sum_{k\geq 0, -\langle \beta^\vee_i,\gamma \rangle}f_i^{(k)}u_k,
\end{equation*}
where $e_iu_k=0$ for all $k\ge 0$. We define the Kashiwara operators $\te_i$ and $\tf_i$ as follows:\begin{equation*}
\begin{split}
\te_i u &=
\begin{cases}
\sum_{k}q_i^{l_k-2k+1}f_i^{(k-1)}u_k, & \text{if $i\in I_{m|0}$},\\
\sum_{k}f_i^{(k-1)}u_k, & \text{if $i \in I_{0|n}$},
\end{cases}
\\ 
\tf_i u &=
\begin{cases}
\sum_{k}q_i^{-l_k+2k+1}f_i^{(k+1)}u_k, & \text{if $i\in I_{m|0}$},\\
\sum_{k}f_i^{(k+1)}u_k, & \text{if $i \in I_{0|n}$},
\end{cases}
\end{split}
\end{equation*}
where  $l_k=\langle \beta^\vee_i,\gamma+k\beta_i \rangle$ for $k\geq 0$. If $i=0$, then we define
\begin{equation*}
\begin{split}
\te_0 u =e_0u, \ \ \ \
\tf_0 u =q_0 f_0t_0^{-1}u .
\end{split}
\end{equation*}\vskip 2mm

Let $\mathbb{A}$ denote the subring of $\mathbb{Q}(q)$ consisting
of all rational functions which are regular at $q=0$.
We call a pair $(L,B)$ a {\it crystal base of $M$} if
\begin{itemize}
\item[(1)] $L$ is an $\mathbb{A}$-lattice of $M$, where  $L=\bigoplus_{\gamma\in P_{m|n}}L_{\gamma}$ with $L_{\gamma}=L\cap
M_{\gamma}$,
\item[(2)] $\tilde{e}_i L\subset L$ and $\tilde{f}_i L\subset L$
for $i\in I_{m|n}$,

\item[(3)] $B$ is a pseudo-basis of $L/qL$ (i.e.
$B=B^{\bullet}\cup(-B^{\bullet})$ for a $\mathbb{Q}$-basis
$B^{\bullet}$ of $L/qL$),

\item[(4)] $B=\bigsqcup_{\gamma\in P_{m|n}}B_{\gamma}$ with
$B_{\gamma}=B\cap(L/qL)_{\gamma}$,

\item[(5)] $\tilde{e}_iB \subset B\sqcup \{0\}$,
$\tilde{f}_i B\subset B\sqcup \{0\}$ for $i\in I_{m|n}$,

\item[(6)] for $b,b'\in B$ and $i\in I_{m|n}$,
$\tilde{f}_i b = b'$ if and only if $b=\tilde{e}_i b'$
\end{itemize}
(see \cite{BKK}).
The set $B / \{\pm 1 \}$ has an $I_{m|n}$-colored oriented graph
structure, where  $b\stackrel{i}{\rightarrow}
b'$ if and only if $\tf_i b=b'$ for $i\in I_{m|n}$ and  $\,b, b' \in B / \{\pm
1\}$. We call $B / \{\pm 1 \}$ the {\it crystal} of $M$.
For $b\in B$ and $i\in I_{m|n}$, we set
$\varepsilon_i(b)= \max\{\,r\in\mathbb{Z}_{\geq
0}\,|\,\tilde{e}_i^r b\neq 0\,\}$ and $\varphi_i(b)=
\max\{\,r\in\mathbb{Z}_{\geq 0}\,|\,\tilde{f}_i^r b\neq 0\,\}$. We denote the weight of $b$ by ${\rm wt}(b)$.

Let $M_i$ $(\,i=1,2\,)$ be a
$U_q(\mf{d}_{m|n})$-module in $\mathcal{O}^{int}_q(m|n)$
with a crystal base $(L_i,B_i)$. Then  $(L_1\otimes L_2,B_1\otimes B_2)$ is a crystal
base of $M_1\otimes M_2$ \cite[Proposition 2.8]{BKK}. The actions of $\te_i$ and $\tf_i$ on $B_1\otimes B_2$ are as follows.

For $i\in I_{0|n}$, we have {\allowdisplaybreaks
\begin{equation}\label{lower tensor product rule}
\begin{split}
&\tilde{e}_i(b_1\otimes b_2)= \begin{cases}
(\tilde{e}_i b_1) \otimes b_2, & \text{if $\varphi_i(b_1)\geq\varepsilon_i(b_2)$}, \\
b_1 \otimes (\tilde{e}_i b_2), & \text{if $\varphi_i(b_1)<\varepsilon_i(b_2)$},\\
\end{cases}
\\
&\tilde{f}_i(b_1\otimes b_2)=
\begin{cases}
(\tilde{f}_i b_1) \otimes b_2, & \text{if $\varphi_i(b_1)>\varepsilon_i(b_2)$}, \\
 b_1 \otimes (\tilde{f}_i b_2), & \text{if $\varphi_i(b_1)\leq\varepsilon_i(b_2)$}.
\end{cases}
\end{split}
\end{equation}}

For $i\in I_{m|0}$, we have
{\allowdisplaybreaks
\begin{equation}\label{upper tensor product rule}
\begin{split}
&\tilde{e}_i(b_1\otimes b_2)= \begin{cases}
 b_1 \otimes (\tilde{e}_ib_2), & \text{if $\varphi_i(b_2)\geq\varepsilon_i(b_1)$}, \\
(\tilde{e}_i b_1) \otimes b_2, & \text{if $\varphi_i(b_2)<\varepsilon_i(b_1)$},\\
\end{cases}
\\
&\tilde{f}_i(b_1\otimes b_2)=
\begin{cases}
b_1 \otimes (\tilde{f}_i b_2), & \text{if $\varphi_i(b_2)>\varepsilon_i(b_1)$}, \\
(\tilde{f}_i b_1) \otimes  b_2, & \text{if $\varphi_i(b_2)\leq\varepsilon_i(b_1)$}.
\end{cases}
\end{split}
\end{equation}}

For $i=0$, we have{\allowdisplaybreaks
\begin{equation}\label{tensor product rule for 0}
\begin{split}
\tilde{e}_0(b_1\otimes b_2)=&
\begin{cases}
\pm b_1\otimes (\tilde{e}_0 b_2), & \text{if }\langle \beta^\vee_0,{\rm
wt}(b_2)\rangle>0, \\ (\tilde{e}_0 b_1)\otimes  b_2, & \text{if
}\langle \beta^\vee_0,{\rm wt}(b_2)\rangle=0,
\end{cases}
\\
\tilde{f}_0(b_1\otimes b_2)=&
\begin{cases}
\pm b_1\otimes (\tilde{f}_0 b_2), & \text{if }\langle \beta^\vee_0,{\rm
wt}(b_2)\rangle>0, \\ (\tilde{f}_0 b_1)\otimes  b_2, & \text{if
}\langle \beta^\vee_0,{\rm wt}(b_2)\rangle=0.
\end{cases}
\end{split}
\end{equation}}

\subsection{Fock space}\label{subsec:Fock space}
Let $\A^+_q$ be an associative $\mathbb{Q}(q)$-algebra with $1$ generated by $\psp_{a}$, $\psm_a$, $\om_a$ and $\om^{-1}_a$ for $a\in \J_{m|n}$
subject to the following relations:
{\allowdisplaybreaks
\begin{align*}
&\om_a \om_b= \om_b\om_a, \ \ \om_a\om_a^{-1}=1, \\
&\om_a \psp_b \om_a^{-1} =
q^{(-1)^{|a|}\delta_{ab}}\psp_b ,\ \ \ \ \ \om_a \psm_b \om_a^{-1} =
q^{-(-1)^{|a|}\delta_{ab}}\psm_b,  \\
&\psp_a\psp_b+(-1)^{|a||b|}\psp_b\psp_a=0,\ \ \ \ \psm_a\psm_b+(-1)^{|a||b|}\psm_b\psm_a=0 , \\
&\psp_a\psm_b+(-1)^{|a||b|}\psm_b\psp_a=0 \ \ \ \ (a\neq b), \\
&\psp_a\psm_a =[q\om_a],\ \ \
\psm_a\psp_a =
(-1)^{1+|a|}[\om_a].
\end{align*}}
\noindent Here $[q^k\om_a^{\pm 1}]=\frac{q^k\om^{\pm 1}_a-q^{-k}\om_a^{\mp 1}}{q-q^{-1}}$ for $k\in\Z$ and $a\in \J_{m|n}$ (cf. \cite{Ha}).  Let 
\begin{equation}\label{Fock space}
\mathscr{V}_q : =\A^+_q  \, |0\rangle
\end{equation}
be the $\A^+_q$-module generated by $|0\rangle$ satisfying $\psm_b|0\rangle =0$ and $\omega_b|0\rangle =q^{-1}|0\rangle$ 
for $a, b\in \J_{m|n}$. 
Let $\bf{B}^+$ be the set of sequences ${\bf m}=(m_a)$ of non-negative integers indexed by $\J_{m|n}$ such that $m_a\leq 1$ for $|a|=1$. For ${\bf m}=(m_a)\in \bf{B}^+$,  let
\begin{equation*}
\begin{split}
&\psi_{{\bf m}}= \overrightarrow{\prod_{a\in\J_{m|n}}}{\psp_a}^{(m_a)},
\end{split}
\end{equation*}
where the product is taken in the order on $\J_{m|n}$ and $\psp_a^{(r)}=(\psp_a)^{r}/[r]!$, ${\psm_a}^{(r)}= (\psm_a)^{r}/[r]!$.
By similar arguments as in \cite[Proposition 2.1]{Ha}, we can check that $\mathscr{V}_q$ is an irreducible $\A^+_q$-module with a $\mathbb{Q}(q)$-linear basis $\{\,\psi_{{\bf m}}|0\rangle\,|\,{\bf m}\in \B^+\,\}$.

It is shown in \cite[Proposition 5.3]{K13} that $\mathscr{V}_q$ has a $U_q(\mf{d}_{m|n})$-module structure, where ${\rm wt}(\psi_{{\bf m}})=\Lambda_{\ov{m}}+\sum_{a\in \J_{m|n}}m_a\delta_a$.
Since $\mathscr{V}_q\in \mc{O}^{int}_q(m|n)$, $\mathscr{V}_q^{\otimes \ell}$ is completely reducible by Theorem \ref{complete reducibility} for $\ell\geq 1$.
Also by \cite[Theorem 5.6]{K13}, $\mathscr{V}_q$ has a crystal base $(\mathscr{L}^+,\mathscr{B}^+)$, where
\begin{equation}\label{crystal of V_q}
\mathscr{L}^+=\sum_{\bf{m}\in \B^+}\mathbb{A} \psp_{\bf{m}}|0\rangle, \ \ \
\mathscr{B}^+=\left\{\, \pm\psp_{\bf{m}}|0\rangle \!\!\!\! \pmod{q\mathscr{L}}\,|\,\bf{m}\in \B^+\, \right\}.
\end{equation}

\section{Ortho-symplectic tableaux of type $D$}\label{ortho-symplectic Tableaux of Type D}

\subsection{Semistandard tableaux}

Let us recall some basic terminologies and notations related with tableaux. We refer the reader to \cite[Section 3.1]{K13}. 
We assume that  $\mc{A}$ is a linearly ordered countable set with a $\mathbb{Z}_2$-grading $\mc{A}=\mc{A}_0\sqcup\mc{A}_1$. When $\mc{A}$ is (a subset of) $\mathbb{Z}_{>0}$, we assume that $\mc{A}_0=\mc{A}$ with the usual linear ordering.  For a skew Young diagram $\lambda/\mu$, we denote by $SST_{\mc{A}}(\lambda/\mu)$ be the set of $\mc{A}$-semistandard tableaux of shape $\lambda/\mu$. 
For $T\in SST_{\mc{A}}(\lambda/\mu)$, ${\rm sh}(T)$ denotes the shape of $T$, ${\rm wt}(T)=(m_a)_{a\in\mc{A}}$ is the weight of $T$, where $m_a$ is the number of occurrences of $a$ in $T$, and $w(T)$ is the word given by reading the entries of $T$ column by column from right to left and from top to bottom in each column.

For $T\in SST_\mc{A}(\lambda)$ and
$a\in \mc{A}$, we denote by $a \rightarrow T$ the tableau obtained by the column insertion of $a$ into $T$ (cf. \cite{BR,Ful}). For a finite word $w=w_1\ldots w_r$ with letters in $\mc{A}$, we define $(w \rightarrow T)=(w_r\rightarrow(\cdots(w_1\rightarrow T)))$. For an $\mc{A}$-semistandard tableau $S$, we define $(S\rightarrow T)=(w(S)\rightarrow T)$.  

For an $\mc{A}$-semistandard tableau $S$ of single-columned shape, we denote by $S(i)$ ($i\geq 1$) the $i$-th entry from the bottom, and by ${\rm ht}(S)$ the height of $S$.

\subsection{Signature on tableaux of two-columned shape}\label{subsec:signature}
Let $\sigma=(\sigma_{1},\ldots,\sigma_{r})$ be a sequence with $\sigma_{i}\in \{\,+\,,\,-\, , \ \cdot\ \}$. We denote by $\tilde{\sigma}$ the sequence obtained from $\sigma$ by signature rule, that is, by applying the process of replacing any adjacent pair $(+,-)$ (ignoring $\cdot$ between them) with $(\,\cdot\,,\,\cdot\,)$ as far as possible.  

Let $U$ and $V$ be two $\mc{A}$-semistandard tableaux of single-columned shapes. Let $w=w_1\ldots w_r$ be the word given by rearranging the entries of $U$ and $V$ in weakly decreasing order. Here we assume that for $a\in \mc{A}$ occurring in both $U$ and $V$, we put the letters $a$ in $U$ to the right (resp. left) of those in $V$ if $a$ is even or $a\in \mc{A}_0$ (resp. $a$ is odd or $a\in \mc{A}_{1}$). Let $\sigma=(\sigma_{1},\ldots,\sigma_{r})$ be a sequence given by $\sigma_i=+$ if $w_i$ comes from $V$ and $-$ if $w_i$ comes from $U$. We define
\begin{equation}
\sigma(U,V)=(p,q),
\end{equation}
where $p$ (resp. $q$) is the number of $-$'s (resp. $+$'s) in $\tilde{\sigma}$, and call it the {\it signature of $(U,V)$}. 

For $a,b,c\in\Z_{\geq 0}$, let $\lambda(a,b,c)=(2^{b+c},1^a)/(1^b)$, which is a skew Young diagram with two columns of heights $a+c$ and $b+c$. For example,\vskip 1mm
$$\resizebox{.18\hsize}{!}{$
\lambda(2,1,3)\ =\ {\def\lr#1{\multicolumn{1}{|@{\hspace{.75ex}}c@{\hspace{.75ex}}|}{\raisebox{-.04ex}{$#1$}}}\raisebox{-.6ex}
{$\begin{array}{cc}
\cline{2-2}
\cdot &\lr{}\\ 
\cline{1-1}\cline{2-2}
\lr{} &\lr{}\\ 
\cline{1-1}\cline{2-2}
\lr{} &\lr{\ \ }\\ 
\cline{1-1}\cline{2-2}
\lr{}&\lr{}\\ 
\cline{1-1}\cline{2-2}
\lr{\ \ }& \cdot \\ 
\cline{1-1} 
\lr{} & \cdot \\ 
\cline{1-1}
\end{array}$}}$}$$\vskip 1mm
Let  $T$ be a tableau of shape $\lambda(a,b,c)$, whose column is $\mc{A}$-semistandard. We denote by $T^{\tt L}$ and $T^{\tt R}$ the left and right columns of $T$, respectively.
Then 
\begin{equation}\label{signature condition}
\text{$T\in SST_{{\mc A}}(\lambda(a,b,c))$\ \ if and only if \ 
 $\sigma(T^{\tt L},T^{\tt R})=(a-p,b-p)$}
\end{equation}
for some $0\leq p\leq \min\{a,b\}$ \cite[Lemma 6.2]{K13}.

\subsection{RSK and signatures}
For $\ell\in \mathbb{Z}_{>0}$, let ${\bf M}_{\mc{A}\times \ell}$ be the set of $\mc{A}\times\ell$ matrices ${\bf m}=[m_{a i}]$ ($a\in \mc{A}$, $i=1,\ldots, \ell$) such that (1) $m_{a i}\in \Z_{\geq 0}$, (2) $m_{a i}\leq 1$ for $|a|=0$, (3) $\sum_{a,i}m_{a i}< \infty$. Let ${\bf m}\in {\bf M}_{\mc{A}\times \ell}$ be given. 
For $1\leq k\leq \ell$, let ${\bf m}^{(k)}=[m_{a k}]$ denote the $k$th column of ${\bf m}$, and $|{\bf m}^{(k)}|=\sum_{a\in \mc{A}}m_{ak}$. We will often identify  each ${\bf m}^{(k)}$ with an $\mc{A}$-semistandard tableau of single-columned shape $(1^{|{\bf m}^{(k)}|})$. Similarly, for $a\in \mc{A}$, let ${\bf m}_{(a)}$ be the $a$th row of ${\bf m}$, and $|{\bf m}_{(a)}|=\sum_{1\leq k\leq \ell}m_{ak}$. We put $|{\bf m}|=\sum_{a}|{\bf m}_{(a)}|=\sum_{k}|{\bf m}^{(k)}|$. We remark that our convention of column indices are increasing from right to left  so that ${\bf m}=[\,{\bf m}^{(\ell)}:\cdots:{\bf m}^{(1)}\,]$.

For $1\leq k\leq \ell$, let 
$P({\bf m})^{(k)}= ({\bf m}^{(k)}\rightarrow (\cdots({\bf m}^{(2)}\rightarrow {\bf m}^{(1)})))$, and let $\lambda^{(k)}={\rm sh}(P({\bf m})^{(k)})$. Put $P({\bf m})=P({\bf m})^{(\ell)}$ and $\lambda=\lambda^{(\ell)}$.  Let $Q({\bf m}) \in SST_{\{1,\ldots,\ell\}}(\lambda')$ be such that the subtableau of shape $\lambda^{(k)}/\lambda^{(k-1)}$ is filled with $k$ for $1\leq k\leq \ell$, where $\lambda'$ is the conjugate of $\lambda$ and $\lambda^{(0)}$ is the empty Young diagram. Then the map ${\bf m}\mapsto (P({\bf m}), Q({\bf m}))$ gives a bijection
\begin{equation}\label{RSK}
{\bf M}_{\mc{A}\times \ell} \longrightarrow \bigsqcup_{\substack{\lambda\in \cP \\ \lambda_1\leq \ell}} SST_{\mc{A}}(\lambda)\times SST_{\{1,\ldots,\ell\}}(\lambda'),
\end{equation}
which is known as the (dual) RSK correspondence.

Suppose that $\mc{A}$ has a minimal element, that is, $\mc{A}=\{\,a_1<a_2<\ldots\,\}$. For $a\in \mc{A}$,  
we identify ${\bf m}_{(a)}$ with a tableau in $SST_{\{1,\ldots,\ell\}}(1^{p})$ (resp.  $SST_{\{1,\ldots,\ell\}}(p)$) if $a\in \mc{A}_0$ (resp. $a\in \mc{A}_1$), where $p=|{\bf m}_{(a)}|$. Hence ${\bf m}_{(a)}$ can be regarded as an element of a $\gl_\ell$-crystal \cite{KashNaka} with respect to Kashiwara operators, say $\tE_i$ and $\tF_i$ for $i=1,\ldots, \ell-1$, and ${\bf m}$ as $\cdots\otimes {\bf m}_{(a_2)}\otimes {\bf m}_{(a_1)}$ following \eqref{lower tensor product rule}. Then the map \eqref{RSK} is an isomorphism of $\gl_{\ell}$-crystals by \cite[Theorems 3.11 and 4.5]{K07}. Note that on the righthand side of \eqref{RSK}, $\tE_i$ and $\tF_i$ act on $SST_{\{1,\ldots,\ell\}}(\lambda')$.

Fix $i\in \{1,\ldots, \ell-1\}$. Let $w=w_1\ldots w_r$ be a finite word with letters in $\{\,1,\ldots,\ell\,\}$.
Let $\sigma=(\sigma_{1},\ldots,\sigma_{r})$ be a sequence with $\sigma_{j}\in \{\,+\,,\,-\, , \ \cdot\ \}$ such that $\sigma_j=+$ if $w_j=i$, $-$ if $w_j=i+1$, and $\cdot$  otherwise. We let $\sigma(w;i)=(a,b)$, where $a$ (resp. $b$) is the number of $-$'s (resp. $+$'s) in $\td{\sigma}$ (see Section \ref{subsec:signature}). 
If we regard $w$ as an element of a $\gl_\ell$-crystal, then $\tE_i w$ (resp. $\tF_iw$) is the word replacing $i+1$ (resp. $i$) corresponding to the right-most $-$ (resp. the left-most $+$) in $\td{\sigma}$ with $i$ (resp. $i+1$).

Now, let $\sigma({\bf m};i)=\sigma(w({\bf m});i)$, where $w({\bf m})= \ldots w({\bf m}_{{(a_2)}}) w({\bf m}_{{(a_1)}})$ is the concatenation of the words $w({\bf m}_{(a_k)})$ ($1\leq k\leq \ell$). Also, for $T\in SST_{\{1,\ldots,\ell\}}(\lambda)$, let $\sigma(T;i)=\sigma(w(T);i)$. Then the action of $\tE_i$ and $\tF_i$ on ${\bf m}$ and $T$ can be described in terms of $\sigma({\bf m};i)$ and $\sigma(T;i)$ as in the above paragraph. Since the bijection \eqref{RSK} commutes with $\tE_i$ and $\tF_i$, we have 
\begin{equation}\label{invariance of signature}
\sigma({\bf m};i)=\sigma(Q({\bf m});i).
\end{equation}

Finally, let $U$ and $V$ be $\mc{A}$-semistandard tableaux of single-columned shapes. 
Let ${\bf m}=[\,{\bf m}^{(2)}:{\bf m}^{(1)}\,]\in {\bf M}_{\mc{A}\times 2}$, where ${\bf m}^{(2)}$ (resp. ${\bf m}^{(1)}$) corresponds to $U$ (resp. $V$). Then it is not difficult to see that 
\begin{equation}\label{signature in matrix}
\sigma(U,V)= \sigma({\bf m};1)=(\max\{\,k\,|\,\tE_1^k{\bf m}\neq 0\,\},\max\{\,k\,|\,\tF_1^k{\bf m}\neq 0\,\}).
\end{equation}

\subsection{Ortho-symplectic tableaux of type $D$}

\begin{df}\label{Definition of T_A(a)}{\rm \mbox{}

(1) For $a \in \Z_{\geq  0}$, we define ${\bf T}^{{\mf d}}_\mc{A}(a)={\bf T}_\mc{A}(a)$  to be  
 the set of  $T=(T^{\tt L}, T^{\tt R})\in SST_{\mc{A}}(\lambda(a,b,c))$ such that
\begin{itemize}
\item[(i)] $b,c \in  2\Z_{\geq 0}$,

\item[(ii)] $\sigma(T^{\tt L}, T^{\tt R})=(a-r,b-r)$ for some $r=0,1$.
\end{itemize}
We denote $r$ in (ii) by ${\mf r}_{T}$. We also define $\ov{\bf T}_{\mc{A}}(0)$ to be set of  $T\in SST_{\mc{A}}(\lambda(0,b,c+1))$ for some $b,c \in 2\Z_{\geq 0}$.

(2) Let ${\bf T}^{\rm sp}_\mc{A}$ be the set of $\mc{A}$-semistandard tableaux of single-columned shape. We define ${\bf T}^{\rm sp\,+}_\mc{A}=\{\,T\in {\bf T}^{\rm sp}_\mc{A}\,|\,{\mf r}_{T}=0\,\}$ and ${\bf T}^{\rm sp\,-}_\mc{A}=\{\,T\in {\bf T}^{\rm sp}_\mc{A}\,|\,{\mf r}_{T}=1\,\}$, where ${\mf r}_{T}$ is defined to be the residue of ${\rm ht}(T)$ modulo 2.
}
\end{df}

\begin{rem}\label{meaning of signature restriction}{\rm Given $T\in SST_{\mc{A}}(\lambda(a,b,c))$, one may regard the pair $(T^{\tt L},T^{\tt R})$ as a (not necessarily $\mc{A}$-semistandard) tableau of shape $\lambda(a-k,b-k,c+k)$ sliding $T^{\tt R}$ by $k$ positions down for $0\leq k\leq \min\{a,b\}$.
Then by \eqref{signature condition}, Definition \ref{Definition of T_A(a)}(ii) means that the pair $(T^{\tt L},T^{\tt R})$ is $\mc{A}$-semistandard of shape $\lambda(a-k,b-k,c+k)$ if and only if $k$ is either $0$ or $1$, and the maximum of such $k$ is ${\mf r}_T$.

}
\end{rem}

\begin{ex}\label{ex:T(a)}{\rm
Suppose that $\mc{A}=\J_{4|\infty}$, and let $T\in {\bf T}_\mc{A}(3)$ be as follows.
$$\resizebox{.9\hsize}{!}{$
T=(T^{\tt L}, T^{\tt R}) \ = \ 
{\def\lr#1{\multicolumn{1}{|@{\hspace{.75ex}}c@{\hspace{.75ex}}|}{\raisebox{-.04ex}{$#1$}}}\raisebox{-.6ex}
{$\begin{array}{cc}
\cline{2-2}
&\lr{\overline{3}}\\ 
 \cline{2-2}
 &\lr{\overline{2}}\\ 
\cline{1-1}\cline{2-2}
\lr{\overline{4}}&\lr{\tfrac{3}{2}}\\ 
\cline{1-1}\cline{2-2}
\lr{\overline{1}}&\lr{\tfrac{5}{2}}\\ 
\cline{1-1}\cline{2-2} 
\lr{\tfrac{1}{2}} \\ 
\cline{1-1}
\lr{\tfrac{3}{2}} \\ 
\cline{1-1}
\lr{\tfrac{3}{2}} \\ 
\cline{1-1}
\end{array}$}}\ \in SST_{\J_{4|\infty}}(\lambda(3,2,2)) \ \ \ \ \ \  \stackrel{\rightsquigarrow}{\text{\tiny sliding  $T^{\tt R}$ down}} \ \ \ \ \ \ \
{\def\lr#1{\multicolumn{1}{|@{\hspace{.75ex}}c@{\hspace{.75ex}}|}{\raisebox{-.04ex}{$#1$}}}\raisebox{-.6ex}
{$\begin{array}{cc}
\cline{2-2}
&\lr{\overline{3}}\\ 
 \cline{2-2}
\cline{1-1}\cline{2-2}
\lr{\overline{4}}&\lr{\overline{2}}\\ 
\cline{1-1}\cline{2-2}
\lr{\overline{1}}&\lr{\tfrac{3}{2}}\\ 
\cline{1-1}\cline{2-2} 
\lr{\tfrac{1}{2}} & \lr{\tfrac{5}{2}}\\ 
\cline{1-1}\cline{2-2}
\lr{\tfrac{3}{2}} \\ 
\cline{1-1}
\lr{\tfrac{3}{2}} \\ 
\cline{1-1}
\end{array}$}} \ \ \in SST_{\J_{4|\infty}}(\lambda(2,1,3))$}
$$ 
\noindent We have ${\mf r}_T=1$ since we also have a $\J_{4|\infty}$-semistandard tableau of shape $\lambda(2,1,3)$ after sliding $T^{\tt R}$ down by one position (the tableau on the right).
}
\end{ex}

For $T\in {\bf T}_{\mc{A}}(a)$,  let us identify $T$ with ${\bf m}\in {\bf M}_{\mc{A}\times 2}$ such that $T^{\tt L}$ (resp. $T^{\tt R}$) correspond to ${\bf m}^{(2)}$ (resp. ${\bf m}^{(1)}$). Then we define  
\begin{equation}
({}^{\tt L}T, {}^{\tt R}T)=\tE_1^{a-{\mf r}_T}T,
\end{equation}
that is, the pair of tableaux corresponding to the matrix $\tE_1^{a-{\mf r}_T}{\bf m}$, and 
\begin{equation}
({T}^{\tt L^\ast}, {T}^{\tt R^\ast})=\tF_1 T,
\end{equation}
when ${\mf r}_T=1$.
Note that 
\begin{equation}
\begin{split}
{\rm ht}({}^{\tt L}T)&={\rm ht}({T}^{\tt L})-a+{\mf r}_T,\ \
{\rm ht}({}^{\tt R}T)={\rm ht}(T^{\tt R})+a-{\mf r}_T.
\end{split}
\end{equation}

One can describe algorithms for $({}^{\tt L}T,{}^{\tt R}T)$ and $({T}^{\tt L^\ast},{T}^{\tt R^\ast})$ more explicitly as follows.  \vskip 2mm

\begin{alg}\label{algorithm-1}{\rm \mbox{}

\begin{itemize}
\item[(1)] Let $\boxed{y}$ be the box  at the bottom of $T^{\tt R}$.

\item[(2)] Slide down $\boxed{y}$ until the entry $x$ of $T^{\tt L}$ in the same row is no greater (resp. smaller) than $y$ if $y$ is even (resp. odd). If $y$ is even (resp. odd) and no entry of $T^{\tt L}$ is greater than (resp. greater than or equal to) $y$, we place $\boxed{y}$ to the right of the bottom entry of $T^{\tt L}$.

\item[(3)] Repeat the process (2) with the entries of $T^{\tt R}$ above $y$ until there is no moving down of the entries in $T^{\tt R}$. 

\item[(4)] Move each box $\boxed{x}$ in $T^{\tt L}$ to the right if its right position is empty. (Indeed the number of such boxes is $a-{\mf r}_T$.)

\item[(5)] Then ${}^{\tt R}T$ is the tableau given by the boxes in $T^{\tt R}$ together with boxes which have come from the left, and ${}^{\tt L}T$ is the tableau given by the remaining boxes on the left.  
\end{itemize}
 
In case of Example \ref{ex:T(a)}, we have
$$\resizebox{.8\hsize}{!}{$
T=(T^{\tt L}, T^{\tt R})\ = \ 
{\def\lr#1{\multicolumn{1}{|@{\hspace{.75ex}}c@{\hspace{.75ex}}|}{\raisebox{-.04ex}{$#1$}}}\raisebox{-.6ex}
{$\begin{array}{cc}
\cline{2-2}
&\lr{\overline{3}}\\ 
 \cline{2-2}
 &\lr{\overline{2}}\\ 
\cline{1-1}\cline{2-2}
\lr{\overline{4}}&\lr{\tfrac{3}{2}}\\ 
\cline{1-1}\cline{2-2}
\lr{\overline{1}}&\lr{\tfrac{5}{2}}\\ 
\cline{1-1}\cline{2-2} 
\lr{\tfrac{1}{2}} \\ 
\cline{1-1}
\lr{\tfrac{3}{2}} \\ 
\cline{1-1}
\lr{\tfrac{3}{2}} \\ 
\cline{1-1}
\end{array}$}}
\ \ \  \rightarrow  \ \ \ 
{\def\lr#1{\multicolumn{1}{|@{\hspace{.75ex}}c@{\hspace{.75ex}}|}{\raisebox{-.04ex}{$#1$}}}\raisebox{-.6ex}
{$\begin{array}{cc}
& \\
\cline{2-2}
&\lr{\overline{3}}\\ 
\cline{1-1}\cline{2-2}
\lr{\overline{4}}& \lr{\overline{2}}\\ 
\cline{1-1}\cline{2-2}
\lr{\overline{1}}\\ 
\cline{1-1}\cline{2-2} 
\lr{\tfrac{1}{2}}&\lr{\tfrac{3}{2}} \\ 
\cline{1-1}\cline{2-2}
\lr{\tfrac{3}{2}} \\ 
\cline{1-1}\cline{2-2}
\lr{\tfrac{3}{2}}&\lr{\tfrac{5}{2}} \\ 
\cline{1-1}\cline{2-2}
\end{array}$}}
\ \ \  \rightarrow  \ \ \ 
{\def\lr#1{\multicolumn{1}{|@{\hspace{.75ex}}c@{\hspace{.75ex}}|}{\raisebox{-.04ex}{$#1$}}}\raisebox{-.6ex}
{$\begin{array}{cc}
& \\
\cline{2-2}
&\lr{\overline{3}}\\ 
\cline{1-1}\cline{2-2}
\lr{\overline{4}}& \lr{\overline{2}}\\ 
\cline{1-1}\cline{2-2}
& \lr{\overline{1}}\\ 
\cline{1-1}\cline{2-2} 
\lr{\tfrac{1}{2}}&\lr{\tfrac{3}{2}} \\ 
\cline{1-1}\cline{2-2}
&\lr{\tfrac{3}{2}} \\ 
\cline{1-1}\cline{2-2}
\lr{\tfrac{3}{2}}&\lr{\tfrac{5}{2}} \\ 
\cline{1-1}\cline{2-2}
\end{array}$}}
\ \ \  \rightarrow  \ \ \ 
{\def\lr#1{\multicolumn{1}{|@{\hspace{.75ex}}c@{\hspace{.75ex}}|}{\raisebox{-.04ex}{$#1$}}}\raisebox{-.6ex}
{$\begin{array}{cc}
& \\ 
\cline{1-1}\cline{2-2}
\lr{\overline{4}}& \lr{\overline{3}}\\ 
\cline{1-1}\cline{2-2}
\lr{\tfrac{1}{2}}& \lr{\overline{2}}\\ 
\cline{1-1}\cline{2-2} 
\lr{\tfrac{3}{2}}&\lr{\overline{1}} \\ 
\cline{1-1}\cline{2-2}
&\lr{\tfrac{3}{2}} \\ 
\cline{2-2}
&\lr{\tfrac{3}{2}} \\ 
\cline{2-2}
&\lr{\tfrac{5}{2}} \\ 
\cline{2-2}
\end{array}$}}\ = \ ({}^{\tt L}T, {}^{\tt R}T)$}
$$\vskip 2mm
}
\end{alg} 
\noindent where we arrange $({}^{\tt L}T, {}^{\tt R}T)$ so that each of the pairs $({}^{\tt L}T,T^{\tt R})$ and $(T^{\tt L}, {}^{\tt R}T)$ shares the same bottom line.

\begin{alg}\label{algorithm-2}{\rm \mbox{}

\begin{itemize}
\item[(1)] Let $\boxed{x}$ be the box at the top of $T^{\tt L}$.

\item[(2)] Slide upward $\boxed{x}$ until the entry $y$ of $T^{\tt R}$ in the same row is no smaller than (resp. no greater or equal to) $x$ if $x$ is even (resp. odd). If $x$ is even (resp. odd)  and no entry of $T^{\tt R}$ is smaller than (resp. greater than or equal to) $x$, we place $\boxed{x}$ to the right of the top entry of $T^{\tt R}$.

\item[(3)] Repeat the process (2) with the next entry of $T^{\tt L}$ below $x$ until there is no moving up of the entries in $T^{\tt L}$. 

\item[(4)] Choose the lowest box $\boxed{y}$ in $T^{\tt R}$ whose left position is empty, and then move it to the left.  (Since ${\mf r}_T=1$, there exists at least one such $\boxed{y}$.)

\item[(5)] Then ${T}^{\tt L^\ast}$ is the tableau given by the boxes in $T^{\tt L}$ together with $\boxed{y}$, and ${T}^{\tt R^\ast}$ is the tableau given by the remaining boxes on the right.  
\end{itemize}

In case of Example \ref{ex:T(a)}, we have
$$\resizebox{.8\hsize}{!}{$
T=(T^{\tt L}, T^{\tt R})\ = \ 
{\def\lr#1{\multicolumn{1}{|@{\hspace{.75ex}}c@{\hspace{.75ex}}|}{\raisebox{-.04ex}{$#1$}}}\raisebox{-.6ex}
{$\begin{array}{cc}
\cline{2-2}
&\lr{\overline{3}}\\ 
 \cline{2-2}
 &\lr{\overline{2}}\\ 
\cline{1-1}\cline{2-2}
\lr{\overline{4}}&\lr{\tfrac{3}{2}}\\ 
\cline{1-1}\cline{2-2}
\lr{\overline{1}}&\lr{\tfrac{5}{2}}\\ 
\cline{1-1}\cline{2-2} 
\lr{\tfrac{1}{2}} \\ 
\cline{1-1}
\lr{\tfrac{3}{2}} \\ 
\cline{1-1}
\lr{\tfrac{3}{2}} \\ 
\cline{1-1}
\end{array}$}}
\ \ \  \rightarrow  \ \ \ 
{\def\lr#1{\multicolumn{1}{|@{\hspace{.75ex}}c@{\hspace{.75ex}}|}{\raisebox{-.04ex}{$#1$}}}\raisebox{-.6ex}
{$\begin{array}{cc}
\cline{1-1}\cline{2-2}
\lr{\overline{4}}&\lr{\overline{3}}\\ 
\cline{1-1}\cline{2-2}
& \lr{\overline{2}}\\ 
\cline{1-1}\cline{2-2}
\lr{\overline{1}} &\lr{\tfrac{3}{2}}\\ 
\cline{1-1}\cline{2-2} 
\lr{\tfrac{1}{2}}&\lr{\tfrac{5}{2}} \\ 
\cline{1-1}\cline{2-2}
\lr{\tfrac{3}{2}} \\ 
\cline{1-1} 
\lr{\tfrac{3}{2}}  \\ 
\cline{1-1} \\
\end{array}$}}
\ \ \  \rightarrow  \ \ \ 
{\def\lr#1{\multicolumn{1}{|@{\hspace{.75ex}}c@{\hspace{.75ex}}|}{\raisebox{-.04ex}{$#1$}}}\raisebox{-.6ex}
{$\begin{array}{cc}
\cline{1-1}\cline{2-2}
\lr{\overline{4}}&\lr{\overline{3}}\\ 
\cline{1-1}\cline{2-2}
\lr{\overline{2}} &\\ 
\cline{1-1}\cline{2-2}
\lr{\overline{1}} &\lr{\tfrac{3}{2}}\\ 
\cline{1-1}\cline{2-2} 
\lr{\tfrac{1}{2}}&\lr{\tfrac{5}{2}} \\ 
\cline{1-1}\cline{2-2}
\lr{\tfrac{3}{2}} \\ 
\cline{1-1} 
\lr{\tfrac{3}{2}}  \\ 
\cline{1-1} \\
\end{array}$}}
\ \ \  \rightarrow  \ \ \ 
{\def\lr#1{\multicolumn{1}{|@{\hspace{.75ex}}c@{\hspace{.75ex}}|}{\raisebox{-.04ex}{$#1$}}}\raisebox{-.6ex}
{$\begin{array}{cc}
& \\
\cline{1-1}\cline{2-2}
\lr{\overline{4}}&\lr{\overline{3}}\\ 
\cline{1-1}\cline{2-2}
\lr{\overline{2}} &\lr{\tfrac{3}{2}}\\ 
\cline{1-1}\cline{2-2}
\lr{\overline{1}} &\lr{\tfrac{5}{2}}\\ 
\cline{1-1}\cline{2-2} 
\lr{\tfrac{1}{2}}& \\ 
\cline{1-1} 
\lr{\tfrac{3}{2}} \\ 
\cline{1-1} 
\lr{\tfrac{3}{2}}  \\ 
\cline{1-1} 
\end{array}$}}\ = \ ({T}^{\tt L^\ast} , {T}^{\tt R^\ast})$}
$$ \vskip 2mm
\noindent Note that each of the pairs $(T^{\tt L^*},T^{\tt L})$ and $(T^{\tt R}, T^{\tt R^*})$ shares the same bottom line.
}
\end{alg}

\begin{df}\label{admissible}{\rm  \mbox{}

(1) For $T \in {\bf T}_\mc{A}(a)$ and $S \in {\bf T}_\mc{A}(a')\cup {\bf T}_{\mc{A}}^{\rm sp}$ with $a\geq a'$,
we write $T\prec S$ if 

\begin{itemize}
\item[(i)]  ${\rm ht}({T}^{\tt R})\leq {\rm ht}(S^{\tt L}) -a'+  2{\mf r}_{S}{\mf r}_{T}$,

\item[(ii)] for $i\geq 1$, we have
\begin{equation*}
\begin{cases}
{T}^{\tt R^\ast}(i)\leq  {}^{\tt L}S(i), & \text{if ${\mf r}_S={\mf r}_T=1$},\\
{T}^{\tt R}(i)\leq  {}^{\tt L}S(i), & \text{otherwise},
\end{cases}
\end{equation*}

\item[(iii)] for $i\geq 1$, we have 
\begin{equation*}
\begin{cases}
{}^{\tt R}T(i+a-a'+\epsilon)\leq {S}^{\tt L^\ast}(i), & \text{if ${\mf r}_S={\mf r}_T=1$},\\
{}^{\tt R}T(i+a-a')\leq {S}^{\tt L}(i), & \text{otherwise},
\end{cases}
\end{equation*}
\end{itemize}
where the equality holds in (ii) and (iii) only if the entries are even, and $\epsilon=1$ if $S\in {\bf T}_{\mc{A}}^{\rm sp -}$   and $0$ otherwise.
Here we assume that $a'={\mf r}_S$, $S=S^{\tt L}={}^{\tt L}S={S}^{\tt L^\ast}$ when $S\in {\bf T}_{\mc{A}}^{\rm sp}$.\vskip 1mm

(2)
For  $T\in {\bf T}_{\mc{A}}(a)$ and $S \in \ov{\bf T}_{\mc{A}}(0)$, define $T\prec S$ if $T\prec S^{\tt L}$ in the sense of (1), where $S^{\tt L}\in {\bf T}^{\rm sp -}_{\mc{A}}$.
\vskip 1mm

(3) 
For $T\in \ov{\bf T}_{\mc{A}}(0)$ and $S \in \ov{\bf T}_{\mc{A}}(0)\cup {\bf T}^{\rm sp -}_{\mc{A}}$, define $T\prec S$ if $(T^{\tt R},S^{\tt L})\in \ov{\bf T}_{\mc{A}}(0)$.\vskip 2mm

We say that the pair $(T,S)$ is {\em admissible} when $T\prec S$.
}\end{df}

\begin{ex}{\rm  
Consider the pair $(T,S)$  
$$
\resizebox{.3\hsize}{!}{$T\ = \
 {\def\lr#1{\multicolumn{1}{|@{\hspace{.75ex}}c@{\hspace{.75ex}}|}{\raisebox{-.04ex}{$#1$}}}\raisebox{-.6ex}
{$\begin{array}{cc}
\cline{2-2}
&\lr{\overline{3}}\\ 
 \cline{2-2}
 &\lr{\overline{2}}\\ 
\cline{1-1}\cline{2-2}
\lr{\overline{4}}&\lr{\tfrac{3}{2}}\\ 
\cline{1-1}\cline{2-2}
\lr{\overline{1}}&\lr{\tfrac{5}{2}}\\ 
\cline{1-1}\cline{2-2} 
\lr{\tfrac{1}{2}} \\ 
\cline{1-1}
\lr{\tfrac{3}{2}} \\ 
\cline{1-1}
\lr{\tfrac{3}{2}} \\ 
\cline{1-1}
\end{array}$}} \  \  \ \ \ \
{\def\lr#1{\multicolumn{1}{|@{\hspace{.75ex}}c@{\hspace{.75ex}}|}{\raisebox{-.04ex}{$#1$}}}\raisebox{-.6ex}
{$\begin{array}{cc}
\cline{2-2}
&\lr{\overline{2}}\\ 
 \cline{2-2}
 &\lr{\overline{1}}\\ 
\cline{1-1}\cline{2-2}
\lr{\overline{1}}&\lr{\tfrac{7}{2}}\\ 
\cline{1-1}\cline{2-2}
\lr{\tfrac{5}{2}}&\lr{\tfrac{9}{2}}\\ 
\cline{1-1}\cline{2-2} 
\lr{\tfrac{7}{2}} &\\ 
\cline{1-1}
\lr{\tfrac{9}{2}} &\\ 
\cline{1-1}\\
\end{array}$}}\ = \ S$}
$$
\noindent where $T$ is as in Example \ref{ex:T(a)}. Note that $T$ and $S$ are arranged so that $T^{\tt R}$ and $S^{\tt R}$ share the same bottom line. First, we have $4={\rm ht}(T^{\tt R})\leq {\rm ht}(S^{\tt L})-2+2{\mf r}_T{\mf r}_S=4$, which satisfies Definition \ref{admissible}(1)(i). Since
$$\resizebox{.45\hsize}{!}{$
 ({}^{\tt L}S, {}^{\tt R}S) \ = \ 
{\def\lr#1{\multicolumn{1}{|@{\hspace{.75ex}}c@{\hspace{.75ex}}|}{\raisebox{-.04ex}{$#1$}}}\raisebox{-.6ex}
{$\begin{array}{cc}
\cline{1-1}\cline{2-2}
\lr{\overline{1}}&\lr{\overline{2}}\\ 
\cline{1-1}\cline{2-2}
\lr{\tfrac{5}{2}}& \lr{\overline{1}}\\ 
\cline{1-1}\cline{2-2}
\lr{\tfrac{7}{2}}& \lr{\tfrac{7}{2} }\\ 
\cline{1-1}\cline{2-2} 
 &\lr{\tfrac{9}{2}} \\ 
 \cline{2-2}
&\lr{\tfrac{9}{2}} \\ 
\cline{2-2}
\end{array}$}}\ \ \  \
({S}^{\tt L^\ast} , {S}^{\tt R^\ast}) \ = \ 
{\def\lr#1{\multicolumn{1}{|@{\hspace{.75ex}}c@{\hspace{.75ex}}|}{\raisebox{-.04ex}{$#1$}}}\raisebox{-.6ex}
{$\begin{array}{cc}
\cline{1-1}\cline{2-2}
\lr{\overline{2}}&\lr{\overline{1}}\\ 
\cline{1-1}\cline{2-2}
\lr{\overline{1}} &\lr{\tfrac{7}{2}}\\ 
\cline{1-1}\cline{2-2}
\lr{\tfrac{5}{2}} &\lr{\tfrac{9}{2}}\\ 
\cline{1-1}\cline{2-2} 
\lr{\tfrac{7}{2}}& \\ 
\cline{1-1} 
\lr{\tfrac{9}{2}} \\ 
\cline{1-1} 
\end{array}$}}$}
$$ 
we have $$\resizebox{.45\hsize}{!}{$
({T}^{\tt R^\ast},{}^{\tt L}S)\ = \
{\def\lr#1{\multicolumn{1}{|@{\hspace{.75ex}}c@{\hspace{.75ex}}|}{\raisebox{-.04ex}{$#1$}}}\raisebox{-.6ex}
{$\begin{array}{cc}
\cline{1-1}\cline{2-2}
\lr{\overline{3}} & \lr{\overline{1}}\\ 
\cline{1-1}\cline{2-2}
\lr{\tfrac{3}{2}}& \lr{\tfrac{5}{2}} \\ 
\cline{1-1}\cline{2-2}
\lr{\tfrac{5}{2}}& \lr{\tfrac{7}{2}} \\ 
\cline{1-1}\cline{2-2} 
\end{array}$}}
\ \ \ \ \ \ \ 
({}^{\tt R}T,S^{\tt L^\ast})\ = \
{\def\lr#1{\multicolumn{1}{|@{\hspace{.75ex}}c@{\hspace{.75ex}}|}{\raisebox{-.04ex}{$#1$}}}\raisebox{-.6ex}
{$\begin{array}{cc}
\cline{1-1}\cline{2-2}
\lr{\overline{3}} & \lr{\overline{2} }\\ 
\cline{1-1}\cline{2-2}
\lr{\overline{2}}& \lr{\overline{1}} \\ 
\cline{1-1}\cline{2-2}
\lr{\overline{1}}& \lr{\tfrac{5}{2}} \\ 
\cline{1-1}\cline{2-2} 
\lr{\tfrac{3}{2}}& \lr{\tfrac{7}{2}} \\ 
\cline{1-1}\cline{2-2}
\lr{\tfrac{3}{2}}& \lr{\tfrac{9}{2}}\\ 
\cline{1-1}\cline{2-2}
\lr{\tfrac{5}{2}}& \\ 
\cline{1-1}
\end{array}$}}$}
$$  
\noindent which are $\J_{4|\infty}$-semistandard, and hence $T\prec S$ by Definition \ref{admissible}(1)(ii) and (iii).
On the other hand, if we have  
$$
\resizebox{.3\hsize}{!}{$T\ = \
 {\def\lr#1{\multicolumn{1}{|@{\hspace{.75ex}}c@{\hspace{.75ex}}|}{\raisebox{-.04ex}{$#1$}}}\raisebox{-.6ex}
{$\begin{array}{cc}
& \\
& \\
\cline{2-2}
&\lr{\overline{3}}\\ 
 \cline{2-2}
 &\lr{\overline{2}}\\ 
\cline{1-1}\cline{2-2}
\lr{\overline{4}}&\lr{\tfrac{3}{2}}\\ 
\cline{1-1}\cline{2-2}
\lr{\overline{1}}&\lr{\tfrac{5}{2}}\\ 
\cline{1-1}\cline{2-2} 
\lr{\tfrac{1}{2}} \\ 
\cline{1-1}
\lr{\tfrac{3}{2}} \\ 
\cline{1-1}
\lr{\tfrac{3}{2}} \\ 
\cline{1-1}
\end{array}$}} \  \  \ \ \ \
{\def\lr#1{\multicolumn{1}{|@{\hspace{.75ex}}c@{\hspace{.75ex}}|}{\raisebox{-.04ex}{$#1$}}}\raisebox{-.6ex}
{$\begin{array}{cc}
\cline{1-1}\cline{2-2}
\lr{\overline{3}}&\lr{\overline{2}}\\ 
\cline{1-1}\cline{2-2}
\lr{\overline{2}}&\lr{\overline{1}}\\ 
\cline{1-1}\cline{2-2}
\lr{\overline{1}}&\lr{\tfrac{1}{2}}\\ 
\cline{1-1}\cline{2-2}
\lr{\tfrac{1}{2}} &\lr{\tfrac{3}{2}}\\ 
\cline{1-1}\cline{2-2}
\lr{\tfrac{3}{2}}&\lr{\tfrac{7}{2}}\\ 
\cline{1-1}\cline{2-2}
\lr{\tfrac{5}{2}}&\lr{\tfrac{9}{2}}\\ 
\cline{1-1}\cline{2-2} 
\lr{\tfrac{7}{2}} \\ 
\cline{1-1}\\
& \\
\end{array}$}}\ = \ S$}
$$
with ${\mf r}_S=0$, then
$$\resizebox{.7\hsize}{!}{$
({}^{\tt L}S, {}^{\tt R}S) \ = \ 
{\def\lr#1{\multicolumn{1}{|@{\hspace{.75ex}}c@{\hspace{.75ex}}|}{\raisebox{-.04ex}{$#1$}}}\raisebox{-.6ex}
{$\begin{array}{cc}
\cline{1-1}\cline{2-2}
\lr{\overline{3}}&\lr{\overline{2}}\\ 
\cline{1-1}\cline{2-2}
\lr{\overline{2}}&\lr{\overline{1}}\\ 
\cline{1-1}\cline{2-2}
\lr{\overline{1}}&\lr{\tfrac{1}{2}}\\ 
\cline{1-1}\cline{2-2}
\lr{\tfrac{1}{2}} &\lr{\tfrac{3}{2}}\\ 
\cline{1-1}\cline{2-2}
\lr{\tfrac{5}{2}}&\lr{\tfrac{3}{2}}\\ 
\cline{1-1}\cline{2-2}
\lr{\tfrac{7}{2}}&\lr{\tfrac{7}{2}}\\ 
\cline{1-1}\cline{2-2} 
&\lr{\tfrac{9}{2}} \\ 
\cline{2-2}
\end{array}$}} \ \ \  \ \ \ \ \ \ 
({T}^{\tt R},{}^{\tt L}S)\ = \
{\def\lr#1{\multicolumn{1}{|@{\hspace{.75ex}}c@{\hspace{.75ex}}|}{\raisebox{-.04ex}{$#1$}}}\raisebox{-.6ex}
{$\begin{array}{cc}
\cline{2-2}
& \lr{\overline{3}}\\ 
\cline{2-2}
& \lr{\overline{2}}\\ 
\cline{1-1}\cline{2-2}
\lr{\overline{3}} & \lr{\overline{1}}\\ 
\cline{1-1}\cline{2-2}
\lr{\overline{2}} & \lr{\tfrac{1}{2}}\\ 
\cline{1-1}\cline{2-2}
\lr{\tfrac{3}{2}}& \lr{\tfrac{5}{2}} \\ 
\cline{1-1}\cline{2-2}
\lr{\tfrac{5}{2}}& \lr{\tfrac{7}{2}} \\ 
\cline{1-1}\cline{2-2} \\
\end{array}$}}
\ \ \ \ \ \ \ \ 
({}^{\tt R}T,S^{\tt L})\ = \
{\def\lr#1{\multicolumn{1}{|@{\hspace{.75ex}}c@{\hspace{.75ex}}|}{\raisebox{-.04ex}{$#1$}}}\raisebox{-.6ex}
{$\begin{array}{cc}
\cline{2-2}
 & \lr{\overline{3} }\\
\cline{2-2}
 & \lr{\overline{2} }\\
\cline{2-2}
 & \lr{\overline{1} }\\
\cline{1-1}\cline{2-2}
\lr{\overline{3}} & \lr{\tfrac{1}{2} }\\ 
\cline{1-1}\cline{2-2}
\lr{\overline{2}}& \lr{\tfrac{3}{2}} \\ 
\cline{1-1}\cline{2-2}
\lr{\overline{1}}& \lr{\tfrac{5}{2}} \\ 
\cline{1-1}\cline{2-2} 
\lr{\tfrac{3}{2}}& \lr{\tfrac{7}{2}} \\ 
\cline{1-1}\cline{2-2}
\lr{\tfrac{3}{2}}&  \\ 
\cline{1-1} 
\lr{\tfrac{5}{2}}& \\ 
\cline{1-1}
\end{array}$}}$}
$$  
\noindent Hence we also have $T\prec S$.

}
\end{ex}

\begin{rem}\label{remark on admissible}
{\rm We can describe equivalent conditions for admissibility in terms of signature $\sigma$, which will be useful in the proof of Theorem \ref{Schur positivity}.
Let $(T,S)$ be as in Definition \ref{admissible}(1). 
The condition (ii) is equivalent to saying  that
\begin{equation}\label{admissible-equiv-1}
\begin{split}
&\text{$(T^{\tt R},{}^{\tt L}S)$ or $(T^{\tt R^\ast},{}^{\tt L}S) \in SST_{\mc{A}}(\lambda(0,b,c))$\ \  with}\\
(b,c)&=
({\rm ht}({S}^{\tt L})-{\rm ht}(T^{\tt R})-a'+{\mf r}_S({\mf r}_T+1), {\rm ht}(T^{\tt R})-{\mf r}_T{\mf r}_S), 
\end{split}
\end{equation}
and hence by \eqref{signature condition} equivalent to
\begin{equation}\label{signature condition-1}
\sigma(T^{\tt R},{}^{\tt L}S)\ \text{ or }\ \sigma(T^{\tt R^\ast},{}^{\tt L}S)=(0, b).
\end{equation} 

In a similar way, the condition (iii) is equivalent to saying that 
\begin{equation}\label{admissible-equiv-2}
\begin{split}
&\text{$({}^{\tt R}T,{S}^{\tt L})$ or $({}^{\tt R}T,{S}^{\tt L^\ast})\in SST_{\mc{A}}(\lambda(a-a'+\epsilon,b,c-\epsilon))$ \  with}\\
(b&,c)=({\rm ht}(S^{\tt L})-{\rm ht}({T}^{\tt R})-a'+{\mf r}_T({\mf r}_S+1), {\rm ht}({T}^{\tt R})+a'-{\mf r}_T),
\end{split}
\end{equation}
or equivalent to
\begin{equation}\label{signature condition-2}
\sigma({}^{\tt R}T,{S}^{\tt L}) \text{ or } \sigma({}^{\tt R}T,{S}^{\tt L^\ast})=(a-a'+\epsilon-p, b-p),
\end{equation} 
for some $p\geq 0$. We remark that the condition (i) is equivalent to ${\rm ht}({T}^{\tt R})\leq {\rm ht}(S^{\tt L}) -a'+  {\mf r}_T({\mf r}_S+1)$ or $b\geq 0$ since ${\rm ht}(T^{\tt L})$ and ${\rm ht}(S^{\tt L})-a'$ are even integers.
We have similar conditions as in \eqref{signature condition-1} and \eqref{signature condition-2} for the pairs $(T,S)$ in Definition \ref{admissible}(2) and (3). }
\end{rem}\vskip 2mm

Let $(\lambda,\ell)\in\cP({\mf d})$ be given. Let $q_\pm$ and $r_\pm$ be non-negative integers such that
\begin{equation*}
\begin{cases}
\ell-2\lambda_1=2q_++r_+, & \text{if $\ell-2\lambda_1\geq 0$,}\\
2\lambda_1-\ell=2q_-+r_-, & \text{if $\ell-2\lambda_1\leq 0$,}\\
\end{cases}
\end{equation*}
where $r_\pm=0,1$. Let $\ov{\lambda}=(\ov{\lambda}_i)_{i\geq 1}\in\cP$ be such that $\ov{\lambda}_1=\ell-\lambda_1$ and $\ov{\lambda}_i=\lambda_i$ for $i\geq 2$. 
Let 
\begin{equation}\label{Length of tuples}
\begin{split}
& \nu=\lambda', \ \ \ov{\nu}=(\ov{\lambda})',\\
& M_+=\lambda_1, \ \ M_-=\ov{\lambda}_1=\ell-\lambda_1,\\
& L=M_\pm+q_\pm. 
\end{split}
\end{equation}
Note that $2L+r_\pm=\ell$.
Put
\begin{equation}\label{product form of T}
\widehat{\bf T}_{\mc{A}}(\lambda,\ell)=\begin{cases}
{\bf T}_\mc{A}(\nu_{1})\times\cdots \times {\bf T}_\mc{A}(\nu_{M_+})\times  \left({\bf T}_\mc{A}(0)\right)^{q_+} \times \left({\bf T}^{\rm sp +}_{\mc{A}}\right)^{r_+}, &  \!\!\!\text{if $\ell-2\lambda_1\geq 0$},\\
{\bf T}_\mc{A}(\ov{\nu}_{1})\times\cdots \times {\bf T}_\mc{A}(\ov{\nu}_{M_-})\times  \left(\ov{\bf T}_\mc{A}(0)\right)^{q_-}\times \left({\bf T}^{\rm sp -}_{\mc{A}}\right)^{r_-}, &  \!\!\!\text{if $\ell-2\lambda_1\leq0$},
\end{cases}
\end{equation}
Here, we assume that $\left({\bf T}^{\rm sp \pm}_{\mc{A}}\right)^{r_\pm}$ is empty if $r_{\pm}=0$.

Now we introduce our main combinatorial object.

\begin{df}\label{def:osptableaux}{\rm
For $(\lambda,\ell)\in\cP({\mf d})$, we define ${\bf T}^{\mf{d}}_\mc{A}(\lambda,\ell)={\bf T}_\mc{A}(\lambda,\ell)$ to be the set of ${\bf T}=(T_L,\ldots,T_1,T_0)$  in $\widehat{\bf T}_{\mc{A}}(\lambda,\ell)$
such that $T_{k+1}\prec T_{k}$  for $0\leq k\leq L-1$. We call ${\bf T}\in {\bf T}_\mc{A}(\lambda,\ell)$ an {\it ortho-symplectic tableau of type $D$ and shape $(\lambda,\ell)$}.
}
\end{df}

\begin{rem}\label{rem:remark on BC type}{\rm
Here, we are using a convention slightly different from the cases of type $B$ and $C$ in \cite{K13}, when we define the notion of  admissibility and ortho-symplectic tableaux of type $D$. But we may still apply Definition \ref{admissible} to ${\bf T}^{\mf g}_{\mc A}(a)$  and ${\bf T}^{\rm sp}_{\mc{A}}$ for $\g=\mf{b}, \mf{c}$ and $a\geq 0$ in \cite{K13}, where all tableaux are of residue $0$, and define ${\bf T}^{\g}_{\mc{A}}(\lambda,\ell)$ for $(\lambda,\ell)\in\cP(\g)$ as in Definition \ref{def:osptableaux} with the order of product of ${\bf T}^{\mf g}_{\mc A}(a)$  and ${\bf T}^{\rm sp}_{\mc{A}}$'s in \cite[Definition 6.10]{K13} reversed as in \eqref{product form of T}. Then we can check without difficulty that all the results in \cite{K13} can be obtained with this version of ortho-symplectic tableaux of type $B$ and $C$.
 
}
\end{rem}

Let $x_{\mc{A}}=\{\,x_a\,|\,a\in\mc{A}\,\}$ be the set of formal commuting variables indexed by $\mc{A}$. For $\lambda\in\cP$, let $s_{\lambda}(x_{\mc{A}})=\sum_{T}x_\mc{A}^T$ be the super Schur function corresponding to $\lambda$, where the sum is over $T\in SST_{\mc{A}}(\lambda)$ and $x_\mc{A}^T=\prod_{a}x_a^{m_a}$ with ${\rm wt}(T)=(m_a)_{a\in \mc{A}}$.
For $(\lambda,\ell)\in\cP({\mf d})$, put 
\begin{equation*}
S_{(\lambda,\ell)}(x_{\mc{A}})=z^\ell\sum_{{\bf T}\in {\bf T}_\mc{A}(\lambda,\ell)}\prod_{k=0}^L x_\mc{A}^{T_k},
\end{equation*}
where $z$ is another formal variable.
First, we have the following Schur positivity of $S_{(\lambda,\ell)}(x_{\mc{A}})$ as in the case of type $B$ and $C$ \cite[Theorem 6.12]{K13}. 

\begin{thm}\label{Schur positivity}
 For $(\lambda,\ell)\in\cP({\mf d})$, we have
\begin{equation*}
S_{(\lambda,\ell)}(x_{\mc{A}})=z^{\ell}\sum_{\mu\in\cP}K_{\mu\,(\lambda,\ell)}s_\mu(x_\mc{A}),
\end{equation*}
for some non-negative integers $K_{\mu\,(\lambda,\ell)}$. Moreover, the coefficients $K_{\mu\,(\lambda,\ell)}$ do not depend on $\mc{A}$.
\end{thm}
\pf 
Let $L$ be  as in \eqref{Length of tuples}.  
Let ${\bf T}=(T_L,\ldots,T_1,T_0)\in {\bf T}_\mc{A}(\lambda,\ell)$ be given. 
Let ${\bf m}=[\,{\bf m}^{(\ell)}:\cdots :{\bf m}^{(1)}\,]$ be the unique matrix in ${\bf M}_{\mc{A}\times \ell}$, where ${\bf m}^{(1)}$ corresponds to $T_0$, and $[\,{\bf m}^{(2k+1)}:{\bf m}^{(2k)}\,]$ corresponds to $T_{k}$ for $1\leq k\leq L$. We assume that $T_0$ is empty and ${\bf m}^{(1)}$ is trivial when $r_\pm=0$.

Put $Q=Q({\bf m})$, which is of $\{1,\ldots,\ell\}$-semistandard and ${\rm wt}(Q)=(m_1,m_2,\ldots, m_{\ell})$ with $m_i=|{\bf m}^{(i)}|$.
For convenience, we put for $1\leq k\leq M_\pm$,
\begin{equation}\label{notations for lengths}
\begin{split}
&m^{\tt L}_k =m_{2q_\pm +2k+1}, \ \ \ 
m^{\tt R}_k=m_{2q_\pm+2k}, \  \ \  r_k={\mf r}_{T_{q_\pm +k}},\\
&a_k= 
\begin{cases}
\nu_{M_++1-k} & \text{if $\ell-2\lambda_1\geq 0$}, \\
\ov{\nu}_{M_-+1-k} & \text{if $\ell-2\lambda_1\leq 0$}. 
\end{cases}
\end{split}
\end{equation}

First, for $1\leq k\leq L$,  we see from $T_{q_\pm+k}\in {\bf T}_{\mc{A}}(a_k)$, \eqref{invariance of signature}, and \eqref{signature in matrix} that
\begin{itemize}
\item[(Q1)] $m^{\tt L}_{k}-a_{k},  m^{\tt R}_{k} \in 2\Z_{\geq 0}$,

\item[(Q2)] $m^{\tt L}_{k}-a_{k}\leq m^{\tt R}_{k}$,

\item[(Q3)] $\sigma(Q;2q_\pm+2k)=(a_k-r_k,m^{\tt R}_{k}- m^{\tt L}_{k}+a_{k}-r_k)$.

\end{itemize}
Put 
\begin{equation*}
 Q^{(k)}=\tF^{r_{k}r_{k+1}}_{2q_\pm +2k+2}\tE_{2q_\pm +2k}^{a_{k}-r_{k}}Q, \ \ \ \ Q^{[k]}=\tF^{r_kr_{k+1}}_{2q_\pm +2k}\tE_{2q_\pm + 2k+2}^{a_{k+1}-r_{k+1}}Q,
\end{equation*}
for $1 \leq k\leq M_\pm-1$.
Since $T_{q_\pm+k+1}\prec T_{q_\pm+k}$ for $1 \leq k\leq M_\pm-1$, we have by Definition \ref{admissible}(1), \eqref{signature condition-1}, and \eqref{signature condition-2} that
\begin{itemize}
\item[(Q4)] $m^{\tt R}_{k+1}\leq m^{\tt L}_{k} -a_{k}+ 2 r_kr_{k+1}$, 

\item[(Q5)] ${\sigma}(Q^{(k)};2q_\pm+2k)=(0,m^{\tt L}_{k}-m^{\tt R}_{k+1}-a_k+r_{k}(r_{k+1}+1))$,

\item[(Q6)]$\sigma(Q^{[k]};2q_\pm+2k)= (a_{k+1}-a_{k}-p_k, m^{\tt L}_{k}-m^{\tt R}_{k+1}-a_{k}+r_{k+1}(r_k+1)-p_k)$ 
for some $p_k\geq 0$.
\end{itemize} 

Next, since $T_{k}\in {\bf T}_{\mc{A}}(0)$ or $\ov{\bf T}_{\mc{A}}(0)$ for $1\leq k\leq q_\pm$, $T_0\in {\bf T}_{\mc{A}}^{\rm sp \pm}$, and $T_{k+1}\prec T_{k}$ for $0\leq k\leq q_\pm-1$, we have by Definition \ref{admissible}(3) that
\begin{itemize}
\item[(Q7)] $m_k\in \Z_{\geq 0}$ for $0\leq k\leq 2q_+$  and $m_k\in \Z_{> 0}$ for $0\leq k\leq 2q_-$, 

\item[(Q8)] $m_{k+1}\leq m_{k}$ for $0\leq k\leq 2q_\pm-1$,

\item[(Q9)] $\sigma(Q;k)=(0,m_{k}-m_{k+1})$ for $0\leq k\leq 2q_\pm-1$.

\end{itemize}

Finally, put $Q^{[0]}=\tE_{2q_\pm + 2}^{a_{1}-r_{1}}Q$.
Since $T_{q_\pm+1}\prec T_{q_\pm}$, we have by Definition \ref{admissible}(1) and (2), \eqref{signature condition-1}, and \eqref{signature condition-2} that
\begin{itemize}
\item[(Q10)]  $m^{\tt R}_{1} \leq m_{2q_\pm+1}-r_0 + 2r_0r_1$,

\item[(Q11)]  $\sigma(Q;2q_\pm+1)=(0,m_{2q_\pm+1}-m^{\tt R}_{1} +r_{0}r_{1})$,

\item[(Q12)] $\sigma(Q^{[0]};2q_\pm+1)=(a_{1}-p_0, m_{2q_\pm+1}-m^{\tt R}_{1}-r_0+r_1(r_{0}+1)-p_0)$ for some $p_0\geq 0$,
\end{itemize} 
where $r_0=1$ if $\ell-2\lambda_1< 0$, and $0$ otherwise

Conversely, for $\mu\in \cP$, let $(P,Q)$ be given where $P\in SST_{\mc{A}}(\mu)$ and $Q\in SST_{\{1,\ldots,\ell\}}(\mu')$ with ${\rm wt}(Q)=(m_1,\ldots, m_{\ell})$ satisfying the conditions (Q1)--(Q12) for some $r_k$ ($1\leq k\leq M_\pm$) and $p_k$ ($0\leq k\leq M_\pm-1$). Note that if such $Q$ exists, then $r_k$ and $p_k$ are uniquely determined by (Q3), (Q6), and (Q12). By \eqref{RSK}, there exists a unique ${\bf m}\in {\bf M}_{\mc{A}\times \ell}$ such that $(P({\bf m}),Q({\bf m}))=(P,Q)$. Then it follows from  \eqref{signature condition}, \eqref{invariance of signature}, \eqref{signature in matrix}, and Remark \ref {remark on admissible} that there exists a unique ${\bf T}\in {\bf T}_{\mc{A}}(\lambda,\ell)$ which corresponds to ${\bf m}$. 
Hence, the map \eqref{RSK} induces a weight preserving bijection
\begin{equation}\label{Pieri}
{\bf T}_\mc{A}(\lambda,\ell) \longrightarrow \bigsqcup_{\mu\in\cP}  SST_{\mc{A}}(\mu)\times  {\bf K}_{\mu\,(\lambda,\ell)},
\end{equation}
where ${\bf K}_{\mu\,(\lambda,\ell)}$ is the set of $Q\in SST_{\{1,\ldots,\ell\}}(\mu')$ with ${\rm wt}(Q)=(m_1,\ldots, m_{\ell})$ satisfying (Q1)--(Q12).
This implies that $S_{(\lambda,\ell)}(x_{\mc{A}})=z^{\ell}\sum_{\mu\in\cP}K_{\mu\,(\lambda,\ell)}s_\mu(x_\mc{A})$, where $K_{\mu\,(\lambda,\ell)}=|{\bf K}_{\mu\,(\lambda,\ell)}|$.
\qed

\section{Character formula of a highest weight  module}\label{crystal structure for ortho-symplectic tableaux m+n}

\subsection{Lie algebra $\mf{d}_{m+n}$}
We assume the following notations for the classical Lie algebra $\mf{d}_{m+n}$ of type $D_{m+n}$ (see \cite{K13} for more details):

\begin{itemize}
\item[$\cdot$] $\J_{m+n}=\{\,\ov{m}< \ldots <\ov{2}< \ov{1} < 1< 2 <\ldots<n\,\}$,

\item[$\cdot$] $P_{m+n}=\bigoplus_{a\in \J_{m+n}}\Z \delta_a \oplus \Z \Lambda_{\ov{m}}$ : the weight lattice, 


\item[$\cdot$] $I_{m+n}=\{\,\ov{m},\ldots,\ov{1},0,1,\ldots,n-1\,\}$,

\item[$\cdot$] $\Pi_{m+n}=\{\,\alpha_i\,|\,i\in I_{m+n}\,\}$ : the set of simple roots, where

\begin{equation*}
\alpha_i=
\begin{cases}
-\delta_{\ov{m}}-\delta_{\ov{m-1}}, & \text{if $i=\ov{m}$},\\
\delta_{\ov{k+1}}-\delta_{\ov{k}}, & \text{if $i=\ov{k}\ (\neq \ov{m})$},\\
\delta_{\ov{1}}-\delta_{1}, & \text{if $i=0$},\\
\delta_i-\delta_{i+1}, & \text{if $i=1,\ldots, n-1$}.
\end{cases}
\end{equation*}
\end{itemize}
\noindent
Here, we assume that $P_{m+n}$ has a symmetric bilinear form $(\,\cdot\,|\,\cdot\,)$ such that $(\delta_a|\delta_b)=\delta_{ab}$ and $(\Lambda_{\ov{m}}|\delta_a)=-\frac{1}{2}$ for $a,b\in \J_{m+n}$. The associated Dynkin diagram is 
\begin{center}
\setlength{\unitlength}{0.16in} \hskip -3cm
\begin{picture}(24,5.8)
\put(6,0){\makebox(0,0)[c]{$\bigcirc$}}
\put(6,4){\makebox(0,0)[c]{$\bigcirc$}}
\put(8,2){\makebox(0,0)[c]{$\bigcirc$}}
\put(10.4,2){\makebox(0,0)[c]{$\bigcirc$}}
\put(14.85,2){\makebox(0,0)[c]{$\bigcirc$}}
\put(17.25,2){\makebox(0,0)[c]{$\bigcirc$}}
\put(19.4,2){\makebox(0,0)[c]{$\bigcirc$}}
\put(21.5,2){\makebox(0,0)[c]{$\bigcirc$}}
\put(6.35,0.3){\line(1,1){1.35}} \put(6.35,3.7){\line(1,-1){1.35}}
\put(8.4,2){\line(1,0){1.55}} \put(10.82,2){\line(1,0){0.8}}
\put(13.2,2){\line(1,0){1.2}} \put(15.28,2){\line(1,0){1.45}}
\put(17.7,2){\line(1,0){1.25}} \put(19.8,2){\line(1,0){1.25}}
\put(21.95,2){\line(1,0){1.4}} 
\put(12.5,1.95){\makebox(0,0)[c]{$\cdots$}}
\put(24.5,1.95){\makebox(0,0)[c]{$\cdots$}}
\put(6,5){\makebox(0,0)[c]{\tiny $\alpha_{\ov{m}}$}}
\put(6,-1.2){\makebox(0,0)[c]{\tiny $\alpha_{\ov{m-1}}$}}
\put(8.2,1){\makebox(0,0)[c]{\tiny $\alpha_{\ov{m-2}}$}}
\put(10.4,1){\makebox(0,0)[c]{\tiny $\alpha_{\ov{m-3}}$}}
\put(14.8,1){\makebox(0,0)[c]{\tiny $\alpha_{\ov{1}}$}}
\put(17.15,1){\makebox(0,0)[c]{\tiny $\alpha_0$}}
\put(19.5,1){\makebox(0,0)[c]{\tiny $\alpha_{1}$}}
\put(21.5,1){\makebox(0,0)[c]{\tiny $\alpha_{2}$}}
\end{picture}
\end{center}\vskip 5mm

For $(\lambda,\ell)\in \cP(\mf{d})$, let
\begin{equation*}
\Lambda_{m+\infty}(\lambda,\ell)=\ell\Lambda_{\ov{m}} + \lambda_1\delta_{\ov{m}}+\cdots+\lambda_m\delta_{\ov{1}}+\lambda_{m+1}\delta_{1}+\lambda_{m+2}\delta_{2}+\cdots.
\end{equation*}
 Put $\cP(\mf{d})_{m+n}=\{\,(\lambda,\ell)\in \cP(\mf{d})\,|\,\Lambda_{m+\infty}(\lambda,\ell)\in P_{m+n}\,\}$. For $(\lambda,\ell)\in \cP(\mf{d})_{m+n}$, we write $\Lambda_{m+n}(\lambda,\ell)=\Lambda_{m+\infty}(\lambda,\ell)$.  
Then $\{\,\Lambda_{m+n}(\lambda,\ell)\,|\,(\lambda,\ell)\in \cP(\mf{d})_{m+n}\,\}$ is the set of dominant integral weights for $\mf{d}_{m+n}$. Let $\Lambda_i$ be the $i$th fundamental weight for $i\in I_{m+n}$.

\subsection{Crystal structure on ${\bf T}_{m+n}(\lambda,\ell)$}\label{subsec:crystal structure for m+n}
Put ${\bf T}_{m+n}(a)={\bf T}_{\J_{m+n}}(a)$, $\ov{\bf T}_{m+n}(0)=\ov{\bf T}_{\J_{m+n}}(0)$, ${\bf T}_{m+n}(\lambda,\ell)={\bf T}_{\J_{m+n}}(\lambda,\ell)$, 
 and ${\bf T}_{m+n}^{\rm sp \pm}={\bf T}^{\rm sp \pm}_{\J_{m+n}}$
for $a\in\Z_{\geq 0}$ and $(\lambda,\ell)\in \cP(\mf{d})_{m+n}$.

Let us define an (abstract)  ${\mf d}_{m+n}$-crystal structure on ${\bf T}_{m+n}(\lambda,\ell)$.
We denote the Kashiwara operators on ${\mf d}_{m+n}$-crystals  by $\td{\mathsf e}_i$ and $\td{\mathsf f}_i$ for $i\in I_{m+n}$, and assume that the tensor product rule follows \eqref{lower tensor product rule}.

Recall that $\J_{m+n}$ has a $\gl_{m+n}$-crystal structure with respect to $\td{\mathsf e}_i$ and $\td{\mathsf f}_i$ for $i\in I_{m+n}\setminus\{\ov{m}\}$ as follows;
\begin{equation*}
\begin{split}
&\ov{m}\ \stackrel{^{\ov{m-1}}}{\longrightarrow}\ \ov{m-1}\ \stackrel{^{\ov{m-2}}}{\longrightarrow}
\cdots\stackrel{^{\ov{1}}}{\longrightarrow}\ \ov{1}\ \stackrel{^{0}}{\longrightarrow}\ 1 \
\stackrel{^1}{\longrightarrow}\ 2\stackrel{^2}{\longrightarrow}\cdots 
\end{split}
\end{equation*}
where ${\rm wt}(a)=\delta_a$ for $a\in \J_{m+n}$.  Applying $\td{\mathsf e}_i$ and $\td{\mathsf f}_i$ to the word of a $\J_{m+n}$-semistandard tableau, we have a $\gl_{m+n}$-crystal structure on $SST_{\J_{m+n}}(\lambda/\mu)$ for a skew Young diagram $\lambda/\mu$ \cite{HK,KashNaka}, where $\varepsilon_i$ and $\varphi_i$ are defined in a usual way. For $\lambda\in\cP$, we denote by $H_{\lambda}$ the highest weight element in $SST_{\J_{m+n}}(\lambda)$. 
\vskip 2mm

Let ${\mc B}$ denote one of ${\bf T}_{m+n}^{\rm sp \pm}$, $\ov{\bf T}_{m+n}(0)$, and ${\bf T}_{m+n}(a)$ for $0\leq a\leq m+n-1$. For $T\in {\mc B}$ with ${\rm wt}(T)=(m_s)_{s\in \J_{m+n}}$, let
\begin{equation*}
{\rm wt}(T)=
\begin{cases}
2\Lambda_{\ov{m}}+\sum_{s\in \J_{m+n}}m_s\delta_s, & \text{if  $\mc{B}= \ov{\bf T}_{m+n}(0)$ or ${\bf T}_{m+n}(a)$},\\
\Lambda_{\ov{m}}+\sum_{s\in \J_{m+n}}m_s\delta_s, & \text{if  $\mc{B}={\bf T}_{m+n}^{\rm sp\pm}$}.\\
\end{cases}
\end{equation*}
Since ${\mc B}$ is a set of $\J_{m+n}$-semistandard tableaux, it is a $\gl_{m+n}$-crystal with respect to $\td{\mathsf e}_i$ and $\td{\mathsf f}_i$ for $i\in I_{m+n}\setminus\{\ov{m}\}$. \vskip 2mm

Let us define $\td{\mathsf e}_{\ov{m}}$ and $\td{\mathsf f}_{\ov{m}}$ on ${\mc B}$. Suppose first that $\mc{B}={\bf T}^{\rm sp \pm}_{m+n}$.
For $T\in {\bf T}^{\rm sp}_{m+n}$, let $t_1$ and $t_2$ be the first two top entries of $T$. If $t_1=\ov{m}$ and $t_2=\ov{m-1}$, then we define $\td{\mathsf e}_{\ov{m}}T$ to be the tableau  obtained by removing the domino
$
\resizebox{.05\hsize}{!}{$
{\def\lr#1{\multicolumn{1}{|@{\hspace{.75ex}}c@{\hspace{.75ex}}|}{\raisebox{-.04ex}{$#1$}}}\raisebox{.1ex}
{$\begin{array}{cc} 
\cline{1-1}
\lr{\ov{m}} \\ 
\cline{1-1}
\lr{\ov{m-1}} \\ 
\cline{1-1}
\end{array}$}}$} 
$
from $T$. Otherwise, we define $\td{\mathsf e}_{\ov{m}}T={\bf 0}$. We define $\td{\mathsf f}_{\ov{m}}T$ in a similar way by adding a domino 
$
\resizebox{.05\hsize}{!}{$
{\def\lr#1{\multicolumn{1}{|@{\hspace{.75ex}}c@{\hspace{.75ex}}|}{\raisebox{-.04ex}{$#1$}}}\raisebox{-.0ex}
{$\begin{array}{cc} 
\cline{1-1}
\lr{\ov{m}} \\ 
\cline{1-1}
\lr{\ov{m-1}} \\ 
\cline{1-1}
\end{array}$}}$} 
$
on top of $T$. 
Next, suppose that $\mc{B}={\bf T}_{m+n}(a)$ for $0\leq a\leq m+n$. We regard ${\bf T}_{m+n}(a)$ as a subset of $({\bf T}^{\rm sp}_{m+n})^{\otimes 2}$ by identifying $T=(T^{\tt L},T^{\tt R}) \in {\bf T}_{m+n}(a)$ with $T^{\tt R}\otimes T^{\tt L}$. Then we apply $\td{\mathsf e}_{\ov{m}}$ and $\td{\mathsf f}_{\ov{m}}$ to $T$ following the tensor product rule \eqref{lower tensor product rule}. 
For $T\in \mc{B}$, put
$\varepsilon_{\ov{m}}(T)=\max\{\,r\in\mathbb{Z}_{\geq 0}\,|\,\td{\mathsf e}_{\ov{m}}^r T\neq 0\,\}$ and
$\varphi_{\ov{m}}(T)={\rm wt}(T)+\varepsilon_{\ov{m}}(T)$.
 
\begin{figure}
        \centering
        \begin{subfigure}[b]{0.3\textwidth}
                \includegraphics[width=2cm, height=15cm]{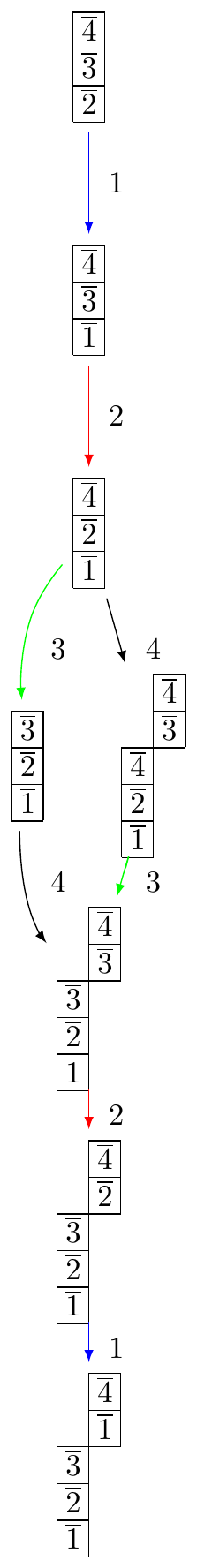}
                \caption{${\bf T}_{4}(3)$}
        \end{subfigure}%
        ~ 
        \begin{subfigure}[b]{0.3\textwidth}
                \includegraphics[width=1.7cm, height=15cm]{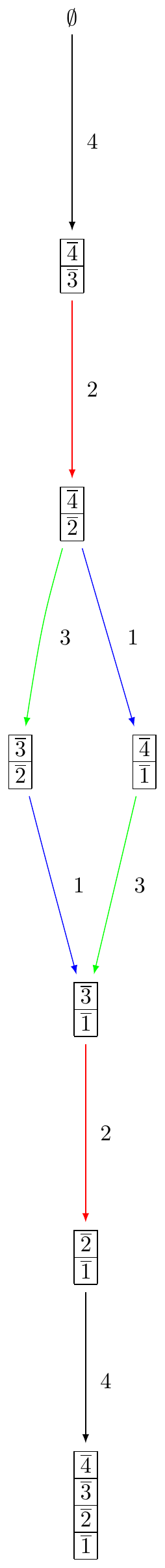}
                \caption{${\bf T}_{4}^{\rm sp +}$}
        \end{subfigure}%
        \begin{subfigure}[b]{0.3\textwidth}
                \includegraphics[width=1.7cm, height=15cm]{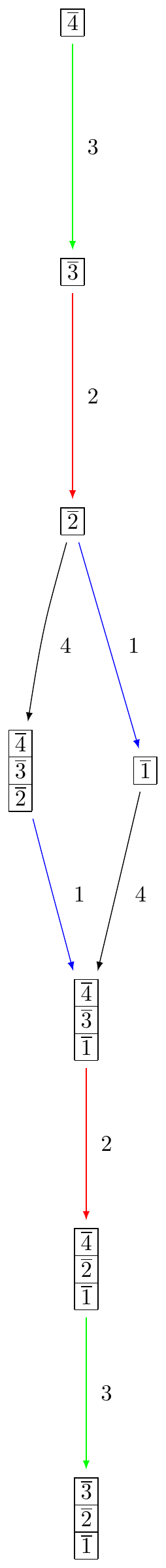}
                \caption{${\bf T}_{4}^{\rm sp -}$}
        \end{subfigure}%
\caption{The crystals of type $D_4$ with $m=4$ and $n=0$.} \label{Graph B}
\end{figure}

\begin{figure}
\includegraphics[width=5.5cm, height=19cm]{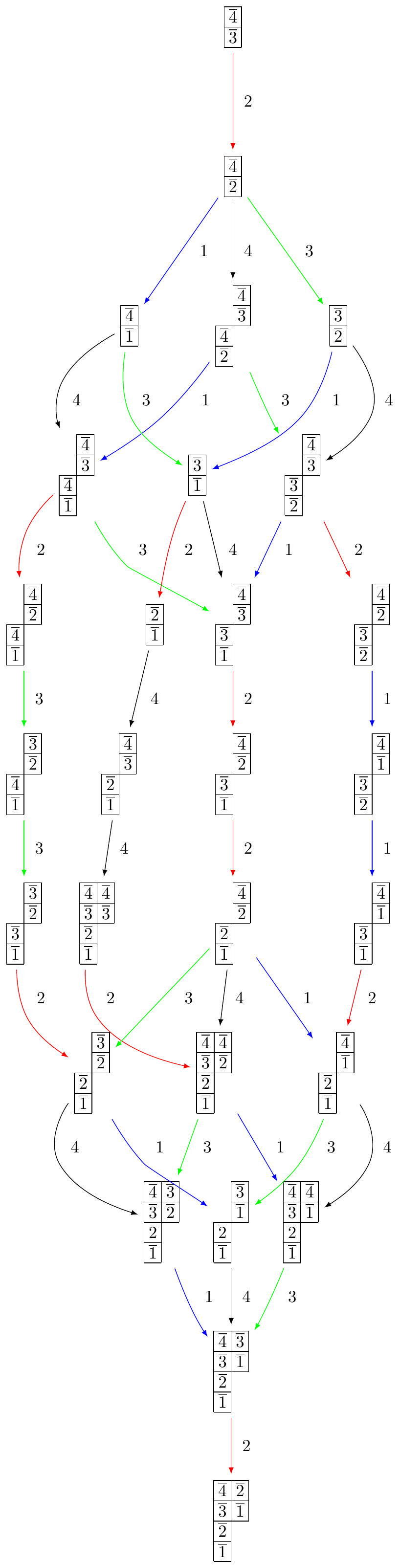}
\caption{The crystal ${\bf T}_{4}(2)$ of type $D_4$ with $m=4$ and $n=0$.} \label{Graph B}
\end{figure}

\begin{lem}\label{Crystal B} 
Under the above hypothesis, $\mc{B}$ is a well-defined $\mf{d}_{m+n}$-crystal with respect to  ${\rm wt}$, $\varepsilon_i$, $\varphi_i$ and $\td{\mathsf e}_{i}$, $\td{\mathsf f}_{i}$ for $i\in I_{m+n}$.
\end{lem}
\pf It is clear that ${\bf T}_{m+n}^{\rm sp\pm}\cup\{{\bf 0}\}$ is invariant under $\td{\mathsf e}_{\ov{m}}$ and $\td{\mathsf f}_{\ov{m}}$, and hence becomes a ${\mf d}_{m+n}$-crystal. So it remains to show that ${\bf T}_{m+n}(a)\cup\{{\bf 0}\}$ or $\ov{\bf T}_{m+n}(0)\cup\{{\bf 0}\}$ is invariant under $\td{\mathsf e}_{i}$, $\td{\mathsf f}_{i}$ for $i\in I_{m+n}$. We will prove the case of ${\bf T}_{m+n}(a)$ since the proof for  $\ov{\bf T}_{m+n}(0)$ is similar.

Let $T\in {\bf T}_{m+n}(a)$ be given with ${\rm sh}(T)=\lambda(a,b,c)$ for some $b,c\in 2\Z_{\geq 0}$. We first observe that $\sigma(T^{\tt L},T^{\tt R})$ is invariant under $\td{\mathsf x}_{i}$  for ${\mathsf x}={\mathsf e}, {\mathsf f}$ and $i\in I_{m+n}\setminus\{\ov{m}\}$ such that $\td{\mathsf e}_{i} T\neq {\bf 0}$ or $\td{\mathsf f}_{i}T\neq {\bf 0}$, since the map \eqref{RSK} is an isomorphism of $(\gl_{m+n},\gl_2)$-bicrystals.

Next, suppose that $\td{\mathsf e}_{\ov{m}}T\neq {\bf 0}$. 
If $\td{\mathsf e}_{\ov{m}}T= T^{\tt R}\otimes (\td{\mathsf e}_{\ov{m}}T^{\tt L})$, then ${\rm sh}(\td{\mathsf e}_{\ov{m}}T)=\lambda(a,b+2,c-2)$. 
Note that the top entry of $T^{\tt R}$ is $\ov{m}$ or $\ov{m-1}$ since otherwise we have $\td{\mathsf e}_{\ov{m}}T = {\bf 0}$ by tensor product rule. 
Then by \eqref{signature condition}   (see also Remark \ref{meaning of signature restriction}) we can check without difficulty that
\begin{equation*}
\sigma(\td{\mathsf e}_{\ov{m}}T)=
\begin{cases}
(a-{\mf r}_T,b+2-{\mf r}_T), & \text{if ${\rm ht}(T^{\tt L})<{\rm ht}(T^{\tt R})$},\\
(a-r,b+2-r), & \text{if ${\rm ht}(T^{\tt L})={\rm ht}(T^{\tt R})$},
\end{cases}
\end{equation*}
for some $r=0,1$. 
Next, if $\td{\mathsf e}_{\ov{m}}T=(\td{\mathsf e}_{\ov{m}} T^{\tt R}) \otimes T^{\tt L}$, then  ${\rm sh}(\td{\mathsf e}_{\ov{m}}T)=\lambda(a,b-2,c)$ and 
\begin{equation*}
\sigma(\td{\mathsf e}_{\ov{m}}T)=
\begin{cases}
(a-{\mf r}_T,b-2-{\mf r}_T), & \text{if ${\rm ht}(T^{\tt L})<{\rm ht}(T^{\tt R})-2$},\\
(a,0), & \text{if ${\rm ht}(T^{\tt L})={\rm ht}(T^{\tt R})-2$}.
\end{cases}
\end{equation*} 
So $\td{\mathsf e}_{\ov{m}}T\in {\bf T}_{m+n}(a)$. Hence ${\bf T}_{m+n}(a)\cup\{{\bf 0}\}$ is invariant under $\td{\mathsf e}_{\ov{m}}$.
By similar arguments, ${\bf T}_{m+n}(a)\cup\{{\bf 0}\}$ is also invariant under $\td{\mathsf f}_{\ov{m}}$.
Therefore, ${\bf T}_{m+n}(a)$ is a subcrystal of $({\bf T}^{\rm sp}_{m+n})^{\otimes 2}$ with respect to ${\rm wt}$, $\varepsilon_i$, $\varphi_i$ and $\td{\mathsf e}_{i}$, $\td{\mathsf f}_{i}$ for $i\in I_{m+n}$. 
\qed

 Let $U_q({\mf d}_{m+n})$ be the quantized enveloping algebra associated to ${\mf d}_{m+n}$ and let $L_q({\mf d}_{m+n},\Lambda)$ be its irreducible highest weight module with highest weight $\Lambda\in P_{m+n}$.
Recall that $\Lambda_{m+n}((0),1)=\Lambda_{\ov{m}}$, $\Lambda_{m+n}((1),1)=\Lambda_{\ov{m-1}}$, and  $\Lambda_{m+n}((1^a),2)
$ represents the other fundamental weights for $2\leq a\leq m+n-1$.

\begin{prop}\label{regularity of level 1 osp crystal}\mbox{}
\begin{itemize}
\item[(1)] ${\bf T}_{m+n}^{\rm sp +}$ is isomorphic to the crystal of $L_q({\mf d}_{m+n},\Lambda_{\ov{m}})$.

\item[(2)] ${\bf T}_{m+n}^{\rm sp -}$ is isomorphic to the crystal of $L_q({\mf d}_{m+n},\Lambda_{\ov{m-1}})$.

\item[(3)] ${\bf T}_{m+n}(a)$ is isomorphic to the crystal of $L_q({\mf d}_{m+n},\Lambda_{m+n}((1^a),2))$ for $0\leq a\leq m+n-1$.
\end{itemize}
\end{prop}
\pf (1) Let $T\in {\bf T}^{\rm sp +}_{m+n}$ be given. Let $(\sigma_a)_{a\in \J_{m+n}}$ be sequence of $\pm$ such that $\sigma_a = -$ if and only if $a$ occurs as an entry of $T$. Then the map sending $T$ to $(\sigma_a)$ is isomorphism of $\mf{d}_{m+n}$-crystals from ${\bf T}^{\rm sp +}_{m+n}$ to the crystal of the spin representation $L_q({\mf d}_{m+n},\Lambda_{\ov{m}})$ (cf. \cite[Section 6.4]{KashNaka}). The proof of (2) is almost the same.

(3) We first claim  that ${\bf T}_{m+n}(a)$ is connected. Let $T\in {\bf T}_{m+n}(a)$ be given. We use induction on the number of boxes in $T\in {\bf T}_{m+n}(a)$, say $|T|$, to show that $T$ is connected to $H_{(1^a)}$, where ${\rm wt}(H_{(1^a)})=\Lambda_{m+n}((1^a),2)$. Suppose that ${\rm sh}(T)=\mu$. Since ${\bf T}_{m+n}(a)$ is a $\gl_{m+n}$-crystal, $T$ is connected to $H_\mu$. If $T^{\tt R}$ is empty, then ${\rm ht}(T)=a$ and $T=H_{(1^a)}$. If $T^{\tt R}$ is not empty, then it has a domino 
$
\resizebox{.05\hsize}{!}{$
{\def\lr#1{\multicolumn{1}{|@{\hspace{.75ex}}c@{\hspace{.75ex}}|}{\raisebox{-.04ex}{$#1$}}}\raisebox{0.5ex}
{$\begin{array}{cc} 
\cline{1-1}
\lr{\ov{m}} \\ 
\cline{1-1}
\lr{\ov{m-1}} \\ 
\cline{1-1}
\end{array}$}}$} 
$\,. 
Hence $\td{\mathsf e}_{\ov{m}}T\neq {\bf 0}$ and $|\td{\mathsf e}_{\ov{m}}T|=|T|-2$, which completes our induction. 

Since ${\bf T}_{m+n}^{\rm sp}$ is a regular crystal and ${\bf T}_{m+n}(a)\subset ({\bf T}_{m+n}^{\rm sp})^{\otimes 2}$, ${\bf T}_{m+n}(a)$ is also regular, which implies that it is isomorphic to the crystal of $L_q({\mf d}_{m+n},\Lambda_{m+n}((1^a),2))$. \qed\vskip 2mm

Let $(\lambda,\ell)\in \cP({\mf d})_{m+n}$ be given. We keep the notations in \eqref{Length of tuples}. We regard ${\bf T}_{m+n}(\lambda,\ell)$ as a subset of
\begin{equation}\label{tensor product of fundamental crystals}
\begin{cases}
({\bf T}^{\rm sp +}_{m+n})^{\otimes r_+}\otimes ({\bf T}_{m+n}(0))^{\otimes q_+}\otimes {\bf T}_{m+n}(\nu_{M_+})\otimes\cdots \otimes {\bf T}_{m+n}(\nu_{1}), & \text{if $\ell-2\lambda_1\geq 0$},\\
({\bf T}^{\rm sp -}_{m+n})^{\otimes r_-}\otimes (\ov{\bf T}_{m+n}(0))^{\otimes q_-}\otimes {\bf T}_{m+n}(\ov{\nu}_{M_-})\otimes\cdots \otimes {\bf T}_{m+n}(\ov{\nu}_{1}), & \text{if $\ell-2\lambda_1\leq 0$},\\
\end{cases}
\end{equation}
by identifying ${\bf T}=(T_L,\ldots, T_0) \in {\bf T}_{m+n}(\lambda,\ell)$ with $T_0\otimes \cdots \otimes T_L$, and apply $\td{\mathsf{e}}_{i}, \td{\mathsf{f}}_{i}$ on ${\bf T}_{m+n}(\lambda,\ell)$ for $i\in I_{m+n}$. We assume that $({\bf T}^{\rm sp \pm}_{m+n})^{\otimes r_\pm}$ is empty or trivial crystal when $r_\pm=0$. Then we have the following.

\begin{thm}\label{crystal invariance of osp tableaux} 
For $(\lambda,\ell)\in \cP({\mf d})_{m+n}$,  
\begin{itemize}
\item[(1)] ${\bf T}_{m+n}(\lambda,\ell)\cup\{{\bf 0}\}$ is invariant under $\td{\mathsf{e}}_{i}$ and $ \td{\mathsf{f}}_{i}$ for $i\in I_{m+n}$,

\item[(2)] ${\bf T}_{m+n}(\lambda,\ell)$ is a connected $\mf{d}_{m+n}$-crystal with highest weight $\Lambda_{m+n}(\lambda,\ell)$.
\end{itemize}
\end{thm}

Theorem \ref{crystal invariance of osp tableaux} immediately implies the following new combinatorial realization of crystal of $L_q({\mf d}_{m+n},\Lambda_{m+n}(\lambda,\ell))$, which plays a crucial role in this paper. The proof of Theorem \ref{crystal invariance of osp tableaux} is given in Section \ref{appendix:A}. 
 
\begin{thm}\label{character formula for m+n}
For $(\lambda,\ell)\in \cP({\mf d})_{m+n}$, ${\bf T}_{m+n}(\lambda,\ell)$ is isomorphic to the crystal of $L_q({\mf d}_{m+n},\Lambda_{m+n}(\lambda,\ell))$.
\end{thm}
\pf By Proposition \ref{regularity of level 1 osp crystal}, ${\bf T}_{m+n}(a)$, $\ov{\bf T}_{m+n}(0)$, and ${\bf T}_{m+n}^{\rm sp \pm}$ are  regular crystals and so is the crystal \eqref{tensor product of fundamental crystals}. By Theorem \ref{crystal invariance of osp tableaux}, ${\bf T}_{m+n}(\lambda,\ell)$ is a regular connected crystal with highest weight $\Lambda_{m+n}(\lambda,\ell)$, and hence it is isomorphic to the crystal of $L_q({\mf d}_{m+n},\Lambda_{m+n}(\lambda,\ell))$.  \qed\vskip 2mm

Let  $\Z[P_{m+n}]$ be a group ring of $P_{m+n}$ with a $\Z$-basis  $\{\,e^\mu\,|\,\mu\in P_{m+n}\,\}$.  Put $z=e^{\Lambda_{\ov{m}}}$ and $x_a=e^{\delta_a}$ for $a\in \J_{m+n}$. By Theorem \ref{character formula for m+n} we have

\begin{cor}\label{character formula for m+n-2}
For $(\lambda,\ell)\in \cP({\mf d})_{m+n}$, we have $${\rm ch}L_q({\mf d}_{m+n},\Lambda_{m+n}(\lambda,\ell))=S_{(\lambda,\ell)}(x_{\J_{m+n}}).$$
\end{cor}

\subsection{Character of a highest weight  module}

Now,  we have the following combinatorial character formula for the irreducible highest weight module with highest weight $\Lambda_{m|n}(\lambda,\ell)$ for $(\lambda,\ell)\in \cP({\mf d})_{m|n}$, which is the main result in this section.

\begin{thm}\label{character formula for m|n}
For $(\lambda,\ell)\in \cP({\mf d})_{m|n}$, we have $${\rm ch}L_q({\mf d}_{m|n},\Lambda_{m|n}(\lambda,\ell))=S_{(\lambda,\ell)}(x_{\J_{m|n}}).$$
That is, the weight generating function of ortho-symplectic tableaux of type $D$ and shape $(\lambda,\ell)$ is equal to the character of $L_q({\mf d}_{m|n},\Lambda_{m|n}(\lambda,\ell))$.
\end{thm}
\pf 
By Corollary \ref{character formula for m+n-2} and Theorem \ref{Schur positivity}, we have
\begin{equation*}
\begin{split}
{\rm ch}L_q({\mf d}_{m+\infty},\Lambda_{m+\infty}(\lambda,\ell))&=S_{(\lambda,\ell)}(x_{\J_{m+\infty}})=z^{\ell}\sum_{\mu\in\cP}K_{\mu\,(\lambda,\ell)}s_\mu(x_{\J_{m+\infty}}).
\end{split}
\end{equation*}
Hence by considering  the classical limit of $L_q({\mf d}_{m|\infty},\Lambda_{m|\infty}(\lambda,\ell))$ (see also \cite[Section 4]{K13}) and super duality \cite[Theorems 4.6 and 5.4]{CLW}, we have
\begin{equation*}
\begin{split}
{\rm ch}L_q({\mf d}_{m|\infty},\Lambda_{m|\infty}(\lambda,\ell))
&=z^{\ell}\sum_{\mu\in\cP}K_{\mu\,(\lambda,\ell)}s_\mu(x_{\J_{m|\infty}})=S_{(\lambda,\ell)}(x_{\J_{m|\infty}}).
\end{split}
\end{equation*}
In particular, ${\rm ch}L_q({\mf d}_{m|n},\Lambda_{m|n}(\lambda,\ell))$ is obtained by specializing $x_a=0$ for $a>n+1$, which is equal to $S_{(\lambda,\ell)}(x_{\J_{m|n}})$. \qed

\section{Crystal base of a highest weight module in $\mc{O}_q^{int}(m|n)$}\label{Crystal base for m|n}
In this section, we prove that $L_q(\mf{d}_{m|n},\Lambda_{m|n}(\lambda,\ell))$ is an irreducible module in $\mc{O}_q^{int}(m|n)$ and it has a unique crystal base for $(\lambda,\ell)\in\cP(\mf{d})_{m|n}$.

\subsection{Crystal structure of ${\bf T}_{m|n}(\lambda,\ell)$}\label{Crystal structure on V_q}
Let $U_q(\mf{gl}_{m|n})=\langle\, e_i, f_i, q^{\pm E_a}\,\big|\,i\in I_{m|n}\setminus\{\ov{m}\}, \ a\in \J_{m|n}\,\rangle$
be the subalgebra of $U_q(\mf{d}_{m|n})$ isomorphic to quantized enveloping algebras associated to general linear Lie superalgebras $\gl_{m|n}$ \cite{Ya}.

We understand $\J_{m|n}$ as the crystal of the natural representation of $U_q(\gl_{m|n})$, where
\begin{equation*}
\begin{split}
&\ov{m}\ \stackrel{^{\ov{m-1}}}{\longrightarrow}\ \ov{m-1}\ \stackrel{^{\ov{m-2}}}{\longrightarrow}
\cdots\stackrel{^{\ov{1}}}{\longrightarrow}\ \ov{1}\ \stackrel{^0}{\longrightarrow}\ \tfrac{1}{2} \
\stackrel{^{\hf}}{\longrightarrow}\ \tfrac{3}{2} \stackrel{{}^{\frac{3}{2}}}{\longrightarrow}\cdots  
\end{split}
\end{equation*}
with ${\rm wt}(a)=\delta_a$ for $a\in \J_{m|n}$ \cite{BKK}, and each non-empty word $w=w_1\cdots w_r$ with letters in $\J_{m|n}$ as $w_1\otimes\cdots\otimes w_r\in (\J_{m|n})^{\otimes r}$.

Then for a skew Young diagram $\lambda/\mu$, $SST_{\J_{m|n}}(\lambda/\mu)$ has an (abstract) $\gl_{m|n}$-crystal structure \cite[Theorem 4.4]{BKK}, where  $\te_i$ and $\tf_i$ are defined via the map $SST_{\J_{m|n}}(\lambda/\mu) \rightarrow \bigsqcup_{r\geq 0}(\J_{m|n})^{\otimes r}$ sending $T$ to $w^{\rm rev}(T)$, the reverse word of $w(T)$.
It is known \cite[Theorem 5.1]{BKK} that 
for $\lambda\in\cP$ with $\lambda_{m+1}\leq n$, $SST_{\J_{m|n}}(\lambda)$ is isomorphic to the crystal of an irreducible polynomial $U_q(\gl_{m|n})$-module with highest weight $\Lambda_{m|n}(\lambda,0)\in P_{m|n}$. We denote by $H^\natural_\lambda$ the highest weight element with weight $\Lambda_{m|n}(\lambda,0)$, called a {\em genuine highest weight element} \cite[Section 4.2]{BKK}.
\vskip 2mm

\begin{rem}\label{crystal for gl type-2} 
{\rm
As in \cite{K13}, our convention for a crystal base of a $U_q(\gl_{m|n})$-module is different from \cite{BKK}. In our setting, it is a upper crystal base as a $U_q(\gl_{m|0})$-module and a lower crystal base as a $U_q(\gl_{0|n})$-module  (see \cite[Remarks 5.1 and 8.1]{K13} for more details).  }
\end{rem}

Put
${\bf T}^{\rm sp}_{m|n}={\bf T}^{\rm sp}_{\J_{m|n}}$ 
and ${\bf T}^{\rm sp \pm}_{m|n}={\bf T}^{\rm sp \pm}_{\J_{m|n}}$, which are clearly $\gl_{m|n}$-crystals.
We also have an $I_{m|n}$-colored oriented graph structure on
${\bf T}^{\rm sp}_{m|n}$,
where $\td{e}_{\ov{m}}$ (resp. $\td{f}_{\ov{m}}$) is defined by adding (resp. removing) an domino 
$
\resizebox{.05\hsize}{!}{$
{\def\lr#1{\multicolumn{1}{|@{\hspace{.75ex}}c@{\hspace{.75ex}}|}{\raisebox{-.04ex}{$#1$}}}\raisebox{0.5ex}
{$\begin{array}{cc} 
\cline{1-1}
\lr{\ov{m}} \\ 
\cline{1-1}
\lr{\ov{m-1}} \\ 
\cline{1-1}
\end{array}$}}$} 
$\, as in the case of ${\bf T}_{m+n}^{\rm sp}$ (see Section \ref{subsec:crystal structure for m+n}).
Then ${\bf T}_{m|n}^{\rm sp \pm}\cup\{{\bf 0}\}$ is invariant under $\te_i$ and $\tf_i$ for $i\in I_{m|n}$.

Let ${\bf m}=(m_a)\in \B^+$ be given (see Section \ref{subsec:Fock space}). Let  $T({\bf m})\in SST_{\J_{m|n}}(1^{d})$ be the unique tableaux such that the entries in $T({\bf m})$ are the $a$'s with $m_a\neq 0$ counting multiplicity as many as $m_a$ times, where $d =\sum_{a\in\J_{m|n}}m_a$. 
Since $\B^+$ may be regarded as a crystal of  a $U_q(\mf{d}_{m|n})$-module $\mathscr{V}_{q}$ by \eqref{crystal of V_q}, we can check  the following (see the proof of \cite[Theorem 5.6]{K13}).

\begin{lem}\label{B=T} The map
$\Psi^+ : \B^+ \longrightarrow {\bf T}^{\rm sp}_{m|n}$ 
given by $\Psi^+({\bf m})=T({\bf m})$ is a bijection which commute with $\te_i$ and $\tf_i$ for $i\in I_{m|n}$.
Hence,  we may regard  ${\bf T}_{m|n}^{\rm sp}$  as a crystal of  $\mathscr{V}_{q}$, where ${\rm wt}$, $\varepsilon_i$ and $\varphi_i$ $(i\in I_{m|n})$ are induced from those on  $\B^+$  via  $\Psi^+$.
\end{lem}

Next, put ${\bf T}_{m|n}(a)={\bf T}_{\J_{m|n}}(a)$, $\ov{\bf T}_{m|n}(0)=\ov{\bf T}_{\J_{m|n}}(0)$, and ${\bf T}_{m|n}(\lambda,\ell)={\bf T}_{\J_{m|n}}(\lambda,\ell)$   
for $a\in\Z_{\geq 0}$ and $(\lambda,\ell)\in \cP(\mf{d})_{m|n}$.
We regard
${\bf T}_{m|n}(a), \ov{\bf T}_{m|n}(0) \subset  ({\bf T}_{m|n}^{\rm sp})^{\otimes 2}$
by identifying $T$ with $T^{\tt L}\otimes T^{\tt R}$ (see Remark \ref{crystal for gl type-2}). 

\begin{lem}\label{Invariance of T m|n (a)}
${\bf T}_{m|n}(a)\cup\{{\bf 0}\}$ $(a\geq 0)$ and $\ov{\bf T}_{m|n}(0)\cup\{{\bf 0}\}$ are invariant under $\td{e}_i$ and $\td{f}_i$ for $i\in I_{m|n}$.
\end{lem}
\pf The proof is almost the same as in Lemma \ref{Crystal B}.
\qed\vskip 2mm

Keeping the notations in \eqref{Length of tuples}, we consider ${\bf T}_{m|n}(\lambda,\ell)$ for $(\lambda,\ell)\in \cP({\mf d})_{m|n}$ as a subset of
\begin{equation*}\label{tensor product of fundamental crystals-2}
\begin{cases}
{\bf T}_{m|n}(\nu_{1})\otimes\cdots \otimes {\bf T}_{m|n}(\nu_{M_+})\otimes  \left({\bf T}_{m|n}(0)\right)^{q_+} \otimes \left({\bf T}^{\rm sp +}_{m|n}\right)^{\otimes r_+}, & \text{if $\ell-2\lambda_1\geq 0$},\\
{\bf T}_{m|n}(\ov{\nu}_{1})\otimes\cdots \otimes {\bf T}_{m|n}(\ov{\nu}_{M_-})\otimes  \left(\ov{\bf T}_{m|n}(0)\right)^{q_-}\otimes \left({\bf T}^{\rm sp -}_{m|n}\right)^{\otimes r_-}, & \text{if $\ell-2\lambda_1\leq0$},
\end{cases}
\end{equation*}
by identifying ${\bf T}=(T_L,\ldots, T_0) \in {\bf T}_{m|n}(\lambda,\ell)$ with $T_L\otimes \cdots \otimes T_0$, and apply $\te_{i}$ and $\tf_{i}$ on ${\bf T}_{m|n}(\lambda,\ell)$ for $i\in I_{m|n}$. We put
\begin{equation}\label{Highest weight element in T m|n lambda}
{\bf H}^{\natural}_{(\lambda,\ell)}= H_{L}\otimes\cdots\otimes H_{0},
\end{equation}
where $H_{k}$ is empty   for $0\leq k\leq q_+$ when $\ell-2\lambda_1\geq 0$, $H_{0}=$ \resizebox{.035\hsize}{!}{$\boxed{\ov{m}}$}\  and $H_{k}=$
$
\resizebox{.07\hsize}{!}{$
{\def\lr#1{\multicolumn{1}{|@{\hspace{.75ex}}c@{\hspace{.75ex}}|}{\raisebox{-.04ex}{$#1$}}}\raisebox{0.25ex}
{$\begin{array}{cc} 
\cline{1-2}
\lr{\ov{m}} &  \lr{\ov{m}} \\
\cline{1-2}
\end{array}$}}$} 
$\, for $1\leq k\leq q_-$ when $\ell-2\lambda_1\leq 0$, and $H_{q_\pm +k}\in SST_{\J_{m|n}}(1^{a_k})$ for $1 \leq k\leq M_\pm$ ($a_k$ as in \eqref{notations for lengths}) are the unique tableaux such that
\begin{equation*}
(H_L \rightarrow (\cdots(H_{2}\rightarrow H_{0})))=H^\natural_\lambda.
\end{equation*}
We remark that 
$H_{q_\pm+k}$ is not necessarily equal to $H^\natural_{(1^{a_k})}$ in $SST_{\J_{m|n}}(1^{a_k})$ unlike the case of $\J_{m+n}$-semistandard tableaux (cf.~\cite[Example 5.8]{KK}). Indeed $H_{L-k+1}$ is the $k$th column of $H^\natural_\lambda$ from the left for $k\geq 1$.

\begin{thm}\label{connectedness of T m|n (lambda,ell)} For $(\lambda,\ell)\in \cP({\mf d})_{m|n}$,
\begin{itemize}
\item[(1)] ${\bf T}_{m|n}(\lambda,\ell)\cup\{{\bf 0}\}$  is invariant under $\td{e}_i$ and $\td{f}_i$ for $i\in I_{m|n}$,

\item[(2)] ${\bf T}_{m|n}(\lambda,\ell)$ is a connected $I_{m|n}$-colored oriented graph with a highest weight element ${\bf H}^{\natural}_{(\lambda,\ell)}$ of weight $\Lambda_{m|n}(\lambda,\ell)$.
\end{itemize}
\end{thm}
\pf (1) Since the proof is similar to that of Theorem \ref{crystal invariance of osp tableaux} in Section \ref{appendix:A}, we give a brief sketch of it.  
In this case, we have ${\bf M}_{\J_{m|n}\times 1}={\bf T}^{\rm sp}_{m|n}$, and hence ${\bf M}_{\J_{m+n}\times \ell}=({\bf T}_{m|n}^{\rm sp})^{\otimes \ell}$ has a $\gl_{m|n}$-crystal structure, where we identify ${\bf m}=[\,{\bf m}^{(\ell)}:\cdots:{\bf m}^{(1)}\,]\in {\bf M}_{\J_{m+n}\times \ell}$ with ${\bf m}^{(\ell)}\otimes \cdots\otimes {\bf m}^{(1)}\in ({\bf T}_{m|n}^{\rm sp})^{\otimes \ell}$. Then $\td{e}_{i}$, $\td{f}_{i}$ on ${\bf T}_{m|n}(\lambda,\ell)$ coincide with those on ${\bf M}_{\J_{m|n}\times \ell}$ for $i\in I_{m|n}\setminus \{\ov{m}\}$ since ${\bf T}_{m|n}(\lambda,\ell)\subset ({\bf T}_{m|n}^{\rm sp})^{\otimes \ell}$.
Note that ${\bf M}_{\J_{m|n}\times \ell}$ is a $(\gl_{m|n},\gl_\ell)$-bicrystal and the map \eqref{RSK} is an isomorphism of bicrystals \cite{K07}. Hence it follows from \eqref{Pieri} that $\td{x}_i{\bf T}_{m|n}(\lambda,\ell)\subset {\bf T}_{m|n}(\lambda,\ell)\cup\{{\bf 0}\}$ for $x=e, f$ and $i\in  I_{m|n}\setminus \{\ov{m}\}$.
The proof for $\td{x}_{\ov{m}}{\bf T}_{m|n}(\lambda,\ell)\subset {\bf T}_{m|n}(\lambda,\ell)\cup\{{\bf 0}\}$  for $x=e, f$ is the same as in Lemmas \ref{lem:eT-1}--\ref{lem nonspin spin}.
 
(2)  Let ${\bf T}=(T_L,\ldots, T_0)\in {\bf T}_{m|n}(\lambda,\ell)$ be given. As in Lemma \ref{lem cnn}, we use induction on $|{\bf T}|=\sum_{k=0}^L|T_k|$ to show that ${\bf T}$ is connected to ${\bf H}^\natural_{(\lambda,\ell)}$.  By \cite[Theorem 4.8]{BKK}, we may assume that ${\bf T}$ is a genuine $\mf{gl}_{m|n}$-highest weight element, that is, $P:=(T_{L}\rightarrow ( \cdots( T_1 \rightarrow T_0)))=H^\natural_\mu$ for some $\mu\in \cP$.
We will show that ${\bf T}= {\bf H}^\natural_{(\lambda,\ell)}$ or $\te_{\ov{m}}{\bf T}\neq {\bf 0}$, which implies $|\te_{\ov{m}} {\bf T}|<|{\bf T}|$.
  
Suppose that $\ell-2\lambda_1\geq 0$ and consider $T_0$ or $T_1$ (if $T_0$ is empty). From the  insertion process for $P:=(T_{L}\rightarrow ( \cdots( T_1 \rightarrow T_0)))$, we observe that each $k$-th entry of $T_0$ from the top lie in the $l$-th row of $P$ with $l\leq k$. If $T_0$ is not empty, then it contains a domino\, 
$
\resizebox{.05\hsize}{!}{$
{\def\lr#1{\multicolumn{1}{|@{\hspace{.75ex}}c@{\hspace{.75ex}}|}{\raisebox{-.04ex}{$#1$}}}\raisebox{0.5ex}
{$\begin{array}{cc} 
\cline{1-1}
\lr{\ov{m}} \\ 
\cline{1-1}
\lr{\ov{m-1}} \\ 
\cline{1-1}
\end{array}$}}$} 
$\, 
 since $P=H^\natural_\mu$. This implies that $\te_{\ov{m}}T_0\neq {\bf 0}$ and hence $\te_{\ov{m}}{\bf T}\neq {\bf 0}$. If $T_0$ is empty, then  $T_i$ is empty for $1\leq i\leq q_+$ since $T_{i}\prec T_{i-1}$ for $1\leq i\leq q_+$. 
Suppose that $\ell-2\lambda_1 \leq 0$. By almost the same argument, we conclude that $\te_{\ov{m}}{\bf T}\neq {\bf 0}$ or $T_{0}=$ \resizebox{.035\hsize}{!}{$\boxed{\ov{m}}$} and $T_i=$ 
$
\resizebox{.07\hsize}{!}{$
{\def\lr#1{\multicolumn{1}{|@{\hspace{.75ex}}c@{\hspace{.75ex}}|}{\raisebox{-.04ex}{$#1$}}}\raisebox{0.5ex}
{$\begin{array}{cc} 
\cline{1-2}
\lr{\ov{m}} &  \lr{\ov{m}} \\
\cline{1-2}
\end{array}$}}$} 
$ 
for $1\leq i\leq q_-$.

Now we may assume that $T_{k}$ is empty for $0\leq k\leq q_+$ when $\ell-2\lambda_1\geq 0$, and  $T_{0}=$ \resizebox{.035\hsize}{!}{$\boxed{\ov{m}}$}\  and $T_{k}=$
$
\resizebox{.07\hsize}{!}{$
{\def\lr#1{\multicolumn{1}{|@{\hspace{.75ex}}c@{\hspace{.75ex}}|}{\raisebox{-.04ex}{$#1$}}}\raisebox{0.25ex}
{$\begin{array}{cc} 
\cline{1-2}
\lr{\ov{m}} &  \lr{\ov{m}} \\
\cline{1-2}
\end{array}$}}$} 
$\, for $1\leq k\leq q_-$ when $\ell-2\lambda_1\leq 0$.

Consider $T_{q_\pm+1}$. If $\te_{\ov{m}}{\bf T}={\bf 0}$, then by the same argument as in the proof of Lemma \ref{lem cnn}, we have $T^{\tt R}_{q_\pm+1}$ is empty, and hence ${\rm ht}(T^{\tt L}_{q_\pm+1})=a_1$. Now we can prove inductively that $T^{\tt R}_{q_\pm+k}$ is empty and hence ${\rm ht}(T^{\tt L}_{q_\pm+k})=a_k$ for $1\leq k\leq M_\pm$. Since $T_{i+1}\prec T_i$ for $0\leq i\leq L-1$, $(T_L,\ldots,T_0)$ itself forms a $\J_{m|n}$-semistandard tableau $H^\natural_\mu$. Since $|{\bf T}|=|{\bf H}^\natural_{(\lambda,\ell)}|$ is minimal, we conclude that $\mu=\lambda$ and ${\bf T}={\bf H}^\natural_{(\lambda,\ell)}$. The proof completes.
\qed

\subsection{Main result}\label{Main Theorem}

\begin{lem}\label{Highest weight vector for fundamental weight of type D}
For $a\geq 0$, there exists ${\bf v}_{a}\in \mathscr{V}_q^{\otimes 2}$ such that
\begin{itemize}
\item[(1)] ${\bf v}_a$ is a $U_q(\mf{d}_{m|n})$-highest weight vector of weight $\Lambda_{m|n}((1^a),2)$,

\item[(2)] ${\bf v}_a\in \mathscr{L}^+\otimes \mathscr{L}^+$ and
${\bf v}_a \equiv \psi_{{\bf m}^+(a)}|0\rangle \otimes |0\rangle \pmod{q\mathscr{L}^+\otimes \mathscr{L}^+}$,
\end{itemize}
where ${\bf m}^+(a)\in \B^+$ is given by $\Psi^+({\bf m}^+(a))=H^{\natural}_{(1^a)}$.
\end{lem}
\pf  The proof is similar to that of \cite[Lemma 5.5]{K13}.  If $a=0$, then it is clear that ${\bf v}_0=|0\rangle \otimes |0\rangle$. We assume that $a\geq 1$.
Let ${\bf M}(a)$ be the set of ${\bf m}=[m_{rs}]\in {\bf M}_{\J_{m|n}\times 2}$ satisfying the following conditions:
\begin{itemize}
\item[(1)] $m_{r 1}+m_{r 2}=1$ for $\ov{m}\leq r\leq \ov{l+1}$ where $l=\max\{m-a,0\}$,

\item[(2)] $m_{rs}=0$ for $r> \hf$ and $s=1,2$,

\item[(3)] $m_{\hf 1} + m_{\hf 2}= \max\{0,a-m\} $,

\item[(4)] $\sum_{r\in \J_{m|n}}m_{r 1}$ is even.
\end{itemize}

Let ${\bf m}=[m_{rs}]\in {\bf M}(a)$ be given. We write ${\bf m} \stackrel{\ov{m}}{\rightsquigarrow} {\bf m}'$ if $m_{\ov{m}\,2}=m_{\ov{m-1}\,2}=1$ and ${\bf m}'$ is obtained from ${\bf m}$ by replacing 
\begin{equation*}
\begin{bmatrix}
m_{\ov{m}\,2} &  m_{\ov{m}\,1} \\
m_{\ov{m-1}\,2}  &  m_{\ov{m-1}\,1}
\end{bmatrix}
=\begin{bmatrix}
1 &  0 \\
1  & 0
\end{bmatrix}\ \ \text{with}\ \
\begin{bmatrix}
0 &  1 \\
0  & 1
\end{bmatrix}.
\end{equation*}
For $i\in \{\ov{m-1},\ldots,\ov{1}\}$, we write ${\bf m} \stackrel{i}{\rightsquigarrow} {\bf m}'$ if $m_{\ov{i+1}\,2}=0$, $m_{\ov{i}\,2}=1$ and ${\bf m}'$ is obtained from ${\bf m}$ by replacing
\begin{equation*}
\begin{bmatrix}
m_{\ov{i+1}\,2} &  m_{\ov{i}\,1} \\
m_{\ov{i}2}\,  &  m_{\ov{i}\,1}
\end{bmatrix}
=\begin{bmatrix}
0 &  1 \\
1  & 0
\end{bmatrix}\ \ \text{with}\ \
\begin{bmatrix}
1 &  0 \\
0  & 1
\end{bmatrix}.
\end{equation*}
Similarly, we write ${\bf m} \stackrel{0}{\rightsquigarrow} {\bf m}'$ if $m_{\ov{1}\,2}=0$, $m_{\hf 2}\geq 1$ and ${\bf m}'$ is obtained from ${\bf m}$ by replacing
\begin{equation*}
\begin{bmatrix}
m_{\ov{1}2} &  m_{\ov{1}1} \\
m_{\hf 2}  &  m_{\hf 1}
\end{bmatrix}
=\begin{bmatrix}
0 &  1 \\
u  & v
\end{bmatrix}\ \ \text{with}\ \
\begin{bmatrix}
1 &  0 \\
u-1  & v+1
\end{bmatrix}.
\end{equation*}

Then we have 
\begin{equation}\label{M and M'}
e_{i}\left(\psi_{{\bf m}^{(2)}}|0\rangle\otimes \psi_{{\bf m}^{(1)}}|0\rangle\right)=Q_{{\bf m},{\bf m}'}(q) e_{i}\left(\psi_{{\bf m'}^{(2)}}|0\rangle\otimes \psi_{{\bf m'}^{(1)}}|0\rangle\right),
\end{equation}
for ${\bf m} \stackrel{i}{\rightsquigarrow} {\bf m}'$, where $Q_{{\bf m},{\bf m}'}(q)$ is a monomial in $q$ of positive degree given by
\begin{equation}
Q_{{\bf m},{\bf m}'}(q)=
\begin{cases}
q, & \text{if $i=\ov{m},\ldots,\ov{1}$},\\
(-1)^{|{\rm wt}({\bf m}^{(2)})|}q^{\langle \beta^\vee_0,{\rm wt}({\bf m}^{(1)})\rangle}, & \text{if $i=0$}.\\
\end{cases}
\end{equation}
Recall $|{\rm wt}({\bf m}^{(2)})|$ denotes the degree or parity of ${\rm wt}({\bf m}^{(2)})$ (cf.  \cite[Remark 3.1]{K13}).

Let ${\bf m}(a)\in {\bf M}(a)$ be such that $m_{r1}=0$ for all $r\in \J_{m|n}$, that is,  ${\bf m}(a)^{(2)}={\bf m}^+(a)$ and ${\bf m}(a)^{(1)}$ is trivial.
Then for ${\bf m}\in {\bf M}(a)$, we have
\begin{equation}\label{path from M(a) to M}
{\bf m}(a)={\bf m}_0\stackrel{i_1}{\rightsquigarrow}{\bf m}_1\stackrel{i_2}{\rightsquigarrow}\cdots \stackrel{i_{r}}{\rightsquigarrow}{\bf m}_r={\bf m},
\end{equation}
for some $r\geq 0$, $i_1,\ldots,i_r\in \{\ov{m},\ldots,\ov{1},0\}$ and ${\bf m}_1,\ldots, {\bf m}_{r-1}\in{\bf {\bf M}}(a)$.
Put
\begin{equation*}
h({\bf m})=r, \ \ \ Q_{\bf m}(q)=\prod_{k=0}^{r-1}Q_{{\bf m}_k,{\bf m}_{k+1}}(q).
\end{equation*} 

Note that ${\bf m}\in {\bf M}(a)$ is completely determined by its second column ${\bf m}^{(2)}$, and under this identification the $\{\ov{m},\ldots,\ov{1},0\}$-colored graph structure on ${\bf M}(a)$ with respect to $\ \stackrel{i}{\rightsquigarrow}\ $ coincides with the ${\mf d}_{m|1}$-crystal structure on ${\bf T}^{\rm sp+}_{m|1}$ (see Section \ref{Crystal structure on V_q}). This implies as in \cite[Lemma 8.6]{K13} that $h({\bf m})$ and $Q_{\bf m}(q)$ are independent of a path  \eqref{path from M(a) to M} from ${\bf m}(a)$ to ${\bf m}$.
Put 
\begin{equation*}
{\bf v}_a=\sum_{{\bf m}\in {\bf M}(a)}(-1)^{h({\bf m})}Q_{\bf m}(q)\psi_{{\bf m}^{(2)}}|0\rangle\otimes \psi_{{\bf m}^{(1)}}|0\rangle.
\end{equation*}
Then ${\bf v}_a\in \mathscr{L}^+\otimes \mathscr{L}^+$ and ${\bf v}_a \equiv \psi_{{\bf m}^+(a)}|0\rangle \otimes |0\rangle \pmod{q\mathscr{L}^+\otimes \mathscr{L}^+}$. 

Consider the pairs $({\bf m},{\bf m}')$ for ${\bf m}, {\bf m}'\in{\bf M}(a)$ with ${\bf m}\stackrel{i}{\rightsquigarrow} {\bf m}'$ for some $i\in I_{m|n}$. We see that any ${\bf m}\in {\bf M}(a)$ with $e_{i}(\psi_{{\bf m}^{(2)}}|0\rangle\otimes \psi_{{\bf m}^{(1)}}|0\rangle) \neq 0$ belongs to one of these pairs.
Since $h({\bf m}')=h({\bf m})+1$ and $Q_{{\bf m}'}(q)=Q_{{\bf m}}(q)Q_{{\bf m},{\bf m}'}(q)$, we have by \eqref{M and M'}
\begin{equation*}
\begin{split}
&e_{i}\left\{(-1)^{h({\bf m})}Q_{{\bf m}}(q)\psi_{{\bf m}^{(2)}}|0\rangle\otimes \psi_{{\bf m}^{(1)}}|0\rangle+(-1)^{h({\bf m}')}Q_{{\bf m}'}(q)\psi_{{\bf m}^{'(2)}}|0\rangle\otimes \psi_{{\bf m}^{'(1)}}|0\rangle\right\}  \\&=(-1)^{h({\bf m})}Q_{{\bf m}}(q)\left\{e_{i}(\psi_{{\bf m}^{(2)}}|0\rangle\otimes \psi_{{\bf m}^{(1)}}|0\rangle)- Q_{{\bf m},{\bf m}'}(q) e_{i}(\psi_{{\bf m}^{'(2)}}|0\rangle\otimes \psi_{{\bf m}^{'(1)}}|0\rangle) \right\}\\ &=0.
\end{split}
\end{equation*}
This implies that $e_{i}{\bf v}_a =0$ for all $i\in I_{m|n}$, and hence it is a $U_q({\mf d}_{m|n})$-highest weight vector. \qed\vskip 2mm

\begin{prop}\label{Crystal base of a fundamental weight of type D} 
For $a\geq 0$, let ${\bf v}_a$ be as in Lemma \ref{Highest weight vector for fundamental weight of type D}. Then  $U_q(\mf{d}_{m|n}) {\bf v}_a$ is isomorphic to $L_q(\mf{d}_{m|n},\Lambda_{m|n}((1^a),2))$, and it has a crystal base $(\mathscr{L}(a),\mathscr{B}(a))$, where
\begin{equation*}
\begin{split}
\mathscr{L}(a)&=\sum\mathbb{A}\td{x}_{i_1}\cdots\td{x}_{i_r}{\bf v}_a,\\
\mathscr{B}(a)&=\{\,\pm\td{x}_{i_1}\cdots\td{x}_{i_r}{\bf v}_a \!\!\!\pmod{q\mathscr{L}(a)} \, \}\setminus\{0\},\\
\end{split}
\end{equation*}
with $r\geq 0$, $i_1,\ldots,i_r\in I_{m|n}$, and $x=e, f$ for each $i_k$. The crystal $\mathscr{B}(a)/\{\pm1\}$ is isomorphic to ${\bf T}_{m|n}(a)$.
\end{prop}
\pf By Lemma \ref{Highest weight vector for fundamental weight of type D} and Theorem \ref{complete reducibility}, we have
\begin{equation*}
U_q(\mf{d}_{m|n}) {\bf v}_a  \cong L_q(\mf{d}_{m|n},\Lambda_{m|n}((1^a),2)).
\end{equation*} 
Also, it follows from the same argument in \cite[Proposition 8.7]{K13} that $(\mathscr{L}(a),\mathscr{B}(a))$ is a crystal base of $L_q(\mf{d}_{m|n},\Lambda_{m|n}((1^a),2))$, and $\mathscr{B}(a)/\{\pm1\}$ is isomorphic to ${\bf T}_{m|n}(a)$.\qed\vskip 2mm

Now we are ready to state and prove our main theorem in this paper.
\begin{thm}\label{Existence of crystal base}
For $(\lambda,\ell)\in \cP({\mf d})_{m|n}$, $L_q(\mf{d}_{m|n},\Lambda_{m|n}(\lambda,\ell))$ is an irreduble $U_q(\mf{d}_{m|n})$-module in $\mc{O}^{int}_q(m|n)$, and it has a unique crystal base up to scalar multiplication, whose crystal is isomorphic to ${\bf T}_{m|n}(\lambda,\ell)$.
\end{thm}
\pf  Let $(\lambda,\ell)\in \cP({\mf d})_{m|n}$ be given  with $L$ as in \eqref{Length of tuples}. 
Let $V_{(\lambda,\ell)}=V_L\otimes \cdots\otimes V_0$ be a $U_q(\gl_{m|n})$-module, where 
\begin{equation}
V_i=
\begin{cases}
U_q(\gl_{m|n}){\bf v}_{\mu_{L-i+1}}, & \text{if $\ell-2\lambda_1\geq 0$ and $q_++1\leq i\leq L$},\\
U_q(\gl_{m|n})|0\rangle\otimes |0\rangle, & \text{if $\ell-2\lambda_1\geq 0$ and $1\leq i\leq q_+$},\\
U_q(\gl_{m|n})|0\rangle  , & \text{if $\ell-2\lambda_1\geq 0$ and $i=0$},\\
U_q(\gl_{m|n}){\bf v}_{\ov{\mu}_{L-i+1}}, & \text{if $\ell-2\lambda_1< 0$ and $q_-+1\leq i\leq L$},\\
U_q(\gl_{m|n})\psi_{\ov{m}}|0\rangle\otimes \psi_{\ov{m}}|0\rangle, & \text{if $\ell-2\lambda_1< 0$ and $1\leq i\leq q_-$},\\
U_q(\gl_{m|n})\psi_{\ov{m}}|0\rangle, & \text{if $\ell-2\lambda_1< 0$ and $i=0$}.
\end{cases}
\end{equation}
Here we assume that $V_0$ is trivial when $r_+$ or $r_-$ is 0.
Then $V_i$ is isomorphic to an irreducible polynomial $U_q(\mf{gl}_{m|n})$-module, and  
$V_{(\lambda,\ell)}$ is a completely reducible  $U_q(\mf{gl}_{m|n})$-module with a crystal base \cite{BKK}. 
Also, by \cite[Theorem 5.6]{K13} and  Proposition \ref{Crystal base of a fundamental weight of type D}, 
we may assume that the crystal lattice  of $V_{(\lambda,\ell)}$ is contained in a tensor product of $\mathscr{L}(a)$'s and $\mathscr{L}^+$'s, say $\mathscr{L}$.
 
The rest of the proof is the same as in  \cite[Theorem 8.8]{K13}, which we refer the reader to for more details.  
First, from the decomposition of $V_{(\lambda,\ell)}$ (cf.\cite{KK}), we can find a unique $U_q(\mf{gl}_{m|n})$-highest weight vector ${\bf v}_{(\lambda,\ell)}$ in $V_{(\lambda,\ell)}$ (up to scalar multiplication) such that
$U_q(\mf{gl}_{m|n}){\bf v}_{(\lambda,\ell)}$ is isomorphic to the irreducible $U_q(\mf{gl}_{m|n})$-module with highest weight $\Lambda_{m|n}(\lambda,\ell)$ and ${\bf v}_{(\lambda,\ell)}\not\equiv 0 \pmod{q\mathscr{L}}$. Since $e_{\ov{m}}V_i=0$ for all $i$ by construction,  ${\bf v}_{(\lambda,\ell)}$ is a $U_q(\mf{d}_{m|n})$-highest weight vector and
$U_q(\mf{d}_{m|n}){\bf v}_{(\lambda,\ell)}\cong L_q(\mf{d}_{m|n},\Lambda_{m|n}(\lambda,\ell))$ with
${\bf v}_{(\lambda,\ell)}\equiv \pm {\bf H}^\natural_{(\lambda,\ell)} \pmod{q\mathscr{L}}$
(see \eqref{Highest weight element in T m|n lambda}).
Next if we put
\begin{equation*}
\begin{split}
\mathscr{L}(\lambda,\ell)&=\sum\mathbb{A}\td{x}_{i_1}\cdots\td{x}_{i_r}{\bf v}_{(\lambda,\ell)}\subset \mathscr{L},\\
\mathscr{B}(\lambda,\ell)&=\{\,\pm\td{x}_{i_1}\cdots\td{x}_{i_r}{\bf v}_{(\lambda,\ell)} \!\!\!\pmod{q\mathscr{L}(\lambda,\ell)}\, \}\setminus\{0\},
\end{split}
\end{equation*}
where $r\geq 0$, $i_1,\ldots,i_r\in I_{m|n}$, and $x=e, f$ for each $i_k$,
then  it follows from Lemma \ref{B=T}, Proposition \ref{Crystal base of a fundamental weight of type D}, and Theorems \ref{character formula for m|n} and \ref{connectedness of T m|n (lambda,ell)} that $(\mathscr{L}(\lambda,\ell),\mathscr{B}(\lambda,\ell))$ is a crystal base of $L_q(\mf{d}_{m|n},\Lambda_{m|n}(\lambda,\ell))$, and the map
\begin{equation}\label{Map from B(lambda,ell)}
\Psi_{(\lambda,\ell)} : \left(\mathscr{B}(\lambda,\ell)/\{\pm1\}\right)\cup\{0\} \longrightarrow {\bf T}_{m|n}(\lambda,\ell)\cup\{{\bf 0}\},
\end{equation}
given by $\td{x}_{i_1}\cdots\td{x}_{i_r}{\bf v}_{(\lambda,\ell)}\longmapsto \td{x}_{i_1}\cdots\td{x}_{i_r} {\bf H}^\natural_{(\lambda,\ell)}$
for $r\geq 0$, $i_1,\ldots,i_r\in I_{m|n}$ and $x=e,f$ is a weight preserving bijection which commutes with $\td{e}_i$ and $\td{f}_i$ for $i\in I_{m|n}$.
Finally, the uniqueness of a crystal base of $L_q(\mf{d}_{m|n},\Lambda_{m|n}(\lambda,\ell))$ follows from Theorem \ref{connectedness of T m|n (lambda,ell)}(2) and \cite[Lemma 2.7(iii) and (iv)]{BKK}.\qed

\begin{cor}
Each $U_q(\mf{d}_{m|n})$-module in $\mc{O}^{int}_q(m|n)$ has a crystal base.
\end{cor}

\begin{cor}
Each highest weight  $U_q(\mf{d}_{m|n})$-module in $\mc{O}^{int}_q(m|n)$ is a direct summand of $\mathscr{V}_q^{\otimes M}$ for some $M\geq 1$.
\end{cor}

\section{Proof of Theorem \ref{crystal invariance of osp tableaux}}\label{appendix:A}

Let $(\lambda,\ell)\in \cP({\mf d})_{m+n}$ be given. Since ${\bf M}_{\J_{m+n}\times 1}={\bf T}^{\rm sp}_{m+n}$, we may understand ${\bf M}_{\J_{m+n}\times \ell}$ as a $\gl_{m+n}$-crystal, where ${\bf m}\in {\bf M}_{\J_{m+n}\times \ell}$ is identified with ${\bf m}^{(1)}\otimes \cdots\otimes {\bf m}^{(\ell)}\in ({\bf T}_{m+n}^{\rm sp})^{\otimes \ell}$. It is known that ${\bf M}_{\J_{m+n}\times \ell}$ is a $(\gl_{m+n},\gl_\ell)$-bicrystal and the map \eqref{RSK} is an isomorphism of bicrystals. Note that $\td{\mathsf{e}}_{i}$, $\td{\mathsf{f}}_{i}$ on ${\bf T}_{m+n}(\lambda,\ell)$ coincide with those on ${\bf M}_{\J_{m+n}\times \ell}$ for $i\in I_{m+n}\setminus \{\ov{m}\}$ since ${\bf T}_{m+n}(\lambda,\ell)\subset ({\bf T}_{m+n}^{\rm sp})^{\otimes \ell}$.

\begin{lem}\label{lem A1}
${\bf T}_{m+n}(\lambda,\ell)\cup\{{\bf 0}\}$ is invariant under $\td{\mathsf{e}}_{i}$ and $\td{\mathsf{f}}_{i}$ for $i\in I_{m+n}\setminus\{\ov{m}\}$, and hence ${\bf T}_{m+n}(\lambda,\ell)$ is a $\gl_{m+n}$-crystal.
\end{lem}
\pf
Let ${\bf T}\in {\bf T}_{m+n}(\lambda,\ell)$ be given and let ${\bf m}\in {\bf M}_{\J_{m+n}\times \ell}$ be the corresponding matrix.  If $\td{\mathsf{x}}_{i}{\bf m}\neq {\bf 0}$ for some $i\in I_{m+n}\setminus\{\ov{m}\}$ and $\mathsf{x}=\mathsf{e}$ or $\mathsf{f}$, then we have $Q({\bf m})=Q(\,\td{\mathsf{x}}_{i}{\bf m}\,)$ since \eqref{RSK} is an isomorphism of bicrystals, and hence $\td{\mathsf{x}}_{i}{\bf T} \in {\bf T}_{m+n}(\lambda,\ell)$ by \eqref{Pieri}. \qed\vskip 2mm

It remains to show that ${\bf T}_{m+n}(\lambda,\ell)\cup\{{\bf 0}\}$ is invariant under $\td{\mathsf{e}}_{\ov{m}}$, $\td{\mathsf{f}}_{\ov{m}}$. For this, we will show that $\td{\mathsf{x}}_{i}(T_2,T_1)$ ($\mathsf{x}=\mathsf{e},\mathsf{f}$) is also admissible, whenever it is not ${\bf 0}$, for any admissible pair $(T_2,T_1)$. We will prove the case when $\mathsf{x}=\mathsf{e}$ since the proof for $\mathsf{x}=\mathsf{f}$ is similar.

First, we need the following two lemmas, which can be checked in a straightforward manner using Algorithms  \ref{algorithm-1} and \ref{algorithm-2} in Section \ref{Definition of T_A(a)}.
\begin{lem}\label{lem:eT-1} Let $T\in {\bf T}_{m+n}(a)$ be given such that $T':=\td{\mathsf{e}}_{\ov{m}}T=(\td{\mathsf{e}}_{\ov{m}}T^{\tt R})\otimes T^{\tt L}\neq {\bf 0}$. 

Suppose that ${\mf r}_T={\mf r}_{T'}$. Then 
\begin{itemize}
\item[(1)] ${}^{\tt L}T'={}^{\tt L}T$, 

\item[(2)] ${}^{\tt R}T$ has a domino
$
\resizebox{.05\hsize}{!}{$
{\def\lr#1{\multicolumn{1}{|@{\hspace{.75ex}}c@{\hspace{.75ex}}|}{\raisebox{-.04ex}{$#1$}}}\raisebox{0.5ex}
{$\begin{array}{cc} 
\cline{1-1}
\lr{\ov{m}} \\ 
\cline{1-1}
\lr{\ov{m-1}} \\ 
\cline{1-1}
\end{array}$}}$} 
$\,, and ${}^{\tt R}T'$ is obtained from ${}^{\tt R}T$ by removing it,

\item[(3)] $T^{'\tt L^\ast}=T^{\tt L^\ast}$, when ${\mf r}_T=1$, 

\item[(4)]  ${T}^{\tt R^\ast}$ has
a domino
$
\resizebox{.05\hsize}{!}{$
{\def\lr#1{\multicolumn{1}{|@{\hspace{.75ex}}c@{\hspace{.75ex}}|}{\raisebox{-.04ex}{$#1$}}}\raisebox{0.5ex}
{$\begin{array}{cc} 
\cline{1-1}
\lr{\ov{m}} \\ 
\cline{1-1}
\lr{\ov{m-1}} \\ 
\cline{1-1}
\end{array}$}}$} 
$\,,
and $T^{'\tt R^\ast}$  is obtained from ${T}^{\tt R^\ast}$ by removing it, when ${\mf r}_T=1$.
\end{itemize}

Suppose that ${\mf r}_T\neq {\mf r}_{T'}$. Then 
\begin{itemize}
\item[(5)] $({\mf r}_T,{\mf r}_{T'})=(1,0)$ with ${\rm ht}(T^{\tt L})-a={\rm ht}(T^{\tt R})-2$,

\item[(6)] $T^{\tt L}$ and ${}^{\tt L}T$ have exactly one of $\ov{m}$ and $\ov{m-1}$,

\item[(7)]  ${}^{\tt L}T'$ is obtained from ${}^{\tt L}T$ by removing its top entry, 

\item[(8)]  ${}^{\tt R}T$ has 
a domino
$
\resizebox{.05\hsize}{!}{$
{\def\lr#1{\multicolumn{1}{|@{\hspace{.75ex}}c@{\hspace{.75ex}}|}{\raisebox{-.04ex}{$#1$}}}\raisebox{0.5ex}
{$\begin{array}{cc} 
\cline{1-1}
\lr{\ov{m}} \\ 
\cline{1-1}
\lr{\ov{m-1}} \\ 
\cline{1-1}
\end{array}$}}$} 
$\,, and ${}^{\tt R}T'$ is obtained from ${}^{\tt R}T$ by removing $\ov{m}$ or $\ov{m-1}$, which is different from the top entry of $T^{\tt L}$,

\item[(9)]  $T^{\tt L^\ast}$ has
a domino
$
\resizebox{.05\hsize}{!}{$
{\def\lr#1{\multicolumn{1}{|@{\hspace{.75ex}}c@{\hspace{.75ex}}|}{\raisebox{-.04ex}{$#1$}}}\raisebox{0.5ex}
{$\begin{array}{cc} 
\cline{1-1}
\lr{\ov{m}} \\ 
\cline{1-1}
\lr{\ov{m-1}} \\ 
\cline{1-1}
\end{array}$}}$} 
$\,, and $T^{'\tt L}=T^{\tt L}$ is obtained from $T^{\tt L^\ast}$ by removing $\ov{m}$ or $\ov{m-1}$, which is different from the top entry of $T^{\tt L}$,

\item[(10)] ${T}^{\tt R^\ast}$ has exactly one of $\ov{m}$ and $\ov{m-1}$ as its entries, and $T^{'\tt R}$ is obtained from ${T}^{\tt R^\ast}$ by removing it.
\end{itemize}
\end{lem}
\qed

\begin{lem}\label{lem:eT-2}
Let $T\in {\bf T}_{m+n}(a)$ be given such that $T':=\td{\mathsf{e}}_{\ov{m}}T= T^{\tt R}\otimes (\td{\mathsf{e}}_{\ov{m}}T^{\tt L})\neq {\bf 0}$. 

Suppose that ${\mf r}_T={\mf r}_{T'}$. Then 
\begin{itemize}


\item[(1)] ${}^{\tt L}T$ has a domino
$
\resizebox{.05\hsize}{!}{$
{\def\lr#1{\multicolumn{1}{|@{\hspace{.75ex}}c@{\hspace{.75ex}}|}{\raisebox{-.04ex}{$#1$}}}\raisebox{0.5ex}
{$\begin{array}{cc} 
\cline{1-1}
\lr{\ov{m}} \\ 
\cline{1-1}
\lr{\ov{m-1}} \\ 
\cline{1-1}
\end{array}$}}$} 
$\,, and ${}^{\tt L}T'$ is obtained from ${}^{\tt L}T$ by removing it,

\item[(2)] ${}^{\tt R}T'={}^{\tt R}T$, 

\item[(3)]  ${T}^{\tt L^\ast}$ has a domino
$
\resizebox{.05\hsize}{!}{$
{\def\lr#1{\multicolumn{1}{|@{\hspace{.75ex}}c@{\hspace{.75ex}}|}{\raisebox{-.04ex}{$#1$}}}\raisebox{0.5ex}
{$\begin{array}{cc} 
\cline{1-1}
\lr{\ov{m}} \\ 
\cline{1-1}
\lr{\ov{m-1}} \\ 
\cline{1-1}
\end{array}$}}$} 
$\,,
and $T^{'\tt L^\ast}$  is obtained from ${T}^{\tt L^\ast}$ by removing it, when ${\mf r}_T=1$,

\item[(4)] $T^{'\tt R^\ast}=T^{\tt R^\ast}$, when ${\mf r}_T=1$. 
\end{itemize}

Suppose that ${\mf r}_T\neq {\mf r}_{T'}$. Then  
\begin{itemize}
\item[(5)] $({\mf r}_T,{\mf r}_{T'})=(0,1)$ with ${\rm ht}(T^{\tt L})-a={\rm ht}(T^{\tt R})$,

\item[(6)] ${T}^{\tt R}$ and  ${}^{\tt L}T$ have exactly one of $\ov{m}$ or $\ov{m-1}$,

\item[(7)]  ${}^{\tt L}T'$ is obtained from ${}^{\tt L}T$ by removing its top entry, 

\item[(8)]  ${}^{\tt R}T$ has a domino
$
\resizebox{.05\hsize}{!}{$
{\def\lr#1{\multicolumn{1}{|@{\hspace{.75ex}}c@{\hspace{.75ex}}|}{\raisebox{-.04ex}{$#1$}}}\raisebox{0.5ex}
{$\begin{array}{cc} 
\cline{1-1}
\lr{\ov{m}} \\ 
\cline{1-1}
\lr{\ov{m-1}} \\ 
\cline{1-1}
\end{array}$}}$} 
$\,,
and ${}^{\tt R}T'$ is obtained from ${}^{\tt R}T$ by removing either $\ov{m}$ or $\ov{m-1}$, which is different from the top entry of $T^{\tt R}$,

\item[(9)] $T^{'\tt L^\ast}$ is obtained from $T^{'\tt L}$ by adding the top entry of $T^{\tt R}$,

\item[(10)] $T^{'\tt R^\ast}$ is obtained from ${T}^{\tt R}$ by removing its top entry.
\end{itemize}
\end{lem}
\qed

\vskip 2mm

Let $T_2 \in {\bf T}_{m+n}(a_2)$ and $T_1 \in {\bf T}_{m+n}(a_1)$ be given  with $a_2\geq a_1$ and $T_2\prec T_1$. Suppose that $\td{\mathsf{e}}_{\ov{m}}(T_2,T_1)\neq {\bf 0}$.
For convenience, we put
\begin{equation*}
\begin{split}
& a=a_2-a_1,\\
&(T'_2,T'_1)=\td{\mathsf{e}}_{\ov{m}}(T_2,T_1),\\
&r_i={\mf r}_{T_i}, \ \ r'_i={\mf r}_{T'_i}, \\
&2x_i={\rm ht}(T^{\tt L}_i)-a_i, \ \ 2y_i={\rm ht}(T^{\tt R}_i), \\
&2x'_i={\rm ht}(T^{' \tt L}_i)-a_i, \ \ 2y'_i={\rm ht}(T^{' \tt R}_i),\\
\end{split}
\end{equation*}
for $i=1,2$. Note that the condition (i) in  Definition \ref{admissible}(1) is equivalent to $2y_2\leq 2x_1+2r_1r_2$. 

\begin{lem}\label{lem T2R}
Suppose that $(T'_2,T'_1)=(\td{\mathsf{e}}_{\ov{m}}T_2,T_1)$ with $\td{\mathsf{e}}_{\ov{m}}T_2=(\td{\mathsf{e}}_{\ov{m}}T^{\tt R}_2)\otimes T^{\tt L}_2$. Then $T'_2\prec T'_1$.
\end{lem}
\pf  We have either $r_2=r'_2$ or $(r_2,r'_2)=(1,0)$ by Lemma \ref{lem:eT-1}(5), and $T_1=T'_1$.

(1)  It is clear that $2y'_2= 2y_2-2\leq 2x_1 =  2x'_1$ since $x'_i=x_i$ ($i=1,2$), $y'_1=y_1$, and$y'_2=y_2-1$.

(2)  If $r_1=1$ and $(r_2,r'_2)=(1,1)$, then by Lemma \ref{lem:eT-1}(4), we have  ${T}^{'\tt R^\ast}_2(i)={T}^{\tt R^\ast}_2(i)\leq {}^{\tt L}T_1(i)={}^{\tt L}T'_1(i)$ for $1\leq i\leq 2y'_2-1$. If $r_1=1$ and $(r_2,r'_2)=(1,0)$, then by Lemma \ref{lem:eT-1}(10), we have ${T}^{'\tt R}_2(i)={T}^{\tt R^\ast}_2(i) \leq {}^{\tt L}T_1(i)={}^{\tt L}T'_1(i)$ for $1\leq i\leq 2y'_2$. If $r_1=0$, then it is clear that ${T}^{'\tt R}_2(i)={T}^{\tt R}_2(i)\leq {}^{\tt L}T_1(i)={}^{\tt L}T'_1(i)$ for $1\leq i\leq 2y'_2$.

(3) Suppose that $r_2=r'_2$.  If $r_1=1$ and $(r_2,r'_2)=(1,1)$, then by Lemma \ref{lem:eT-1} (2) we have  ${}^{\tt R}T'_2(a+i)={}^{\tt R}T_2(a+i)\leq {T}_1^{\tt L^\ast}(i)={T}^{'\tt L^\ast}_1(i)$ for $1\leq i\leq 2y_2'+a_1-1$. Otherwise, we also have by Lemma \ref{lem:eT-1}(2)  ${}^{\tt R}T'_2(a+i)={}^{\tt R}T_2(a+i)\leq {T}_1^{\tt L}(i)={T}_1^{'\tt L}(i)$ for $1\leq i\leq 2y'_2+a_1-r'_2$.
  
Suppose that $(r_2,r'_2)=(1,0)$. If $r_1=0$, then by Lemma \ref{lem:eT-1}(8) we have ${}^{\tt R}T'_2(a+i)\leq {}^{\tt R}T_2(a+i)\leq T^{\tt L}_1(i)={T}_1^{'\tt L}(i)$ for $1\leq i\leq 2y'_2+a_1$.  So we assume that $r_1=1$.
Let $u_i=T_2^{\tt R}(i)$ for $1\leq i\leq 2y_2$, and let $u'_i={}^{\tt R}T_2(i)$ for $1\leq i\leq 2y_2+a_2-1=:N$. Note that 
\begin{equation}\label{inequality-1}
u'_{r}\leq u_{r-a_2+1}
\end{equation}
for $a_2\leq r\leq N-2$ by definition of ${}^{\tt R}T_2$, where $u'_{N-1}=\ov{m-1}$ and $u'_{N}=\ov{m}$ by Lemma \ref{lem:eT-1}(8).

Let $v_i=T_1^{\tt L}(i)$ for $1\leq i\leq 2x_1+a_1$ and $v^\ast_i=T_1^{\tt L^\ast}(i)$ for $1\leq i\leq 2x_1+a_1+1$.  Then there exists $p$ such that 
\begin{itemize}
\item[$\cdot$] $v^\ast_i=v_i$ for $1\leq i\leq p$,

\item[$\cdot$] $v^\ast_{p+1}=w$ for some entry $w$ in $T_1^{\tt R}$,

\item[$\cdot$] $v^\ast_i=v_{i-1}$ for $p+2 \leq i\leq 2x_1+a_1+1$,
\end{itemize}
where $p+1\geq a_1$ since $\sigma(T_1^{\tt L},T_1^{\tt R})=(a_1-1,2y_1-2x_1-1)$. Let $v'_i={}^{\tt L}T_1(i)$ for $1\leq i\leq 2x_1+1$. Then we have
\begin{itemize}
\item[$\cdot$] $v'_1=v_{i_1},\ldots, v'_{p-a_1+1}=v_{i_{p-a_1+1}}$ for some $1\leq i_1,\ldots, i_{p-a_1+1}\leq p$, 


\item[$\cdot$] $v'_{p-a_1+k}=v_{p+k-1}$ for $k\geq 2$.

\end{itemize}
Since $T_2\prec T_1$, we have by Definition \ref{admissible}(1)(ii)  
\begin{equation}\label{inequality-2}
u_{r-a_1+1}\leq v'_{r-a_1+1}=v_r
\end{equation}
for $p+1\leq r\leq 2y_2+a_1-3$, and $T_2^{\tt R^\ast}(2y_2-1)\leq v_{2y_2+a_1-2}$, where $T_2^{\tt R^\ast}(2y_2-1)=T_2^{\tt L}(2x_2)$. Combining \eqref{inequality-1} and \eqref{inequality-2}, we get 
\begin{equation}\label{inequality-3}
u'_r\leq v_{r-a}
\end{equation}
for $p+a+1\leq r\leq N-2$.  Since ${}^{\tt R}T'_2(i)={}^{\tt R}T_2(i)$ for $1\leq i\leq N$ by Lemma \ref{lem:eT-1}(8), we have by \eqref{inequality-3} 
\begin{equation}\label{inequality-3-1}
{}^{\tt R}T'_2(i+a)=u'_{i+a}\leq v_i=T_1^{\tt L}(i)
\end{equation}
for $p+1 \leq i\leq N-a-2=2y'_2+a_1-1$.   Note that  $T_1^{\tt L^\ast}(i)=T_1^{\tt L}(i)$ for $1\leq i\leq p$ and  hence
\begin{equation}\label{inequality-3-2}
{}^{\tt R}T'_2(i+a)={}^{\tt R}T_2(i+a) \leq T_1^{\tt L}(i)
\end{equation}
for $1\leq i\leq  p$. Also, we have ${}^{\tt R}T'_2(N-1)=T_1^{\tt L}(2x_1)=T_2^{\tt R^\ast}(2y_2-1)\leq v_{2y_2+a_1-2}=T_1^{\tt L}(2y_2+a_1-2)=T_1^{\tt L}(N-a-1)$. Hence by \eqref{inequality-3-1} and \eqref{inequality-3-2} we conclude that ${}^{\tt R}T'_2(i+a)\leq T_1^{'\tt L}(i)=T_1^{\tt L}(i)$ for $1\leq i\leq 2y'_2+a_1$.
 
Therefore, we have $T'_2=(\mathsf{e}_{\ov{m}}T_2)\prec T_1=T'_1$ by (1), (2) and (3).\qed
\vskip 2mm

\begin{lem}\label{lem T2L}
Suppose that $(T'_2,T'_1)=(\td{\mathsf{e}}_{\ov{m}}T_2,T_1)$ with $\td{\mathsf{e}}_{\ov{m}}T_2=T^{\tt R}_2\otimes(\td{\mathsf{e}}_{\ov{m}}T^{\tt L}_2)$. Then $T'_2\prec T'_1$.
\end{lem}
\pf We have either $r_2=r'_2$ or $(r_2,r_2')=(0,1)$ by Lemma \ref{lem:eT-2}(5), and $T_1=T'_1$.

(1)  It is clear that $2y'_2\leq 2x'_1 +2r'_1r'_2$ since  $y'_i=y_i$ ($i=1,2$), $x'_1=x_1$, $x'_2=x_2-1$, and $r'_1=r_1$.

(2) If $r_1=1$ and $(r_2,r'_2)=(1,1)$, then ${T'_2}^{\tt R^\ast}={T}^{\tt R^\ast}_2$ by Lemma \ref{lem:eT-2}(4) and hence ${T}^{'\tt R^\ast}_2(i)={T}^{\tt R^\ast}_2(i)\leq {}^{\tt L}T_1(i)={}^{\tt L}T'_1(i)$ for $1\leq i\leq 2y'_2-1$. If $r_1=1$ and $(r_2,r'_2)=(0,1)$, then by Lemma \ref{lem:eT-2}(10), ${T}^{'\tt R^\ast}_2(i)={T}^{\tt R}_2(i)\leq {}^{\tt L}T_1(i)={}^{\tt L}T'_1(i)$ for $1\leq i\leq 2y'_2-1$. If $r_1=0$, then it is clear that ${T'_2}^{\tt R}(i)={T}^{\tt R}_2(i)\leq {}^{\tt L}T_1(i)={}^{\tt L}T'_1(i)$ for $1\leq i\leq 2y'_2$.

(3) Suppose that $r_2=r'_2$. Then we have ${}^{\tt R}T_2={}^{\tt R}T'_2$ by Lemma \ref{lem:eT-2}(2) and hence ${}^{\tt R}T'_2(a+i)= {}^{\tt R}T_2(a+i)\leq {T}^{\tt L}_1(i)={T}^{'\tt L}_1(i)$ or ${}^{\tt R}T'_2(a+i)\leq{T}^{\tt L^\ast}_1(i)={T}^{'\tt L^\ast}_1(i)$ for $1\leq i\leq 2y'_2+a_1-r'_2$. 

Suppose that $(r_2,r'_2)=(0,1)$. By Lemma \ref{lem:eT-2}(8), we have ${}^{\tt R}T'_2(a+i)\leq {}^{\tt R}T_2(a+i)\leq {T}^{\tt L}_1(i)={T}^{'\tt L}_1(i)$ or ${}^{\tt R}T'_2(a+i)\leq {T}^{\tt L}_1(i)\leq {T}^{'\tt L^\ast}_1(i)$ for $1\leq i\leq 2y'_2+a_1$.

Therefore, we have $T'_2=(\td{\mathsf{e}}_{\ov{m}}T_2)\prec T_1=T'_1$ by (1), (2) and (3).
\qed

\begin{lem}\label{lem T1R}
Suppose that $(T'_2,T'_1)=(T_2,\td{\mathsf{e}}_{\ov{m}}T_1)$ with $\td{\mathsf{e}}_{\ov{m}}T_1=(\td{\mathsf{e}}_{\ov{m}}T^{\tt R}_1)\otimes T^{\tt L}_1$. Then $T'_2\prec T'_1$.
\end{lem}
\pf  We have either $r_1=r'_1$ or $(r_1,r'_1)=(1,0)$  by Lemma \ref{lem:eT-1}(5), and $T_2=T'_2$.

(1) Note that $x'_i=x_i$ ($i=1,2$), $y'_1=y_1-1$, $y'_2=y_2$, and $r'_2=r_2$. If  $r_1=r'_1$, then it is clear that $2y'_2\leq 2x'_1+2r'_1r'_2$.

Suppose that  $(r_1,r'_1)=(1,0)$. If $y_2\leq x_1$, then we have $2y'_2\leq 2x'_1=2x'_1+2r'_1r'_2$.
Now, we claim that we have a contradiction when $y_2>x_1$ (or $y_2=x_1+1$), that is, $2y_2=2x_2+2r_1r_2$ with $r_1=r_2=1$. Since $(r_1,r'_1)=(1,0)$, we have $y_1=x_1+1$. 
By Lemma \ref{lem:eT-1}(6) and Definition \ref{admissible}(1)(ii),  the top entry of $T_2^{\tt R^\ast}$ is no greater than $\ov{m-1}$. On the other hand, by  Lemma \ref{lem:eT-1}(9) and Definition \ref{admissible}(1)(iii), ${}^{\tt R}T_2$ has a domino 
$
\resizebox{.05\hsize}{!}{$
{\def\lr#1{\multicolumn{1}{|@{\hspace{.75ex}}c@{\hspace{.75ex}}|}{\raisebox{-.04ex}{$#1$}}}\raisebox{0.5ex}
{$\begin{array}{cc} 
\cline{1-1}
\lr{\ov{m}} \\ 
\cline{1-1}
\lr{\ov{m-1}} \\ 
\cline{1-1}
\end{array}$}}$} 
$\,, which also implies that the top entry of $T_2^{\tt R}$ is $\ov{m}$.
If $x_2+1<y_2$, then $T_2^{\tt R}$ has a domino
$
\resizebox{.05\hsize}{!}{$
{\def\lr#1{\multicolumn{1}{|@{\hspace{.75ex}}c@{\hspace{.75ex}}|}{\raisebox{-.04ex}{$#1$}}}\raisebox{0.5ex}
{$\begin{array}{cc} 
\cline{1-1}
\lr{\ov{m}} \\ 
\cline{1-1}
\lr{\ov{m-1}} \\ 
\cline{1-1}
\end{array}$}}$} 
$\,, and  $\td{\mathsf{e}}_{\ov{m}}(T_2,T_1)=(\td{\mathsf{e}}_{\ov{m}}T_2,T_1)$, which is a contradiction. 
Next, assume that $x_2+1=y_2$. Put $y=2y_2$. Consider $T_2^{\tt R^\ast}(y-1)$, the top entry of $T_2^{\tt R^\ast}$.  If $T_2^{\tt R^\ast}(y-1)=T_2^{\tt R}(y)=\ov{m}$, then $T_2^{\tt L}(y-2)= T_2^{\tt R}(y)= \ov{m}$. But this implies that the first two top entries of $T_2^{\tt R}$ are equal to those of ${}^{\tt R}T_2$, which are $\ov{m}$ and $\ov{m-1}$.  So we have a contradiction $\td{\mathsf{e}}_{\ov{m}}(T_2,T_1)=(\td{\mathsf{e}}_{\ov{m}}T_2,T_1)$.  If $T_2^{\tt R^\ast}(y-1)=T_2^{\tt R}(y-1)$, then $T_2^{\tt R}(y)=\ov{m}$ and $ T_2^{\tt R}(y-1)=\ov{m-1}$ since $T_2^{\tt R^\ast}(y-1)\leq \ov{m-1}$, which also yields a contradiction $\td{\mathsf{e}}_{\ov{m}}(T_2,T_1)=(\td{\mathsf{e}}_{\ov{m}}T_2,T_1)$. This proves our claim.

(2) Suppose that $r_1=r'_1$. Then by Lemma \ref{lem:eT-1}(1) we have ${}^{\tt L}T_1={}^{\tt L}T'_1$ and hence ${T}_2^{'\tt R}(i)={T}_2^{\tt R}(i)\leq {}^{\tt L}T_1(i)={}^{\tt L}T'_1(i)$ or ${T}^{'\tt R^\ast}_2(i)={T}^{\tt R^\ast}_2(i)\leq {}^{\tt L}T_1(i)={}^{\tt L}T'_1(i)$ for $1\leq i\leq 2y'_2-r'_2$.  

Suppose that $(r_1,r'_1)=(1,0)$, where we have $y_2\leq x_1$ by (1). If $r_2=0$, then by Lemma \ref{lem:eT-1}(7) we have $T_2^{'\tt R}(i)={T}^{\tt R}_2(i)\leq {}^{\tt L}T_1(i)={}^{\tt L}T'_1(i)$ for $1\leq i\leq 2y'_2$. So we assume that $r_2=1$.

Let $u_i=T^{\tt R}_2(i)$ for $1\leq i\leq 2y_2$ and  $u^\ast_i=T^{\tt R^\ast}_2(i)$ for $1\leq i\leq 2y_2-1$. There exists $p\geq 1$ such that
\begin{itemize}
\item[$\cdot$] $u_i^\ast=u_i$ for $1\leq i\leq p$,

\item[$\cdot$] $u_i^\ast=u_{i+1}$ for $p+1\leq i\leq 2y_2-1$
\end{itemize}
Let $u'_i={}^{\tt R}T_2(i)$ for $1\leq i\leq 2y_2+a_2-1$. Then we have 
\begin{equation}\label{inequality-4}
u'_{p+a_2+i-1}=u_{p+i}
\end{equation}
for $1\leq i \leq 2y_2-p$. Let $v_i=T_1^{\tt L}(i)$ for $1\leq i\leq 2x_1+a_1$ and $v^\ast_i=T_1^{\tt L^\ast}(i)$ for $1\leq i\leq 2x_1+a_1-1$.
Then we see that for $1\leq i\leq 2x_1+a_1-1$,
\begin{equation}\label{inequality-5}
v^\ast_i=v_i
\end{equation}
while $v^\ast_{a_1+2x_1}=\ov{m-1}$, $v^\ast_{a_1+2x_1+1}=\ov{m}$, and
$v_{a_1+2x_1}$ is either $\ov{m}$ or $\ov{m-1}$ by Lemma \ref{lem:eT-1}(9).
Since $T_2\prec T_1$, we have by Definition \ref{admissible}(1)(iii) $u'_{a+i}\leq v^\ast_i$ for $1\leq i\leq 2y_2+a_1-1$.
Since $y_2\leq x_1$, we have by \eqref{inequality-5}
\begin{equation}\label{inequality-6}
u'_{a+i}\leq v_i
\end{equation}
for $1\leq i\leq 2y_2+a_1-1$.  On the other hand, let $v'_i={}^{\tt L}T_1(i)$ for $1\leq i\leq 2x_1+1$. Since $r_1=1$, we have $v_{a_1+i-1}\leq v'_i$ for $1\leq i \leq 2x_1-1$.

Now consider $T^{'\tt R}_2=T^{\tt R}_2$ and ${}^{\tt L}T'_1$. By Lemma \ref{lem:eT-1}(7) we have 
${}^{\tt L}T_1(i)={}^{\tt L}T'_1(i)$ 
for $1\leq i\leq 2x_1$. Since $u_i^\ast=u_i$ for $1\leq i\leq p$, we have 
\begin{equation}\label{inequality-6-1}
T_2^{'\tt R}(i)\leq {}^{\tt L}T'_1(i)
\end{equation}
for $1\leq i\leq p$. By \eqref{inequality-4} and \eqref{inequality-6}, we have 
\begin{equation}\label{inequality-6-2}
T_2^{'\tt R}(p+i)=u_{p+i}\leq v_{a_1+p+i-1}\leq v'_{p+i}={}^{\tt L}T'_1(p+i)
\end{equation}
for $1\leq i\leq 2y_2-p$, which implies that $T_2^{'\tt R}(i)\leq {}^{\tt L}T'_1(i)$ for $p+1\leq i\leq 2y_2$. Therefore, $T_2^{'\tt R}(i)\leq {}^{\tt L}T'_1(i)$ for $1\leq i\leq 2y_2$ by \eqref{inequality-6-1} and \eqref{inequality-6-2}.

(3) If $(r_1,r'_1)=(1,1)$ and $r_2=1$, then by Lemma \ref{lem:eT-1}(3) we have $T_1^{\tt L^\ast}=T_1^{'\tt L^\ast}$ and hence ${}^{\tt R}T'_2(a+i) \leq T_1^{\tt L^\ast}(i)=T_1^{'\tt L^\ast}(i)$ for $1\leq i\leq 2y'_2+a_1-1$. If $(r_1,r'_1)=(1,0)$ and $r_2=1$, then by  Lemma \ref{lem:eT-1}(9) we have  
${}^{\tt R}T'_2(a+i) \leq T_1^{\tt L^\ast}(i)={T'_1}^{\tt L}(i)$ for $1\leq i\leq 2y'_2+a_1-1$. If $r_2=0$, then it is clear that ${}^{\tt R}T'_2(a+i) \leq T_1^{\tt L}(i)={T'_1}^{\tt L}(i)$ for $1\leq i\leq 2y'_2+a_1$.

Therefore, we have $T'_2=T_2\prec (\td{\mathsf{e}}_{\ov{m}}T_1)=T'_1$ by (1), (2) and (3).\qed\vskip 2mm

\begin{lem}\label{lem T1L}
Suppose that $(T'_2,T'_1)=(T_2,\td{\mathsf{e}}_{\ov{m}}T_1)$ with $\td{\mathsf{e}}_{\ov{m}}T_1=T^{\tt R}_1\otimes (\td{\mathsf{e}}_{\ov{m}}T^{\tt L}_1)$. Then $T'_2\prec T'_1$.
\end{lem}
\pf   We have either $r_1=r'_1$ or $(r_1,r'_1)=(0,1)$ by Lemma \ref{lem:eT-2}(5), and $T_2=T'_2$.

(1) Note that  $y'_i=y_i$ ($i=1,2$), $x'_1=x_1-1$, $x'_2=x_2$, and $r'_2=r_2$. If $y_2\leq x_1-1=x'_1$, then we have $2y'_2\leq 2x'_1+2r'_1r'_2$. So we assume that $y_2\geq x_1$, that is, $y_2= x_1$ or $y_2=x_1+1$.

(i) Suppose that $y_2=x_1$. If $r_1=r'_1=0$,  then  by Lemma \ref{lem:eT-2}(1) and Definition \ref{admissible}(1)(ii), $T_2^{\tt R}$ has a domino
$
\resizebox{.05\hsize}{!}{$
{\def\lr#1{\multicolumn{1}{|@{\hspace{.75ex}}c@{\hspace{.75ex}}|}{\raisebox{-.04ex}{$#1$}}}\raisebox{0.5ex}
{$\begin{array}{cc} 
\cline{1-1}
\lr{\ov{m}} \\ 
\cline{1-1}
\lr{\ov{m-1}} \\ 
\cline{1-1}
\end{array}$}}$} 
$\,, which implies that $\td{\mathsf{e}}_{\ov{m}}(T_2,T_1)=(\td{\mathsf{e}}_{\ov{m}}T_2,T_1)$, a contradiction. So we have $(r_1,r'_1)=(0,1)$ or $(1,1)$.  

Now, suppose that $r_2=0$. If $x_2<y_2$, then the first two top entries of ${T}_2^{\tt R}$ and ${}^{\tt R}T_2$ are the same, and they are $\ov{m}$ and $\ov{m-1}$ by Definition \ref{admissible}(1)(iii). But this implies that $\td{\mathsf{e}}_{\ov{m}}(T_2,T_1)=(\td{\mathsf{e}}_{\ov{m}}T_2,T_1)$, which is a contradiction. So we have $x_2=y_2$. Now consider the first two top entries of $T_2^{\tt L}$ and $T_2^{\tt R}$. 
Put $x=2x_2$, and $w_1=T_2^{\tt L}(x)$, $w_2=T_2^{\tt L}(x-1)$, $w_3=T_2^{\tt R}(x)$, $w_4=T_2^{\tt R}(x-1)$.
First, we have ${}^{\tt R}T_2(x)=\ov{m}$, ${}^{\tt R}T_2(x-1)=\ov{m-1}$ by Definition \ref{admissible}(1)(iii), which implies that  $w_1= \ov{m}$. Second, we have $w_3\leq \ov{m-1}$ since $w_3\leq {}^{\tt L}T_1(x)\leq \ov{m-1}$ when  $r_1=0$, and $w_3\leq{}^{\tt L}T_1(x-1)= \ov{m-1}$ when $r_1=1$. This implies that $w_2=\ov{m-1}$, and hence $\td{\mathsf{e}}_{\ov{m}}(T_2,T_1)=(\td{\mathsf{e}}_{\ov{m}}T_2,T_1)$, which is also a contradiction. So we should have $r_2=r'_2=1$. Hence, it follows that  $2y'_2=2y_2=2x_1=2(x_1-1)+2= 2x'_1+2r'_1r'_2$.

(ii) Suppose that $y_2=x_1+1$ with $r_1=r_2=1$. Since $r_1=r'_1=1$, ${}^{\tt L}T_1$ has a domino
$
\resizebox{.05\hsize}{!}{$
{\def\lr#1{\multicolumn{1}{|@{\hspace{.75ex}}c@{\hspace{.75ex}}|}{\raisebox{-.04ex}{$#1$}}}\raisebox{0.5ex}
{$\begin{array}{cc} 
\cline{1-1}
\lr{\ov{m}} \\ 
\cline{1-1}
\lr{\ov{m-1}} \\ 
\cline{1-1}
\end{array}$}}$} 
$\, by Lemma \ref{lem:eT-2}(1). But then $T_2^{\tt R^\ast}$ and hence $T_2^{\tt R}$ has a domino
$
\resizebox{.05\hsize}{!}{$
{\def\lr#1{\multicolumn{1}{|@{\hspace{.75ex}}c@{\hspace{.75ex}}|}{\raisebox{-.04ex}{$#1$}}}\raisebox{0.5ex}
{$\begin{array}{cc} 
\cline{1-1}
\lr{\ov{m}} \\ 
\cline{1-1}
\lr{\ov{m-1}} \\ 
\cline{1-1}
\end{array}$}}$} 
$\,, which gives a contradiction $\td{\mathsf{e}}_{\ov{m}}(T_2,T_1)=(\td{\mathsf{e}}_{\ov{m}}T_2,T_1)$. 
 
Therefore, it follows from (i) and (ii) that $2y'_2\leq 2x'_1+2r'_1r'_2$, when $y_2\geq x_1$.

(2) Suppose that $r_1=r'_1$.  Note that $y_2\leq x_1$ by (1) (ii).
If $y_2<x_1$, then by Lemma \ref{lem:eT-2}(1), we have ${T}_2^{'\tt R}(i)=T_2^{\tt R}(i)\leq {}^{\tt L}T_1(i)={}^{\tt L}T'_1(i)$  for $1\leq i\leq 2y'_2$ or ${T}^{'\tt R^\ast}_2(i)=T_2^{\tt R^\ast}(i)\leq {}^{\tt L}T_1(i)={}^{\tt L}T'_1(i)$ for $1\leq i\leq 2y'_2-1$. If $y_2=x_1$, then $r_1=r_2=1$ by (1) (i), and we also have  ${T}^{'\tt R^\ast}_2(i)={T}^{\tt R^\ast}_2(i)\leq {}^{\tt L}T'_1(i)={}^{\tt L}T_1(i)$ for $1\leq i\leq 2y'_2-1$ by Lemma \ref{lem:eT-2}(1). 

Suppose that  $(r_1,r'_1)=(0,1)$.  
Then by Lemma \ref{lem:eT-2}(7), we have ${T}^{'\tt R}_2(i)={T}^{\tt R}_2(i)\leq {}^{\tt L}T_1(i)={}^{\tt L}T'_1(i)$  for $1\leq i\leq 2y'_2$ when $r_2=0$, and ${T}^{'\tt R^\ast}_2(i)\leq{T}^{\tt R}_2(i)\leq {}^{\tt L}T_1(i)={}^{\tt L}T'_1(i)$ for $1\leq i\leq 2y'_2-1$  when $r_2=1$.

(3) Suppose that $r_1=r'_1$. If $(r_1,r'_1)=(1,1)$ and $r_2=1$, then we have by Lemma \ref{lem:eT-2}(3)  ${}^{\tt R}T'_2(a+i)={}^{\tt R}T_2(a+i) \leq {T}^{'\tt L^\ast}_1(i)={T}^{\tt L^\ast}_1(i)$  for $1\leq i\leq 2y'_2+a_1-1$. Otherwise it is clear that ${}^{\tt R}T'_2(a+i)={}^{\tt R}T_2(a+i) \leq {T}^{'\tt L}_1(i)={T}^{\tt L}_1(i)$   for $1\leq i\leq 2y'_2+a_1-r'_2$. 

Suppose that  $(r_1,r'_1)=(0,1)$.  If $r_2=0$, then it is clear that ${}^{\tt R}T'_2(a+i)={}^{\tt R}T_2(a+i) \leq {T}^{'\tt L}_1(i)={T}^{\tt L}_1(i)$  for $1\leq i\leq 2y'_2+a_1$. If $r_2=1$ with $y_2<x_1$, then we have by Lemma \ref{lem:eT-2}(9) ${}^{\tt R}T'_2(a+i)={}^{\tt R}T_2(a+i) \leq {T}^{'\tt L^\ast}_1(i)={T}^{\tt L}_1(i)$  for $1\leq i\leq 2y'_2+a_1-1$. Suppose that $r_2=1$ and $y_2=x_1$. Put $x=2x_1$. Since $T_2^{\tt R}(x)\leq {}^{\tt L}T_1(x)$ and ${}^{\tt L}T_1(x)=T_1^{\tt R}(x)$, we have $T_2^{\tt R}(x)\leq T_1^{\tt R}(x)$. Using this fact and Lemma \ref{lem:eT-2}(9), we can check that ${}^{\tt R}T_2(a+x-1)=T_2^{\tt R}(x)\leq T_1^{\tt R}(x)=T_1^{'\tt L^\ast}(x-1)$ and hence ${}^{\tt R}T_2(a+i)\leq T_1^{'\tt L^\ast}(i)$  for $1\leq i\leq 2y'_2+a_1-1$.

Therefore, we have $T'_2=T_2\prec (\td{\mathsf{e}}_{\ov{m}}T_1)=T'_1$ by (1), (2) and (3).\qed\vskip 3mm

\begin{lem}\label{lem nonspin spin}
Suppose that $T_2 \in {\bf T}_{m+n}(a_2)$ and  $T_1 \in {\bf T}^{\rm sp}_{m+n}$ with $a_2\geq a_1:={\mf r}_{T_1}$ and $T_2\prec T_1$. 
 If $\td{\mathsf{e}}_{\ov{m}}(T_2,T_1)=(T_2',T_1')\neq {\bf 0}$, then $T_2' \prec T_1'$.
\end{lem}
\pf   Put $\ov{T}_1=(T_1,H_{(1^N)})$ for a sufficiently large even integer $N$. Then $\ov{T}_1\in {\bf T}_{m+n}(a_1)$, where  
\begin{itemize}
\item[$\cdot$] ${\mf r}_{\ov{T}_1}={\mf r}_{T_1}$,

\item[$\cdot$] $\ov{T}_1^{\tt L}=T_1$, ${}^{\tt L}\ov{T}_1=T_1$,

\item[$\cdot$] $\ov{T}_1^{\tt L^\ast}$ is obtained by adding the largest entry of $H_{(1^N)}$ at the bottom of $T_1$ when ${\mf r}_{\ov{T}_1}=1$.
\end{itemize}
It is not difficult to see that  $T_2\prec T_1$ if and only if $T_2\prec \ov{T}_1$. Now applying Lemmas \ref{lem T2R}, \ref{lem T2L}, and \ref{lem T1L} to the pair $(T_2,\ov{T}_1)$, we conclude that $T'_2\prec T'_1$.
\qed\vskip 3mm
 
For the admissible pairs $(T_2,T_1)$ in Definition \ref{admissible}(2) and (3), we can check without difficulty that if $\td{\mathsf{e}}_{i}(T_2,T_1)=(T_2',T_1')\neq {\bf 0}$ for some $i\in I_{m+n}$, then $T_2' \prec T_1'$.
  
Hence by Lemmas \ref{lem A1}--\ref{lem nonspin spin}, we conclude that ${\bf T}_{m+n}(\lambda,\ell)\cup\{{\bf 0}\}$ is invariant under $\td{\mathsf{x}}_i$ for $\mathsf{x}=\mathsf{e}, \mathsf{f}$ and $i\in I_{m+n}$, which proves Theorem \ref {crystal invariance of osp tableaux} (1). 
\begin{lem}\label{lem cnn}
${\bf T}_{m+n}(\lambda,\ell)$ is a connected ${\mf d}_{m+n}$-crystal with highest weight $\Lambda_{m+n}(\lambda,\ell)$.
\end{lem}
\pf  Let ${\bf H}_{(\lambda,\ell)}=(T_L,\ldots, T_0)\in {\bf T}_{m+n}(\lambda,\ell)$ be such that 
\begin{itemize}
\item[$\cdot$] $T_{k}$ is empty for $0\leq k\leq q_+$ when $\ell-2\lambda_1\geq 0$,

\item[$\cdot$] $T_{0}=$ \resizebox{.035\hsize}{!}{$\boxed{\ov{m}}$}\  and $T_{k}=$
$
\resizebox{.07\hsize}{!}{$
{\def\lr#1{\multicolumn{1}{|@{\hspace{.75ex}}c@{\hspace{.75ex}}|}{\raisebox{-.04ex}{$#1$}}}\raisebox{0.25ex}
{$\begin{array}{cc} 
\cline{1-2}
\lr{\ov{m}} &  \lr{\ov{m}} \\
\cline{1-2}
\end{array}$}}$} 
$\, for $1\leq k\leq q_-$ when $\ell-2\lambda_1\leq 0$,

\item[$\cdot$] $T_{q_\pm +k}=H_{(1^{a_k})}\in {\bf T}_{m+n}(a_k)$ for $1\leq k\leq M_\pm$, where $a_k$ is as in \eqref{notations for lengths}.

\end{itemize}

We claim that any ${\bf T}\in {\bf T}_{m+n}(\lambda,\ell)$ is connected to ${\bf H}_{(\lambda,\ell)}$ under $\td{\mathsf{e}}_i$ for $i\in I_{m+n}$, where ${\rm wt}({\bf H}_{(\lambda,\ell)})=\Lambda_{m+n}(\lambda,\ell)$. We use induction on $|{\bf T}|=\sum_{k=0}^L |T_k|$. Note that $|{\bf H}_{(\lambda,\ell)}|=\sum_{i\geq 1}\lambda_i< |{\bf T}|$ for all ${\bf T}\in {\bf T}_{m+n}(\lambda,\ell)\setminus\{{\bf H}_{(\lambda,\ell)}\}$.

Suppose that ${\bf T}$ is given. We may assume that $\td{\mathsf{e}}_i {\bf T}={\bf 0}$ for $i\in I_{m+n}\setminus\{\ov{m}\}$ since $|\td{\mathsf{e}}_i {\bf T}|=|{\bf T}|$ whenever $\td{\mathsf{e}}_i {\bf T}\neq {\bf 0}$ for $i\in I_{m+n}\setminus\{\ov{m}\}$.
So, it is enough to show that ${\bf T} = {\bf H}_{(\lambda,\ell)}$ or $\td{\mathsf{e}}_{\ov{m}}{\bf T}\neq {\bf 0}$, which implies $|\td{\mathsf{e}}_{\ov{m}} {\bf T}|<|{\bf T}|$. 

{\em Step 1}. Suppose that $\ell-2\lambda_1\geq 0$. By Definition \ref{admissible}(1), $(T_{q_+},\ldots, T_{0})\in SST_{\J_{m+n}}(\alpha)$, where $\alpha=((q_++r_+)^k)/\nu$ for some $k\in 2\Z_{\geq 0}$ and $\nu\in \cP$ such that each column of $\nu$ is also of even length. Since $\td{\mathsf{e}}_i {\bf T}={\bf 0}$ for $i\in I_{m+n}\setminus\{\ov{m}\}$, $(T_{q_+},\ldots, T_{0})$ is a $\gl_{m+n}$-highest weight element, and 
hence each of $T_0$, $T_i^{\tt L}$, and $T_i^{\tt R}$ ($1\leq i\leq q_+$) is a $\gl_{m+n}$-highest weight element $H_{(1^d)}$ for some $d\in 2\Z_{\geq 0}$.  If $T_i$ is not empty for some $0\leq i\leq q_+$, then $\td{\mathsf{e}}_{\ov{m}}(T_{q_+},\ldots, T_{0})\neq {\bf 0}$, and hence $\td{\mathsf{e}}_{\ov{m}} {\bf T}\neq{\bf 0}$ by tensor product rule. 

Suppose that $\ell-2\lambda_1\leq 0$.  By Definition \ref{admissible}(2), $(T_{q_-},\ldots,T_{0})\in SST_{\J_{m+n}}(\beta)$, where $\beta=((q_-+r_-)^k)/\nu$ for some $k\in 1+2\Z_{\geq 0}$ and $\nu\in \cP$ such that each column of $\nu$ is of even length. As in  (1), each of $T_0$, $T_i^{\tt L}$, and $T_i^{\tt R}$ ($1\leq i\leq q_-$) is $H_{(1^d)}$ for some $d\in 1+2\Z_{\geq 0}$.
If $(T_{q_-},\ldots,T_{0})$ has a column of height greater than 1, then we have  $\td{\mathsf{e}}_{\ov{m}}(T_{q_-},\ldots, T_{0})\neq {\bf 0}$ and hence $\td{\mathsf{e}}_{\ov{m}} {\bf T}\neq{\bf 0}$. Otherwise $T_i=$ \resizebox{.035\hsize}{!}{$\boxed{\ov{m}}$} or   $
\resizebox{.07\hsize}{!}{$
{\def\lr#1{\multicolumn{1}{|@{\hspace{.75ex}}c@{\hspace{.75ex}}|}{\raisebox{-.04ex}{$#1$}}}\raisebox{0.25ex}
{$\begin{array}{cc} 
\cline{1-2}
\lr{\ov{m}} &  \lr{\ov{m}} \\
\cline{1-2}
\end{array}$}}$} 
$\  for $0\leq i\leq q_-$.

By {\em Step 1}, we may assume from now on that $T_{i}$ is empty for $0\leq i\leq q_+$ when $\ell-2\lambda_1\geq 0$, and $T_{0}=$ \resizebox{.035\hsize}{!}{$\boxed{\ov{m}}$}\  and $T_{i}=$
$
\resizebox{.07\hsize}{!}{$
{\def\lr#1{\multicolumn{1}{|@{\hspace{.75ex}}c@{\hspace{.75ex}}|}{\raisebox{-.04ex}{$#1$}}}\raisebox{0.25ex}
{$\begin{array}{cc} 
\cline{1-2}
\lr{\ov{m}} &  \lr{\ov{m}} \\
\cline{1-2}
\end{array}$}}$} 
$\, for $1\leq i\leq q_-$ when $\ell-2\lambda_1\leq 0$. 

{\em Step 2}. Consider $T_{q_\pm +1}$.
Suppose that $\ell-2\lambda_1\geq 0$. Then $T^{\tt R}_{q_+ +1}$ is empty by Definition \ref{admissible}(1) since  $T_{q_+}\in {\bf T}^{\rm sp +}_{m+n}$ is empty. Also, we have $T^{\tt L}_{q_+ +1}=H_{(1^{a_1})}$ since $(T_{q_++1},\ldots, T_0)$ is a $\gl_{m+n}$-highest weight element with  $T_{i}$ empty for $0\leq i\leq q_+$.
 Hence, $T_{q_\pm +1}$ is a ${\mf d}_{m+n}$-highest weight element. 
 
Suppose that  $\ell-2\lambda_1< 0$. Then ${\rm ht}(T^{\tt R}_{q_- +1})\leq 2$ by Definition \ref{admissible}(1)(i). Suppose that ${\rm ht}(T^{\tt R}_{q_-+1})=2$ (with $r_{T_{q_-+1}}=1$), and let $x=T^{\tt R}_{q_- +1}(2)$ and $y=T^{\tt R}_{q_- +1}(1)$.  By Definition \ref{admissible}(1)(ii), we have $x=\ov{m}$. If $y>\ov{m-1}$, then $\td{\mathsf{e}}_i (T_{q_- +1},\ldots, T_0)\neq {\bf 0}$ for some $i\in I_{m+n}\setminus\{\ov{m}\}$, which is a contradiction. So we have $y=\ov{m-1}$. Then $\td{\mathsf{e}}_{\ov{m}}T_{q_-}\neq {\bf 0}$ and hence $\td{\mathsf{e}}_{\ov{m}}{\bf T}\neq {\bf 0}$ since $T_i=$ \resizebox{.035\hsize}{!}{$\boxed{\ov{m}}$} or   $
\resizebox{.07\hsize}{!}{$
{\def\lr#1{\multicolumn{1}{|@{\hspace{.75ex}}c@{\hspace{.75ex}}|}{\raisebox{-.04ex}{$#1$}}}\raisebox{0.2ex}
{$\begin{array}{cc} 
\cline{1-2}
\lr{\ov{m}} &  \lr{\ov{m}} \\
\cline{1-2}
\end{array}$}}$} 
$ for $0\leq i\leq q_-$. If $T^{\tt R}_{q_- +1}$ is empty, then by similar arguments as in (2), we have $T^{\tt L}_{q_- +1}=H_{(1^{a_1})}$.

{\em Step 3}. By {\em Step 2} we may assume that $T_{q_\pm+1}$ is a ${\mf d}_{m+n}$-highest weight element.  Suppose that there exists $k\geq 1$ such that $T_{q_\pm+i}$ is a ${\mf d}_{m+n}$-highest weight element  for $1\leq i\leq k$. 
 Consider $T_{q_\pm+k+1}$. By Definition \ref{admissible}(1)(i), we have $T_{q_\pm+k+1}^{\tt R}$ is empty. 
By Definition \ref{admissible}(1)(iii), the first $a_{k}$ entries of $T^{\tt L}_{q_\pm+k+1}$ from the top are $\ov{m}, \ldots, \ov{m-a_{k}+1}$. Since ${\bf T}$ is a $\gl_{m+n}$-highest weight element, we have $T^{\tt L}_{q_\pm+k+1}=H_{(1^{a_{k+1}})}$, and hence $T_{q_\pm+k+1}$ is a ${\mf d}_{m+n}$-highest weight element. Applying this argument inductively, we conclude that ${\bf T}={\bf H}_{(\lambda,\ell)}$.
\qed
\vskip 2mm

This proves Theorem \ref {crystal invariance of osp tableaux}(2), and completes the proof of the theorem.

\end{document}